\documentclass{article} 

\usepackage{booktabs}       
\usepackage{amsfonts}       
\usepackage{nicefrac}       
\usepackage{microtype}      
\usepackage{fncylab,enumerate,wrapfig,placeins}
\usepackage{bbm}
\usepackage{graphicx}
\usepackage{enumitem}  
\usepackage{amsmath,amsthm,amsfonts,amssymb}  
\usepackage{pifont}
\usepackage{geometry}
\usepackage[usenames,dvipsnames]{xcolor}
\usepackage[colorlinks = true, 
linkcolor = RoyalPurple,
urlcolor  = CadetBlue,
citecolor = blue,
anchorcolor = blue]{hyperref}

\pdfoutput=1  

\usepackage{hyperref}       
\usepackage{url}            

\usepackage{graphicx}  
\usepackage{relsize} 
\usepackage{accents}  
\usepackage{caption}
\usepackage{textgreek,upgreek,bm}

\usepackage{titlesec}

\newtheorem{theorem}{Theorem}[section]
\newtheorem{lemma}[theorem]{Lemma} 
\newtheorem{proposition}[theorem]{Proposition} 
\newtheorem{corollary}[theorem]{Corollary}

\usepackage{amsmath,amsfonts,amssymb}  
\usepackage{bbm}
\usepackage{xcolor}

\usepackage{hyperref} 
\usepackage{url}            
\usepackage{booktabs}       
\usepackage{nicefrac}       
\usepackage{microtype}      
\usepackage{fncylab,enumerate,wrapfig,placeins}
\usepackage{graphicx}
\usepackage{enumitem}  
\usepackage{pifont}
\usepackage{relsize} 
\usepackage{accents}  

\def\urls#1{{\footnotesize\url{#1}}}

\def\mindex#1{\index{#1}}

\def\ocp{*}   

\DeclareFontFamily{U}{mathx}{\hyphenchar\font45}
\DeclareFontShape{U}{mathx}{m}{n}{<-> mathx10}{}
\DeclareSymbolFont{mathx}{U}{mathx}{m}{n}
\DeclareMathAccent{\widebar}{0}{mathx}{"73}

\def\barUpupsilon{\widebar{\Upupsilon}}

\def\SAtime{\uptau}

\def\Obj{\Upgamma}  

\def\Tdiff{\mathcal{D}}

\def\tPR{\text{\tiny\sf  PR}}

\def\SigmaTheta{\Sigma_{\uptheta}}
\def\SigmaPR{\SigmaTheta^{\tPR}}

\def\thetaPR{\theta^{\tPR}}
\def\tilthetaPR{\tilde{\theta}^{\tPR}}

\def\ODEstate{\Uptheta} 

\def\odestate{\upvartheta}

\def\elig{\zeta}

\def\disc{\gamma}

\newcommand{\bbblot}{\raise1pt\hbox{\vrule height .4ex width .4ex depth .05ex}}

\long\def\defbox#1{\framebox[.9\hsize][c]{\parbox{.85\hsize}{%
\parindent=0pt
\baselineskip=12pt plus .1pt      
\parskip=6pt plus 1.5pt minus 1pt 
 #1}}}

\long\def\beginbox#1\endbox{\subsection*{}%
\hbox{\hspace{.05\hsize}\defbox{\medskip#1\bigskip}}%
\subsection*{}}

\def\endbox{}

 \def\archival#1{} 

\def\FRAC#1#2#3{\genfrac{}{}{}{#1}{#2}{#3}}

\def\ddt{{\mathchoice{\FRAC{1}{d}{dt}}%
{\FRAC{1}{d}{dt}}%
{\FRAC{3}{d}{dt}}%
{\FRAC{3}{d}{dt}}}}

\def\ddtp{{\mathchoice{\FRAC{1}{d^{\hbox to 2pt{\rm\tiny +\hss}}}{dt}}%
{\FRAC{1}{d^{\hbox to 2pt{\rm\tiny +\hss}}}{dt}}%
{\FRAC{3}{d^{\hbox to 2pt{\rm\tiny +\hss}}}{dt}}%
{\FRAC{3}{d^{\hbox to 2pt{\rm\tiny +\hss}}}{dt}}}}

\def\ddyp{{\mathchoice{\FRAC{1}{d^{\hbox to 2pt{\rm\tiny +\hss}}}{dy}}%
{\FRAC{1}{d^{\hbox to 2pt{\rm\tiny +\hss}}}{dy}}%
{\FRAC{3}{d^{\hbox to 2pt{\rm\tiny +\hss}}}{dy}}%
{\FRAC{3}{d^{\hbox to 2pt{\rm\tiny +\hss}}}{dy}}}}

\def\half{{\mathchoice{\FRAC{1}{1}{2}}%
{\FRAC{1}{1}{2}}%
{\FRAC{3}{1}{2}}%
{\FRAC{3}{1}{2}}}}

\def\darrow{\buildrel{\rm d}\over\longrightarrow}

\def\limsup{\mathop{\rm lim{\,}sup}}

\def\argmin{\mathop{\rm arg{\,}min}}

\def\state{{\sf X}}

\def\bx{{{\cal B}(\state)}}

\def\bx{{{\cal B}(\state)}}

\def\bfmath#1{{\mathchoice{\mbox{\boldmath$#1$}}%
{\mbox{\boldmath$#1$}}%
{\mbox{\boldmath$\scriptstyle#1$}}%
{\mbox{\boldmath$\scriptscriptstyle#1$}}}}

\def\bfPhi{\bfmath{\Phi}}
\def\bfPsi{\bfmath{\Psi}}

\def\bfDelta{\bfmath{\Delta}}

\def\bfmN{\bfmath{N}}

\def\bfmX{\bfmath{X}}

\def\bfmY{\bfmath{Y}}

\def\bfmhhaY{\bfmath{\hhaY}} 
\def\bfmhhaY{\hbox to 0pt{$\widehat{\bfmY}$\hss}\widehat{\phantom{\raise 1.25pt\hbox{$\bfmY$}}}}

\def\haf{{\hat f}}
\def\hag{{\hat g}}

\def\hagamma{{\hat\gamma}}

\def\haA{\widehat A}

\def\tiltheta{{\tilde \theta}}

\def\tilg{\tilde g}

\def\clB{{\cal B}}

\def\clE{{\cal E}}
\def\clF{{\cal F}}
\def\clG{{\cal G}}

\def\clL{{\cal L}}

\def\clR{{\cal R}}
\def\clS{{\cal S}}
\def\clU{{\cal U}}
\def\clV{{\cal V}}
\def\clW{{\cal W}}
\def\clX{{\cal X}}
\def\clY{{\cal Y}}
\def\clZ{{\cal Z}}

\def\clL{{\cal L}}

\def\eqdef{\mathbin{:=}}

\def\Prob{{\sf P}}

\def\Expect{{\sf E}}

\def\Cov{\hbox{\sf Cov}}

\def\lgmath#1{{\mathchoice{\mbox{\large #1}}%
{\mbox{\large #1}}%
{\mbox{\tiny #1}}%
{\mbox{\tiny #1}}}}

\def\Zero{{\mathchoice{\lgmath{\sf 0}}%
{\mbox{\sf 0}}%
{\mbox{\tiny \sf 0}}%
{\mbox{\tiny \sf 0}}}}

\def\ind{\bbbone}
  
 \def\epsy{\varepsilon}

\def\varble{\,\cdot\,}

\def\formtmp#1#2{{\vskip12pt\noindent\fboxsep=0pt\colorbox{#1}{\vbox{\vskip3pt\hbox to \textwidth{\hskip3pt\vbox{\raggedright\noindent\textbf{#2\vphantom{Qy}}}\hfill}\vspace*{3pt}}}\par\vskip2pt%
\noindent\kern0pt}}

\def\barf{{\widebar{f}}}

\def\derbarf{\bar{A}} 

\def\haderf{\hat{A}}

\def\barg{{\widebar{g}}}

\def\barH{{\bar{H}}}

\def\barL{{\bar{L}}}

{\end{list}}

\def\ass(#1:#2){(#1\ref{#1:#2})}

\def\ritem#1{
\item[{\sf \ass(\current_model:#1)}]
}

\newenvironment{recall-ass}[1]{%
\begin{description}
\def\current_model{#1}}{
\end{description}
}

\def\sq{\hbox{\rlap{$\sqcap$}$\sqcup$}}
\def\qed{\ifmmode\sq\else{\unskip\nobreak\hfil
\penalty50\hskip1em\null\nobreak\hfil\sq
\parfillskip=0pt\finalhyphendemerits=0\endgraf}\fi}

\newcommand{\blot}{\vrule height 1.1ex width .9ex depth -.1ex }
\def\qedb{\ifmmode\blot\else{\vspace{-.2cm}\unskip\nobreak\hfil
\penalty50\hskip1em\null\nobreak\hfil\blot
\parfillskip=0pt\finalhyphendemerits=0\endgraf}\fi}

\newcounter{rmnum}

\newcounter{anum}

\newcommand{\field}[1]{\mathbb{#1}}

\def\Re{\field{R}}

\def\Prob{{\sf P}}

\def\Expect{{\sf E}}

\def\transpose{{\intercal}}

\def\diag{\hbox{\rm diag\thinspace}}

\def\argmin{\mathop{\rm arg\, min}}
\def\ind{\hbox{\large \bf 1}}

\def\trace{\hbox{\rm trace\,}}  
\def\epsy{\varepsilon}
\def\varble{\,\cdot\,}

\def\haH{\widehat{H}}

\def\haL{\widehat L}

\def\haY{\widehat{Y}}

\def\haDelta{\widehat{\Delta}}

\def\hhaY{\hbox to 0pt{$\haY$\hss}\widehat{\phantom{\raise 1.25pt\hbox{Y}}}}

\def\hab{\widehat b}

\def\haA{\widehat A}

\def\haM{{\widehat M}}

\def\haY{\widehat Y}

\def\bfPhi{\bfmath{\Phi}}

\newlength{\dhatheight}

\def\barUpupsilon{\widebar{\Upupsilon}}

\newsavebox{\junk}
\savebox{\junk}[1.6mm]{\hbox{$|\!|\!|$}}

\def\limsup{\mathop{\rm lim\ sup}}

\def\argmin{\mathop{\rm arg\, min}}

\def\clA{{\cal A}}
\def\clB{{\cal B}}

\def\clD{{\cal D}}
\def\clE{{\cal E}}
\def\clF{{\cal F}}
\def\clG{{\cal G}}

\def\clL{{\cal L}}

\def\clR{{\cal R}}
\def\clS{{\cal S}}
\def\clT{{\cal T}}
\def\clU{{\cal U}}
\def\clV{{\cal V}}
\def\clW{{\cal W}}
\def\clX{{\cal X}}
\def\clY{{\cal Y}}
\def\clZ{{\cal Z}}

\makeatletter
\newcommand\gobblepars{%
    \@ifnextchar\par%
        {\expandafter\gobblepars\@gobble}%
{}}
\makeatother

\def\whamrm#1{\smallbreak\pagebreak[3]%
	\noindent\text{\rm#1}\ \ \gobblepars}
	
\def\whamit#1{\smallbreak\pagebreak[3]%
	\noindent\textit{#1}\ \ \gobblepars}

\def\wham#1{\smallbreak\pagebreak[3]%
	\noindent\textbf{#1}\ \ \gobblepars}

\def\bigupgamma{\scalebox{1.2}{\raise1pt\hbox{$\upgamma$}}}
\newtheorem{mytheorem}{Theorem}[section]

\def\barfalpha{\bar{f}_{\!\alpha}}
\def\Aalpha{\derbarf_{\!\alpha}}
\def\Athree{\text{\rm(A3${}^\circ$)}}
\def\AthreeV{\text{\rm(A3${}^\bullet$)}}

\def\sbullet{{\scalebox{0.75}{\textbullet}}}
\def\barftwo{{\bar{f}}}
\def\bargamma{{\bar{\gamma}}}

\def\Ztheta{Z_\uptheta^*}
\def\ZthetaAlpha{Z_\uptheta^\alpha}
\def\ZthetaPR{Z_{\tPR}^*}   %
\def\ZMD{Z_{\MD^*}^*}

\def\YthetaPR{Y_{\tPR}}

\def\thbias{\upbeta_\uptheta}
\def\Zthbias{\upzeta_\uptheta^*}

\def\Probe{\Xi}
\def\tProbe{{\text{\tiny$\Xi$}}}

\def\tilodestate{\tilde\odestate}

\def\alphaTmp#1{\alpha_0^{\text{\rm\tiny\ref*{#1}}}}
\def\bdd#1{b^{\text{\rm\tiny\ref*{#1}}}}

\def\bdde#1{\varrho_{\text{\rm\tiny\ref*{#1}}}}
\def\bddepsy#1{\epsy_{\text{\rm\tiny\ref*{#1}}}}
\def\bddMF{b_\circ}
\def\rhoMF{\varrho_\circ}
\def\Vfor{V_0} 

\def\tilh{\tilde{h}}

\def\tTheta{\text{\tiny$\Theta$}}
\def\tState{\text{\tiny$\state$}}
\def\tilclG{\widetilde{\mathcal{G}}}  

\def\barM{\widebar{M}}

\def\tilG{\widetilde{G}}

\def\haclG{\widehat{\clG}}
\def\barG{\bar{G}}

\def\taucpl{\tau_{\lilc}}
\def\bcpl{b_{\lilc}}
\def\rhocpl{\varrho_{\lilc}}
\def\rglCLT{\rangle_{\text{\tiny \sf CLT}}}
\def\rglLt{\rangle_{L_2}}
\def\ScaledST{\mathcal{P}}
\def\DetST{\mathcal{P}} 
\def\SigmaCLT{\Sigma_{\text{\tiny \sf CLT}}}

\def\Obj{\Upgamma}  %

\def\whamb{\wham{$\bullet$}}

\def\Proof{\whamit{Proof.}}

\def\rd#1{{\color{red}#1}} 
\def\bl#1{{\color{blue}#1}}

\def\state{{\sf X}}
\def\eqdef{\mathbin{:=}}
\def\trace{\hbox{\rm trace\,}}

\def\Cov{\textup{\textsf{Cov}}}  
\def\diag{\,\textup{\rm diag}\,}   

\def\Spx{\textsf{S}}
\def\ind{\field{I}}

\def\Re{\field{R}}

\def\sens{s}
\def\Sens{\clS}

\def\lilsens{\textsf{\textup{\tiny S}}}
\def\lilA{\textsf{\textup{\tiny A}}}
\def\lilc{\textsf{\textup{\tiny c}}}
\def\lile{\textsf{\textup{\scriptsize e}}}

\def\lilb{\textsf{\textup{\tiny b}}}
\def\lilTS{{a}}

\def\lilY{\textsf{\textup{\tiny Y}}}

\def\lilObj{\hbox{\tiny$\Obj$}} 

\def\DistSens{\clD^{\,\lilsens}}
\def\Sense{\Sens^{\lile}}
\def\DistSensND{\clD^{\, \lilsens_a}}
\def\DistSensTS{\clD^{\, \lilsens_b}}
\def\DistSensY{\clD^{\,\lilY}}
\def\DistSensYa{\clD^{\,\lilY a}}

\def\MD{\clW}

\def\Oops{\Upupsilon}
\def\haOops{\widehat{\Upupsilon}}

\def\barOops{\barUpupsilon}
\def\MDA{\MD^{\lilA}}
\def\uppsiA{\uppsi^{\lilm}}

\def\MDb{\MD^{\lilb}} 
\def\uppsib{\uppsi^{\lilb}}

\def\clTA{\clT^{\lilA}}
\def\OopsA{\Upupsilon^{\lilA}}
\def\uppsiA{\uppsi^{\lilA}}

\def\MDSens{\MD^{\,\lilsens}}

\def\clTSens{\clT^{\lilsens}}
\def\OopsSens{\Upupsilon^{\lilsens}}

\newlength{\noteWidth}
\long\def\notes#1{\ifinner
{\footnotesize #1}
\else 
\marginpar{\parbox[t]{\noteWidth}{\raggedright\tiny#1}}  
\fi\typeout{#1}}

\def\rd#1{{\color{red}#1}} 
\def\bl#1{{\color{blue}#1}}

\usepackage{cleveref}

\Crefname{corollary}{Corollary}{Corollaries}
\Crefname{eqnarray}{eq.}{eqs.}
\Crefname{equation}{eq.}{eqs.}
\Crefname{figure}{Fig.}{Figs.}
\Crefname{tabular}{Tab.}{Tabs.}
\Crefname{table}{Tab.}{Tabs.}
\Crefname{proposition}{Prop.}{Propositions}
\Crefname{theorem}{Thm.}{Thms.}
\Crefname{mytheorem}{Thm.}{Thms.}
\Crefname{definition}{Def.}{Defs.} 
\Crefname{section}{Section}{Sections}
\Crefname{lemma}{Lemma}{Lemmas}
\Crefname{assumption}{Assumption}{Assumptions} 

\graphicspath{{./BiasFigures/}}

\def\whamit#1{\smallbreak\pagebreak[3]%
\noindent\textit{#1}\ \ \gobblepars}

\def\wham#1{\smallbreak\pagebreak[3]%
\noindent\textbf{#1}\ \ \gobblepars}

\title{The case for and against fixed step-size:
Stochastic approximation algorithms in optimization
and machine learning}

\author{Caio Kalil Lauand\thanks{Department of ECE, University of Florida, Gainesville.  		{\tt caio.kalillauand@ufl.edu.}} \and Ioannis Kontoyiannis\thanks{
	Statistical Laboratory, University of Cambridge, UK.
	{\tt yiannis@maths.cam.ac.uk.}}\thanks{ 		I.K. 
	was supported in part by 
	the EPSRC funded INFORMED-AI project EP/Y028732/1.
} 	
\and
Sean  Meyn\thanks{Department of ECE, University of Florida, Gainesville.  		{\tt meyn@ece.ufl.edu.}} \thanks{C.K.L. and S.M. received financial support from ARO awards W911NF2010055 and  W911NF2410389,     and NSF award    CCF~2306023. } 
}

\begin{document}

\maketitle

\begin{abstract}%
	Theory and application of stochastic approximation (SA) have become 
	increasingly relevant due in part to applications in optimization
	and reinforcement learning.   This paper takes a new look at the 
	algorithm with  constant step-size  $\alpha>0$,   
	defined by the recursion,
	\[
	\theta_{n+1} = \theta_{n}+ \alpha f(\theta_n,\Phi_{n+1}) \, ,   \qquad \theta_0\in\Re^d\,,
	\]
	in which $\bfPhi = \{ \Phi_{n}\;;\;n\geq 0 \}$ is a Markov chain.   SA algorithms are designed to approximate the solution to the root finding problem $\barf(\theta^*) =0$, where $\barf(\theta) = \Expect[ f(\theta,\Phi) ]$
	and $\Phi$ has the steady-state distribution of $\bfPhi$.
	
	The following conclusions are obtained under an ergodicity assumption 
	on the Markov chain,   compatible assumptions on $f$, and for $\alpha>0$ 
	sufficiently small: 
	
	\wham{1.}  The pair process $\{ (\theta_n,\Phi_n ) \; ;\;  n\ge 0  \}$  
	is geometrically ergodic in a topological 
	sense.
	
	\wham{2.}   
	For every $p\ge 1$,   there is a constant $b_p$ such that   $\limsup_{n\to\infty} \Expect[ \| \theta_n - \theta^*  \|^p] \le  b_p \alpha^{p/2}$ for each initial condition.     In particular, the mean squared error is   of order $O(\sqrt{\alpha})$.  
	
	\wham{3.}  
	The Polyak-Ruppert-style
	averaged estimates $\thetaPR_n = n^{-1} \sum_{k=  1 }^{n} \theta_k$ 
	converge to a limit $\thetaPR_\infty$ 
	almost surely and in mean square,  
	which satisfies $ \thetaPR_\infty  =\theta^* + \alpha \barOops^* 
	+O(\alpha^2)$ for an identified
	non-random $ \barOops^* \in\Re^d$.     Moreover,   the covariance 
	is approximately optimal:  The limiting covariance matrix 
	of $\thetaPR_n$
	is approximately minimal in a matricial sense.
	
	\smallskip
	
	The two main take-aways for practitioners depend on the application.   
	It is argued that, in applications to optimization, constant gain algorithms 
	may be preferable even when the objective has multiple local minima. 
	This is because the average target bias  
	$\frac{1}{n}  \sum_{k=1}^{n}  \barf (\theta_k)  $ 
	vanishes, and   its covariance   is small for large $n$,
	provided exploration is chosen with care.
	In applications to reinforcement learning,   
	the  bias is typically of order $O(\alpha)$, with or without averaging.      In such cases, a   vanishing gain algorithm is preferable.   
\end{abstract}

%
%

\bigskip

\noindent
\textbf{This arXiv version represents a major extension of the results in prior versions:}
Geometric ergodicity of the parameter-noise process is now established for non-linear SA, leading to stronger main results. Moreover, some of the theory is now extended to  non-convex optimization applications.

\bigskip 

\textbf{MSC subject classifications:} 62L20, 68T05

\clearpage

\tableofcontents
\clearpage

\section{Introduction}
The stochastic approximation (SA) algorithm of Robbins and Monro is designed to solve the root finding problem  
\begin{equation}
	\text{$\barf(\theta^*) =0$,  \ in which  }  \ 
	\barf(\theta) \eqdef  \Expect[f(\theta, \Phi)] \,, \ \ \theta\in\Re^d\,,
	\label{e:barf_def}
\end{equation}
for a random variable $\Phi$ taking values in some space $\state$, 
and for a given function $f\colon\Re^d\times\state\to\Re^d$.    
The   SA algorithm is the $d$-dimensional recursion,
\begin{equation}
	\theta_{n+1} = \theta_{n}+ \alpha_{n+1} f(\theta_n,\Phi_{n+1}) \, , 
	\quad 
	n\ge 0 \, ,
	\label{e:SA_recur}
\end{equation}
with initial condition $\theta_0 \in \Re^d$, a nonnegative step-size 
sequence $\{\alpha_n\}$, and under the assumption that $\Phi_n$
converges in distribution to $\Phi\sim\uppi$, as $n\to\infty$.
The sequence  $\bfPhi = \{\Phi_n\}$ is often assumed to be Markovian, 
which is the setting of the present paper.    Convergence theory is 
typically based 
on comparison of 	\eqref{e:SA_recur} with an ODE known as the \textit{mean flow}
ODE:
\begin{equation}
	\ddt \odestate_t  =   \barf(\odestate_t) \,. \quad
	\label{e:meanflow}
\end{equation}

\smallskip

\noindent
{\bf Vanishing step-size.} Much of the relevant
convergence theory is for the case of a vanishing step-size sequence,  with  $\alpha_n = n^{-\rho}$ a common choice.    The constraint  $\rho \in (1/2,1]$ is imposed so that the step-size is square summable, but $\sum_n \alpha_n =\infty$ \cite{kusyin97,benmetpri12,bor20a}.  
Then, almost sure convergence of $\{\theta_n\}$ to $\theta^*$ holds under minimal assumptions on $\barf$ and $\bfPhi$ \cite{bor20a}.  
It is shown in  \cite{laumey24a} that, under mild assumptions,  convergence holds for any value $\rho\in (0,1]$.

Theory for convergence rates remains a research frontier. It is known that   the mean squared error (MSE) decays slowly:  For vanishing gain algorithms, one expects the bound $\Expect [\| \theta_n - \theta^* \|^2] = O(\alpha_n)$.
One approach to speed up convergence is to employ Polyak-Ruppert (PR) averaging \cite{rup88,pol90},
\begin{equation}
	\thetaPR_N := \frac{1}{N-  N_0} \sum_{k= N_0 +1 }^{N} \theta_k \,,
	\quad  
	0 < N_0 < N,
	\label{e:PRdefSA}
\end{equation} 
with $N_0>0$ introduced to discard transients.  

\begin{subequations}
	Subject to conditions on $f$, $\bfPhi$, and the step-size sequence $\{\alpha_n\}$,  averaging achieves two benefits:  The MSE decays at rate $O(N^{-1})$, and the  \textit{asymptotic covariance} is minimal.    That is, 
	with
	\begin{align}
		\SigmaPR & \eqdef 	\lim_{n\to\infty} n \Cov\,(\thetaPR_n) 
		\label{e:SigmaPR}
		\\
		\mbox{and}\quad
		\SigmaTheta^* &=   G^* \Sigma_{\MD^*} {G^*}^\transpose ,
		\label{e:SigmaPRopt}
	\end{align}
	we have $\SigmaPR = \SigmaTheta^*$, 
	in which $G^* \eqdef  - [\derbarf^*]^{-1}$  with $\derbarf^* =\partial_\theta \barf\, (\theta^*)$,  and $ \Sigma_{\MD^*}$ is the ``noise covariance matrix''  defined in  \Cref{t:sensLinear}.   The matrix $G^*  $   defines the
	matrix gain in the stochastic Newton-Raphson algorithm of Ruppert \cite{rup85}.
	The matrix $ \SigmaTheta^* $ is minimal in the matricial sense \cite{rup88,pol90}.       		
	\label{s:SigmaPReqns}
\end{subequations}

%

\smallskip

\noindent
{\bf Fixed step-size.}
Despite this attractive theory for SA with vanishing step-size,
many practitioners advocate a  fixed step-size,    $\alpha_n = \alpha>0$
for all $n$.
There is little hope for convergence of $\{\theta_n\}$  in this case, but bounds on bias and variance can be obtained once boundedness of the recursion is established. In particular, similar to the vanishing step-size case, the MSE is determined by the step-size, $\limsup_{n \to \infty} \Expect [\| \theta_n - \theta^* \|^2] = O(\alpha)$; see \cite{bormey00a} for sufficient conditions.   

\smallskip

\noindent
{\bf Connections and differences.} 
Recent research provides a bridge between theory and practice.   
For example, in the case of fixed step-size,
the papers \cite{bacmou13,durmounausam24} obtain not only convergence of the averaged estimates $\{\thetaPR_N\}$ to $\theta^*$, but also establish the optimal   $O(N^{-1})$ convergence rate for the MSE. 
But these positive results come with a price -- it is assumed that $ \bfPhi$ is an independent and   identically distributed (i.i.d.) sequence. In this paper 
we find that many of these conclusions extend to the case of 
Markovian $\bfPhi$ and fixed step-size when the noise is additive:
\begin{equation}
	\theta_{n+1}    =  \theta_{n} + \alpha  [  \barf(\theta_n) +  \Probe_{n+1}    ]   \,,  \quad n\ge 0,
	\label{e:WarsawStochasticQSA}
\end{equation}
where $\Probe_{n+1} \eqdef f(\theta_n,\Phi_{n+1}) - \barf(\theta_n)$ 
is independent of $\theta_n$, for each $n$. 
This special case  arises frequently in the optimization literature; see \Cref{s:opt}.

However, in general,  for the fixed step-size algorithm there is bias.   
Under suitable conditions, in this work we establish convergence,
\begin{equation}
	\thetaPR_\infty = \lim_{n\to\infty} \thetaPR_n, \quad \mbox{a.s.},
	\label{e:thetaPRinfty}
\end{equation}
and identify the bias $\thbias \eqdef \thetaPR_\infty-\theta^*$.
We show that, 
\begin{equation}
	\thetaPR_N \approx \sqrt{\frac{\alpha}{N}} W_\infty + \thetaPR_\infty\, , \quad W_\infty \sim N(0,\SigmaPR)\,, 
	\label{e:CLTapprox}
\end{equation}
where the approximation is in distribution and the asymptotic covariance   $\SigmaPR$ is approximately equal to $ \SigmaTheta^*$.

\smallskip

\noindent
{\bf Main contributions.} The three main contributions 
of this work are summarised below. Parts A. and B. 
outline the context and contents of the results in
\Cref{t:sens} and \Cref{t:optAssumptions}, respectively.
Part C. describes more precise results for the 
asymptotic bias and variance of the SA estimates
in the special case of linear SA.

\smallskip

\noindent
{\em A. Topological ergodicity.}  
The pair process $\bfPsi=\{ \Psi_n = (\theta_n ,  \Phi_n )\;;\;  n\ge 0\}$ is 
a time-homogeneous Markov chain.  
Under the assumptions of  \Cref{t:sens}, 
it is shown that there is  $\alpha_0>0$ such that the joint process 
is  geometrically ergodic in a topological sense, 
with unique invariant measure   $\upvarpi$,   
provided  the step-size satisfies $0<\alpha\le \alpha_0$.
In particular,
for each 
initial condition $\Psi_0$,
there is a joint realization  
$\{ (\Psi_n,\Psi^\infty_n)  \;;\;  n\ge 0\}$ 
in which $ \Psi^\infty_n  = (\theta^\infty_n,\Phi_n^\infty)$  for each $n$ with $\{\Phi^\infty_n\}$ a stationary process 
with marginal distribution $\uppi$, and
\begin{equation}
	\lim_{n\to\infty} \|  \theta_n - \theta_n^\infty  \| = 0 \quad \mbox{a.s.},
	\label{e:Coupling}
\end{equation}
where the rate of convergence is geometric.  

\smallskip

\noindent
{\em B. Bias and variance.}  
Under the assumptions
of \Cref{t:sens} and 
with $\alpha_0>0$ as before,
there is a fixed constant  $\bdd{t:sens}<\infty$  such 
that the following conclusions hold for $0<\alpha \leq \alpha_0$:
\begin{itemize} 
	\item[(i)]
	\textit{The bias is of order $O(\alpha)$}:  
	For each $\theta_0\in\Re^d$, $\Phi_0\in\state$,
	the target bias is
	\begin{equation}
		\lim_{n\to\infty}
		\Expect[  \barf(\theta_n) ] 
		=
		\Expect_\upvarpi[  \barf(\theta_0^\infty) ] =  	\alpha    \barOops_\alpha\,,
		\label{e:TargetBias}
	\end{equation}
	in which $\theta_0^\infty$ is distributed according to the first marginal of $\upvarpi$,  and   $ \barOops_\alpha\to
	\barOops^*\in\Re^d$ as $\alpha\downarrow 0$ where $\barUpupsilon^* \in \Re^d$ is identified in \eqref{e:OopsApprox}.  Consequently, the 
	parameter bias,
	$\lim_{n\to\infty}
	\Expect[\theta_n]-\theta^*=\Expect_\upvarpi[\theta_0^\infty] -\theta^*$,
	is also of order $O(\alpha)$.

	\item[(ii)]
	\textit{Moment bounds}:    $  \Expect[\|\theta _n   \|^4 ]  $ is uniformly bounded for any initial condition and   
	\begin{equation}
		\lim_{n\to\infty}
		\Expect[\|\theta _n - \theta^* \|^4 ]  =	\Expect_\upvarpi[\| \theta_0^\infty - \theta^* \|^4 ]      \le \bdd{t:sens} \alpha^{2}.
		\label{e:BM2p3}
	\end{equation}
	The steady-state error covariance admits the approximation
	\begin{equation}
		\lim_{\alpha\downarrow 0}  \frac{1}{\alpha} \Expect_{\upvarpi} [  
		(\theta_n-\theta^*) (\theta_n-\theta^*)^\transpose]  =  \Ztheta     \,, 
		\label{e:Sigma0}
	\end{equation} 
	in which  $ \Ztheta $  is identified in  \Cref{t:sens}.

	The next contribution does not require a unique solution to $\barf(\theta) =0$.  In this case we consider 
	the \textit{empirical target bias},  defined by  
	\begin{equation}
		\upbeta_N^\barftwo   \eqdef    \frac{1}{N}  \sum_{n=0}^{N-1}  \barf (\theta_n)  \,,  \quad N\ge 1.
		\label{e:empTargetBias}
	\end{equation}
	
	\item[(iii)]
	\textit{Target bias and variance without convergence}: Under mild conditions, $  \Expect[\|\theta _n   \|^4 ]  $ is uniformly bounded for any initial condition even when the mean flow \eqref{e:meanflow} is not asymptotically stable. Moreover,    
	\Cref{t:optAssumptions} establishes the following
	for the additive noise model \eqref{e:WarsawStochasticQSA}:
	\begin{subequations} 
		\begin{align}
			\lim_{N\to 0}   \upbeta_N^\barftwo &=0 \quad \mbox{a.s.}
			\label{e:NullTargetBias}
			\\
			\lim_{N\to 0}  N  \Expect[ \|  \upbeta_N^\barftwo \|^2  \mid \Psi_0 =z ]  &   =  \trace(\SigmaCLT^{\Probe}).
			\label{e:NullTargetBiasMSE}
		\end{align}
		The matrix $\SigmaCLT^\Probe$ is identified in \Cref{t:optAssumptions},
		and
		conditions are provided under which $\SigmaCLT^{\Probe} =0$.  
		\label{e:TargetBiasResults}
	\end{subequations} 
\end{itemize}

\noindent
{\em C. Finer theory for linear SA.}  The conclusions of  \Cref{t:sensLinear} concern  the linear SA recursion,
\begin{equation}
	\theta_{n+1} = \theta_n + \alpha [A_{n+1} \theta_n - b_{n+1}] ,
	\label{e:LinSAintro}
\end{equation}
in which $A_n$ and $b_n$ are functions of $\Phi_n$, for each $n$.   It is shown that the coupling result \eqref{e:Coupling} holds in the mean square sense,
\begin{equation}
	\lim_{n\to\infty} \Expect[  \|  \theta_n - \theta_n^\infty  \|^2]   = 0,
	\label{e:CouplingLin}
\end{equation} 
where again the rate of convergence is geometric; see \eqref{e:L2LyapExp}.
The target bias limit \eqref{e:TargetBias} remains valid, 
and the parameter bias admits the representation,  
\begin{equation}
	\thbias \eqdef 	\lim_{n\to\infty}   \Expect[  \theta_n    ]   - \theta^*   =	\lim_{n\to\infty}   \Expect[  \thetaPR_n    ]   - \theta^*   
	=  \alpha \Zthbias + O(\alpha^2) ,
	\label{e:Bias}
\end{equation}
in which  $\Zthbias  \in \Re^d$  is identified in \Cref{t:sensLinear}.
Moreover,  the  $d\times d$ matrix $\ZthetaPR$  is identified in the 
approximation:
\begin{equation}
	\SigmaPR \eqdef 	\lim_{n \to\infty}  n \Cov\,(\thetaPR_n) = \SigmaTheta^*   +   \alpha   \ZthetaPR  + O(\alpha^2).
	\label{e:PRCov_Markov}
\end{equation}


\wham{Prior work and related approaches.}
The results of the present paper are most closely related 
to \cite{huakonmey02a}, which also treats SA with constant step-size. 
The conditions imposed are difficult to compare to those in the present work. 
The assumptions on $\bfPhi$ are weaker: simple 
geometric ergodicity rather than the drift condition (DV3).
The conditions on $f$ are far stronger:  it is assumed 
that $  f(\theta,x) =\barf(\theta) + \Delta(\theta,x)$,   in which 
$ \Delta(\theta,x)  \le B(x)$ for all $\theta, x$ and for a function $B$  
satisfying appropriate bounds associated with a Lyapunov function 
for $\bfPhi$.      
In the case of linear SA 
as in \eqref{e:LinSAintro}, the assumption on $f$ is
relaxed somewhat:
it is assumed that $f(\theta,x) =A(x)\theta-b(x)$,
where the matrix-valued function $A$ has entries that are uniformly  
bounded. Under these assumptions, the coupling 
result  \eqref{e:CouplingLin}  is established, and extended to 
convergence of $p$th moments for any $p\ge 2$.
Strong assumptions are indeed necessary to obtain 
such moment bounds 
-- see \cite{borchedevkonmey25} for a counterexample based 
on linear SA with vanishing gain.  

The paper \cite{huakonmey02a} obtains bounds on the MSE
but does not consider bias. 
A bias representation for SA with constant gain first appeared in 
\cite{laumey22b,laumey22c}.  
The more recent work in  \cite{huochexie22,huozhachexie24} also treats fixed step-size SA with Markovian noise but restricts to the linear case, while \cite{huozhachexie24,allgas24} allow for nonlinear recursions. 
Although \cite{allgas24} concerns parameter dependent disturbances (i.e., $\Phi_{n+1}$ depends upon $\theta_n$) as found in $\epsilon$-greedy exploration policies for reinforcement learning, this  paper assumes stability of the general recursion in terms of boundedness of parameters. The work \cite{huochexie22} proves stability of the general recursion, but this conclusion is subject to a stronger \textit{uniform ergodicity} assumption on $\bfPhi$. The stronger condition is imposed (and is necessary)
when finite-time bounds for the algorithm are sought.

Literature concerning bias for algorithms with vanishing step-size is much more scarce. We are not aware of any other work apart from \cite{laumey24a}, which mostly studies linear SA. The only result for the nonlinear case is also from \cite{laumey24a} and  concerns target bias.


Analysis of the SA recursion \eqref{e:SA_recur} commonly begins with 
its interpretation as a ``noisy'' Euler approximation of the mean flow \eqref{e:meanflow}.
With ``noise'' or ``disturbance''  $\Delta_{n+1} \eqdef f(\theta_n,\Phi_{n+1}) - \barf(\theta_n)  $, we may write 
\begin{equation}
	\theta_{n+1} = \theta_{n}+ \alpha_n [\barf(\theta_n) + \Delta_{n+1}].
	\label{e:noisyEuler}
\end{equation}
Our starting point of analysis in this work is the decomposition  of  {M\'etivier} and {Priouret}~\cite{metpri87},  
\begin{equation}
	\Delta_{n+1} =  \MD_{n+1} - \clT_{n+1} + \clT_{n} - 
	\alpha_n\Oops_{n+1} \, ,
	\label{e:DeltaDecomp}
\end{equation}
in which $\{  \MD_{n+1} \}$ is a martingale difference sequence
and the processes $\{\clT_n\}$ and $\{\Oops_n\}$ are
defined in Lemma~\ref{t:noise-decomp}.  
This prior work considered the case of vanishing step-size, and based on this representation established convergence of the SA algorithm.   The idea has been applied in many other papers, and in particular leads to a functional central limit theorem under suitable conditions.  The weakest conditions to date are found in \cite{borchedevkonmey25},  and this recent prior work is a major foundation of the present paper.    Key lemmas from \cite{borchedevkonmey25} extend to the setting of this paper,  which form  components of the proof of  \cref{e:BM2p3,e:TargetBias}. 
Finer results, such as the final set of conclusions in \cref{e:Bias,e:PRCov_Markov}, require multiple applications of the disturbance decomposition \eqref{e:DeltaDecomp} to general functions of $(\theta_n,\Phi_{n+1})$. 

Many of the main results concerning bias and 
variance may be cast as first-order Taylor series approximations for the invariant measure $\upvarpi$ described above, viewed as a function of $\alpha$.       For a given function $g\colon \Re^d\times\state \to \Re$,   let $\eta_\alpha =  \Expect[ g(\Psi_n)]$ where $\Psi_n \sim \upvarpi$ and the SA recursion uses step-size $\alpha>0$.    A representation for the derivative is obtained 
from    \cite{sch68} in terms of a family of zero-mean solutions $\hat{g}_\alpha$ to \emph{Poisson's equation}:
\begin{equation}
	\Expect[    \hag_\alpha  (\Psi_{n+1} )  -\hag_\alpha  (\Psi_{n} )  \mid  \clF_n]   =       -   g(\Psi_n)   +  \eta_\alpha,
	\label{e:fish68}
\end{equation}
where $\clF_m = \sigma\{  \Psi_k : k\le m \}$ for  $m\ge 0$.  

Given
the SA recursion \eqref{e:SA_recur},  for a large class of functions 
$H\colon \Re^d\times\state \to \Re$,  and $z=(\theta,x)$,
we have:
\[
\begin{aligned}
	P_\alpha H\, (z) & \eqdef  \Expect[    H  (\Psi_{n+1} )    \mid   \Psi_n =z]  =    \Expect \big[    H  (\theta+ \alpha f(\theta,\Phi_{n+1})    , \Phi_{n+1} )    \mid   \Phi_n =x \big]   \,.
\end{aligned}
\]
Consequently, there is a large class of functions $H$ for which the 
derivative of $P_\alpha$ with respect to $\alpha$ exists 
and can be expressed,
\[
\begin{aligned}
	P'_\alpha H\, (z)  & \eqdef    \tfrac{d}{d\alpha}   \Expect[    H  (\Psi_{n+1} )     \mid  \Psi_n =z ]   
	\\
	&	 =  \Expect \big[  [  \partial_\theta H  (\theta+ \alpha f(\theta,\Phi_{n+1}) , \Phi_{n+1} )   ] f(\theta,\Phi_{n+1})   \mid   \Phi_n =x \big].
\end{aligned}
\]
Schweitzer's ceberated sensitivity theorem~\cite{sch68} gives the following formula:
\begin{equation}
	\tfrac{d}{d\alpha} \eta_\alpha   =  \Expect\big[ P'_\alpha \hag_\alpha \, (\Psi_n) \big ] .
	\label{e:sch68}
\end{equation}
This formula is behind the proof that the actor-critic algorithm can be designed to construct unbiased gradient samples \cite{kontsi03a,CSRL}.

Given the special structure of $P'_\alpha$, Schweitzer's formula has an appealing form at $\alpha =0$:
\begin{equation}
	\tfrac{d}{d\alpha} \eta_\alpha \big|_{\alpha=0}   =  \Expect\big[   [ \partial_\theta\hag_\alpha \, (\theta_n,\Phi_{n+1}) ]  f (\theta_n,\Phi_{n+1})    \big ]  \big|_{\alpha=0} .
	\label{e:sch68zero}
\end{equation}
However,
the precise approach of  Schweitzer cannot be employed here,
for two important reasons.   The first is that we are unable in general 
to ensure that solutions $\{\hag_\alpha   \}$ to Poisson's
equation exist for the functions $\{g_\alpha\}$ of interest.   
The second is the lack of uniqueness of $\upvarpi$ when $\alpha =0$,  
which makes \eqref{e:sch68zero} ill-defined.  

We modify Schweitzer's approach, considering a   function  $ \hag $ independent of $\alpha$  satisfying, for each fixed  $\theta$,  
\begin{equation}
	\Expect[    \hag(\theta,\Phi_{n+1} ) - \hag(\theta,\Phi_{n} )   \mid  \clF_n]   =      -   g(\theta, \Phi_n)   + \barg(\theta) ,
	\label{e:fish68-MP}
\end{equation}
with $\barg(\theta) = \int g(\theta, x) \, \uppi(dx)$,
where $\uppi$ is the unique invariant measure of $\bfPhi$.
For each $\theta$,   (\ref{e:fish68-MP})
may be regarded as the solution to Poisson's equation in the form  \eqref{e:fish68}
with $\alpha =0$ and $\theta_n\equiv \theta$.

Poisson's equation in the form 
\eqref{e:fish68-MP}
is in fact the starting point in the derivation of the disturbance decomposition  \eqref{e:DeltaDecomp},   using $g=f$.    
It will be used in this paper for many different choices of $g$ to obtain many approximations in analogy with \eqref{e:sch68}.   
For example, this approach is used to approximate $\SigmaPR - \SigmaTheta^* $ as defined in 	\eqref{s:SigmaPReqns}.
The error term $\barOops_\alpha$ appearing in  \eqref{e:TargetBias} is constructed in this way,
and its approximation  $\barOops_\alpha  = \barOops^* + O(\alpha)$ results in an expression similar to \eqref{e:sch68zero}:
\begin{equation}
	\barOops^*\eqdef  - \Expect_\uppi\big[  [ \partial_\theta   \uppsi\, (  \theta^*, \Phi_{n+1}) ] f(\theta^*, \Phi_{n})    \big].
	\label{e:OopsApprox}
\end{equation}
with $\uppsi = f - \haf$. See   \Cref{t:sens} and its proof for details.

\wham{Organization.} 
The remainder of the paper is organized in four sections. 
\Cref{s:curse} and \Cref{s:nocon}
summarize the assumptions and notation used throughout the paper, 
and contain the results summarized in contributions A--C above. 
The theory establihsed in these sections is
illustrated through several numerical examples in 
\Cref{s:ex}.
In particular, in application to stochastic gradient descent (SGD) for non-convex optimization, it is seen that the theory of this paper provides insight for exploration in constant-gain algorithms, leading to better performance. However, in application to temporal difference (TD) learning, the theory developed implies that diminishing step-size algorithms might be preferred to combat the negative effects of memory. 
Conclusions and open paths for future research are discussed in 
\Cref{s:Concl}.
Technical proofs   are contained in the Appendix.

\section{Main Results}
\label{s:curse}


We first set the stage for analysis,  going over assumptions and notation. 

\subsection{Preliminaries}
\label{s:assump}

It is assumed that $\bfPhi \eqdef \{\Phi_n\}$ in \eqref{e:SA_recur} is a geometrically ergodic Markov chain,  whose    state space $\state$ is a locally compact and separable metric space equipped with its Borel
$\sigma$-algebra $\bx$. Its transition kernel is denoted $P$ and its unique invariant measure $\uppi$,  so that $\barf(\theta) = \Expect_\uppi[f(\theta, \Phi_k)]$. The subscript in the expectation defining $\barf$ indicates that it is taken in steady-state: $\Phi_k \sim \uppi$.

The following  assumptions are in place throughout:

\wham{(A1)} 
The SA recursion \eqref{e:SA_recur}
is considered with $\alpha_n \equiv \alpha >0$. 

\wham{(A2i)}    There is a function $L\colon\state\to\Re$ satisfying $\| f(\theta,x) - f(\theta',x) \| \le L(x) \|\theta -\theta'\|$   and $	\| f(0, x)\| \leq L(x)$   for all $x\in\state$ and $\theta,\theta'\in \Re^d$.    

\wham{(A2ii)}   
$\bfPhi$ is  an  aperiodic Markov chain satisfying \textbf{(DV3)}:
\\
For  functions $V\colon\state\to\Re_+$,  $ W\colon\state\to [1, \infty)$, 
a small set $C$, $b>0$ and all $x\in\state$,
\begin{equation}
	\Expect\bigl[  \exp\bigr(  V(\Phi_{k+1})      \bigr) \mid \Phi_k=x \bigr]  
	\le  \exp\bigr(  V(x)  - W(x) +  b \ind_C(x)  \bigl).
	\label{e:DV3}
\end{equation}
In addition, 
$S_W(r)  := \{ x :  W(x)\le r \}$ is either small or empty
and $ \sup\{ V(x) :  x\in S_W(r) \}  <\infty $ for each $r\ge 1$;
see \cite{MT} for  the definition of a small set.

Assumption~(A3)  in either of the forms below implies ultimate boundedness of the mean flow.   The latter condition is more common in the
control theory literature.

Letting $\barf_{\infty}(\theta)\eqdef \lim_{c \to \infty} \frac{1}{c} \barf(c\theta)$, we define the  ODE@$\infty$ as $\ddt \odestate^\infty_t = \barf_{\infty}(\odestate^\infty_t)$.

\wham{(A3${}^\circ$)}   
The scaled vector field $\barf_{\infty}(\theta)$ exists for each $\theta \in \Re^d$ and the ODE@$\infty$ is globally asymptotically stable.

\wham{(A3${}^\bullet$)}   
For a globally Lipschitz continuous, $C^1$ function $\Vfor \colon\Re^d\to\Re_+$ and  constant $\varrho_v >0$, the following two bounds hold for any solution to the mean flow,  and any $t$ for which   $\|\odestate_t\|\ge  {1}/{\varrho_v }$:
\begin{equation}  
	\Vfor (\odestate_t)    \ge \varrho_v  \|\odestate_t \|
	\quad\textit{and}\quad
	\frac{d^+}{dt} \Vfor (\odestate_t) 
	\leq  -  \varrho_v   \Vfor (\odestate_t),
	\label{e:ddt_bound_LyapfunTmp}
\end{equation}
where  ``$+$'' indicates right derivative.

\wham{(A4)} $\displaystyle 
\lim_{r \to\infty} \sup_{x \in \state} \frac{L(x)}{\max\{r,W(x)  \}}  =0$.

\wham{(A5)}   The mean flow \eqref{e:meanflow} is globally asymptotically stable, with $\theta^\ocp \in\Re^d$ its unique stationary point;
$\barf:\Re^d \to \Re^d$ is continuously differentiable in $\theta$, and $\derbarf = \partial \barf$ is  uniformly bounded and uniformly Lipschitz continuous.     
Moreover,  $  \derbarf^\ocp \eqdef \derbarf(\theta^\ocp)$ is Hurwitz.

Boundedness of the Jacobian matrix $\derbarf $ is  immediate under (A2).

The mean flow \eqref{e:meanflow} is  exponentially asymptotical stable under (A3) and (A5)  \cite[Prop. A.11]{laumey25a} (further results are obtained in  \cite{vid22} under stronger conditions).   Geometric ergodicity of $\bfPhi$ follows from (A2)---see \cite{konmey03a,konmey05a}  or \cite[Chapter 20.1]{MT}. 

A relaxation of (A4) may be found in  \cite{borchedevkonmey25}.

The two versions of (A3) are nearly equivalent.

\begin{proposition}
	\label[proposition]{t:ODEatInftAndV}
	Suppose that $\barf$ is Lipschitz continuous. If \Athree\ holds then there exists a solution to \AthreeV.   Conversely,  if 	\AthreeV\ holds, and the scaled vector field $\barf_{\infty}(\theta)$ exists for each $\theta \in \Re^d$,
	then \Athree\ holds.  
\end{proposition}

\Proof  
The implication \Athree\ $\Rightarrow$ \AthreeV\ is a small extension of the proof of \cite[Prop.~4.22]{CSRL}.   
The Lyapunov function can be taken to be 
\[
V_0(x) =  \int_0^T  \sqrt{1 +  \|\odestate_t \|^2} \,  dt  \,,  \qquad x = q_0\in \Re^d\,,
\]
where $T>0$ is chosen so that  $\| q_T \| \le \half \|q_0\| = \|x\|$  for all sufficiently large $\|x\|$.  The existence of $T$ is established in \cite[Prop.~4.22]{CSRL}. 

The converse follows similar arguments.   If \AthreeV\ holds then there exists $\varrho>0$ and $B<\infty$ such that   
$ \|\odestate_t \|  \le \exp(-\varrho t)  \|\odestate_0 \|  + B$ for any  $\odestate_0$ and any $t$.    From the definitions it follows that a similar bound holds for the ODE@$\infty$:
$ \|\odestate^\infty_t \|  \le \exp(-\varrho t)  \|\odestate^\infty_0 \| $ for any  $\odestate^\infty_0$ and any $t$, which implies  \Athree.
\qed 

\smallskip

\wham{Assumption (A2) and uniform ergodicity:} A Markov chain is called 
uniformly ergodic if it is $\psi$-irreducible, aperiodic, and the state 
space   $\state $ is small (the assumption imposed for linear SA in  
\cite{huochexie22}).    An example is a finite state space Markov chain 
that is aperiodic with a single recurrent class \cite{MT}.     
Uniform ergodicity is a stronger assumption, rarely satisfied
in applications. It implies geometric ergodicity, and in particular the 
existence of a unique invariant probability measure $\uppi$.
It also implies the Lyapunov condition
(DV3) used in this work.


A  solution to (DV3)  is obtained with the Lyapunov function   $V(x)= 1$  for $x \in \state$.     
Then  $ \Expect\bigl[  \exp\bigr(  V(\Phi_{k+1})      \bigr) \mid \Phi_k= x \bigr]   = 1 $   for any  $x\in\state$.   
Hence (DV3) holds with $W= 1$,   $b = 1$,  and $C=\state$.

\wham{Notation:} When an invariant measure $\upvarpi$ exists for the Markov chain $\bfPsi := \{   \Psi_n =  (\theta_n,\Phi_n): n \geq 0\}$,  its second marginal is the invariant measure 
$\uppi$ for $\bfPhi$.  
All functions are assumed to be measurable with respect to the Borel 
$\sigma$-algebra over their domain. 
For functions $g,h:\Re^d \times \state \to \Re^d$, we denote $\tilg(z) = g(z)-\upvarpi(g)$,  
with $\upvarpi(g) = \int g(z) \upvarpi(dz)$,  and
\[
\barg(\theta) = \int g(\theta,x) \uppi(dx)\, .   
\]
We adopt the following  $L_2$ notation:
\begin{equation}
	\begin{aligned}
		\langle g, h \rglLt 
		& = 
		\Expect_\upvarpi[g(\Psi_0) h(\Psi_0)  ^\transpose]\, , 
		&\Sigma_{g} &= \langle g, g \rglLt,
		\\
		\langle g, h \rglCLT 
		& = 
		\sum_{k = -\infty}^\infty \Expect_\upvarpi[\tilg(\Psi_0) \tilh(\Psi_k)  ^\transpose] \, ,
		&\SigmaCLT^{g} &= 		\langle g, g \rglCLT. 
	\end{aligned}
	\label{e:SigCLTdef}
\end{equation}
When the infinite sum exists, then  $\SigmaCLT^{g} \ge 0$ is known as the \textit{asymptotic variance}---this is the variance appearing in the central limit 
theorem for partial sums of $\{ g(\Psi_k) : k\ge 0\}$    \cite{MT}.   Under the conditions to be imposed we have the alternative representation,
\[
\SigmaCLT^{g}  =\lim_{n\to\infty} \frac{1}{ n } \Cov(S_n^g )  \,,\quad S_n^g  = \sum_{k=1}^n g(\Psi_k).
\]

The representation  \eqref{e:DeltaDecomp}  first appeared in  \cite{metpri87}, based on the following construction:
for each fixed $\theta$,   let    $\haf(\theta,\varble)$ denote the solution to  Poisson's equation with forcing function $f$.     Precisely as in  \eqref{e:fish68-MP},
\begin{equation}
	\Expect[\haf(\theta, \Phi_{n+1}) - \haf(\theta, \Phi_{n})\mid \Phi_n = x] 
	= 
	- f(\theta, x) + \barf(\theta).
	\label{e:fish-f}
\end{equation} 
The solution    $\haf$ is unique up to an additive constant under the assumptions imposed on $f$ and $\bfPhi$.    

\begin{proposition}
	\label[proposition]{t:bounds-H}	If {\em (A2)} holds then    $\haf \colon\Re^d\times\state\to\Re^d$ exists solving 	\eqref{e:fish-f},  with   
	$\Expect_\uppi[\haf(\theta,\Phi_n)] =\Zero$ for each   $\theta \in\Re^d$, 
	and  for a constant $b_f$  and all $\theta,\theta', x$:
	
	\whamrm{(i)}  
	$\| \haf(\theta, x)\| \leq b_f[1+V(x)] \bigl[ 1+\|\theta\| \bigr]$

	\whamrm{(ii)}  
	$\| \haf(\theta, x)  - \haf(\theta', x)  \| \leq b_f[1+V(x)]  \|\theta - \theta'\| $  .
\end{proposition}

\Proof  (DV3)  together with Jensen's inequality gives
$
\Expect\bigl[    V(\Phi_{k+1})      \mid \Phi_k=x \bigr]  \le   V(x)  - W(x) + b \ind_C(x)  
$.  
The result then follows from   \cite[Thm.~17.4.2]{MT}.  \qed

\smallskip

\begin{subequations}
	
	Denote   $\uppsi(\theta,x) \eqdef  f(\theta,x) - \haf(\theta,x)$  for $\theta\in\Re^d$,  $x\in\state$.  
	\begin{lemma}
		\label[lemma]{t:noise-decomp}
		Equation~\eqref{e:DeltaDecomp} 
		with $\alpha_n=\alpha$ holds under (A2), with   
		\begin{align}
			\MD_{n+1} &  \eqdef  \haf(\theta_n, \Phi_{n+1})  - \Expect[  \haf(\theta_n, \Phi_{n+1})   \mid \clF_{n}]    
			\label{e:MD}
			\\[0.5em]
			\clT_{n+1} &\eqdef   \uppsi(\theta_{n+1}, \Phi_{n+1}) 
			\label{e:telescope}
			\\
			\Oops_{n+1} &\eqdef   \frac{1}{\alpha}  \big[   \uppsi(\theta_{n+1}, \Phi_{n+1})   -   \uppsi(\theta_{n}, \Phi_{n+1})     \big] \,, 
			\label{e:Upsilon}
		\end{align}
		where $\clF_m = \sigma\{  \Psi_k : k\le m \}$ for  $m\ge 0$.
		\qed
	\end{lemma}
	
	\label{e:noise-decomp}
\end{subequations}

Approximations of $\SigmaCLT^\Delta$ and $\SigmaPR$ involve the two noise covariances,  
\begin{equation} 
	\Sigma_{\MD}  =  \Expect_\upvarpi [ \MD_k \MD_k^\transpose]  \,, \quad   
	\Sigma_{\MD^*}  = 
	\Expect_\uppi [ \MD^*_k {\MD_k^*}^\transpose] ,
	\label{e:SigmaW}
\end{equation} 
where $\MD^*_{n+1}   \eqdef  \haf(\theta^* , \Phi_{n+1})  - \Expect[  \haf(\theta^*, \Phi_{n+1})   \mid \clF_{n}]    $.   
Because these are martingale difference sequences, the covariance coincides with the asymptotic covariance: $ \SigmaCLT^{\MD} =   \Sigma_{\MD}$ and $ \SigmaCLT^{\MD^*} =   \Sigma_{\MD^*}$.

The martingale difference sequence  $\{ \MD_{n+1} \}$ dominates the variance of $\{\thetaPR_n \}$ for small $\alpha$, while  the  sequence $\{ \Oops_{n+1} \}$ introduces bias  and impacts variance to a lesser degree---see \Cref{t:sensLinear}. In implementations where the ``noise'' sequence $\bfDelta$ is a Markov chain independent of the parameter, as in the additive noise model \eqref{e:WarsawStochasticQSA}, $\haf$ is only a function of $\bfPhi$, yielding  $\Oops\equiv 0$.



\subsection{Ergodicity and Lyapunov exponents}
\label{s:ergo}

In view of the fact that $\bfPsi$ is a Feller Markov chain,  if 
$\sup_n  	\Expect[\|\theta _n  \|^p ]  <\infty$ for just one initial condition $\Psi_0$ then there exists  a stationary version of the pair process $\bfPsi^\infty = \{ \Psi^\infty_n  : -\infty <n < \infty\}$  \cite[Theorem 12.1.2]{MT}.    

Moment bounds  on the error sequence $\{\tiltheta_n = \theta_n -\theta^*  : n\ge 0\}$ 
are obtained in \cite{borchedevkonmey25} for vanishing step-size algorithms by establishing a Lyapunov drift condition for the function $\clV$, defined by 
\begin{equation}
	\begin{aligned} 
		\clV(\theta,x)   &=  (1 +  \beta_0  \|\theta\|^4) v_+(x),
		\\
		\text{with} \ \ 	v_+( x)   &= \Expect[    \exp\bigl( V( \Phi_{k+1} ) + \epsy_0  W(\Phi_k )  \bigr)     \mid \Phi_k =x   ],
	\end{aligned} 
	\label{e:clV}
\end{equation}
and   $\beta_0 , \epsy_0 \in (0,1)$ small constants.    The drift condition is extended to fixed step-size SA in \Cref{t:results_CLT}.  

To establish uniqueness of $\bfPsi^\infty $ and a form of ergodicity for the 
pair process $\bfPsi$ we apply theory associated with
the \textit{sensitivity process},    defined by
$\{\Sens_n :=  \partial_{\theta_0} \theta_n \;;\;n\geq 0\}$, 
in which   $ \theta_n $ is viewed as a smooth function of the initial condition $\theta_0$.
Its definition requires common randomness  in \eqref{e:SA_recur} for each initial condition,  in which case we obtain the evolution equations  
\begin{equation}
	\Sens_{n+1}  =    \Sens_n +   \alpha A_{n+1} \Sens_n   \,, \quad   A_{n+1} \eqdef \partial_\theta f \, (\theta_n,\Phi_{n+1}).
	\label{e:Sens1}
\end{equation}
The $L_p$-Lyapunov exponent $\Uplambda_p$ is   the growth rate,   
\begin{equation}
	\Uplambda_p \eqdef\lim_{n\to\infty} \frac{1}{n} \log( \Expect[ \|  \Sens_n \|_F^p]^{1/p} ),
	\label{e:LpLyapExp}
\end{equation}
where $\|\varble \|_F$ denotes the Frobenious norm.   If for each initial condition this limit exists and is negative, then parameter sequences from distinct initial conditions converge to a steady-state in a topological sense.   

We can view the joint process $\{(\theta_n,\Sens_n)\}$  as an SA recursion with fixed step-size.   The mean flow  associated with the sensitivity process is the matrix ODE:  
\begin{equation}
	\ddt \sens_t = \derbarf(\odestate_t) \sens_t  \,,\qquad \sens_0 = I \,.
	\label{e:sensODE}
\end{equation}
Unfortunately,  outside of the special case of linear SA,  the right hand side of 
\eqref{e:sensODE}  is not jointly Lipschitz continuous  in  $(\odestate_t ; \sens_t )$ unless   $ \derbarf(\odestate_t) \equiv \derbarf^*$  (as holds in linear SA).   For nonlinear SA we obtain bounds on an almost-sure version of the Lyapunov exponent, which is a major component of the proof of our first main result:    

\begin{mytheorem}
	\label[mytheorem]{t:sens}  Suppose that {\em (A1)--(A5)} hold.   
	Then, there is $\alpha_0 \in (0,1)$ and $\bdd{t:sens}<\infty$  such that the following hold when 
	the  step-size satisfies $0 < \alpha\le \alpha_0$,   for any initial condition  $z=(\theta;x) =\Psi_0$:
	
	\begin{subequations} 
		
		\whamrm{(i)}  There is a unique invariant measure $\upvarpi$ for the 
		Markov chain $\bfPsi$ satisfying,
		\begin{align}
			\limsup_{n\to\infty} \Expect[  \|    \tiltheta_n  \|^4]  &    \le \alpha^2 \bdd{t:sens} 
			\label{e:theta4th-moment}
			\\
			\lim_{n\to\infty} \Expect[       \tiltheta_n \tiltheta_n^\transpose ]  & =   \Expect_\upvarpi[    \tiltheta_0\tiltheta_0^\transpose ]    = \alpha  \Ztheta +O(\alpha^2)\,,
			\label{e:theta2th-moment}
		\end{align}
		with $\Ztheta\ge 0$ the   unique positive semidefinite solution to the Lyapunov equation,
	\end{subequations} 
	\begin{subequations} 
		\begin{align}
			0 &=    \derbarf^* \Ztheta  +  \Ztheta { \derbarf^*}^\transpose +    \SigmaCLT^{\Delta^*}   
			\label{e:Sigma0LyapEqn}
			\\[.5em]
			&\text{in which} \ \   
			\SigmaCLT^{\Delta^*}  =   \Expect_{\upvarpi}  \big[ 
			\Delta_n^* \haDelta_n^{* \transpose}    
			+    \haDelta^*_n    \Delta_n^{* \transpose}    
			- \Delta_n^*  { \Delta_n^* }^\transpose  
			\big ],
			\label{e:Ztheta}
		\end{align}  
		and $  \haDelta^*_n = 
		\haf(\theta^*, \Phi_n )$.
	\end{subequations} 
	
	\whamrm{(ii)} 
	The bias satisfies $\| \thbias \| \le \alpha \bdd{t:sens}$.   Moreover, 
	the limit \eqref{e:TargetBias} holds   with 
	$\displaystyle
	\barOops_\alpha \eqdef       \Expect_\upvarpi  [ \Oops_{n+1} ] = \barOops^* + O(\alpha)$,  where $ \barOops^*$ is defined in \eqref{e:OopsApprox}.

	\whamrm{(iii)}    There is a   realization of $\bfPsi$ together with  a stationary version $\bfPsi^\infty = \{ \Psi^\infty_n  \;;\; -\infty <n < \infty\}$ for which  \eqref{e:Coupling}
	holds in the following strong sense:   there is $\bdde{t:sens} >0$, independent of $\Psi_0$, such that
	\begin{equation}
		\limsup_{n\to\infty} \frac{1}{n} \log  \|  \theta_n - \theta_n^\infty  \|   \le - \bdde{t:sens} \alpha     \quad \mbox{a.s.}
		\label{e:asExp}
	\end{equation}
\end{mytheorem}

Slightly stronger conclusions are obtained under more restrictive noise conditions. 
\begin{corollary}
	\label[corollary]{t:corr_sens}
	Suppose that the assumptions of \Cref{t:sens} hold. Suppose moreover that either of the following conditions hold:
	{\rm (a)} The noise is additive, as in \eqref{e:WarsawStochasticQSA}, or 
	{\rm (b)}  The process $\bfDelta$ is a martingale difference sequence.
	
	Then, the conclusions of \Cref{t:sens} hold. Moreover, the right hand side of \eqref{e:TargetBias} is zero.
	\qed
\end{corollary}

\subsection{Linear stochastic approximation}
\label{s:lin}

Consider the SA recursion composed of $\{\theta_n \}$ and a second $d$-dimensional sequence $\{\Sense_n \}$ that will be related to the sensitivity process 	\eqref{e:Sens1}.  For small $\delta_s>0$, denote   $M_\delta = M_\delta(\alpha) \eqdef \alpha^{-1} [  \exp(\delta_s\alpha)-1]I $,  and   for $(\theta_0;\Sense_0) \in\Re^{2d}$,
\begin{subequations}
	\begin{align}
		\theta_{n+1} &= \theta_n + \alpha [A_{n+1} \theta_n - b_{n+1}]  
		\label{e:LinSAgen}
		\\
		\Sense_{n+1} & =    \Sense_n +    \alpha \bigl[   M_\delta +  A_{n+1} \bigr]\Sense_n .
		\label{e:LinSens2}
	\end{align}
	We have  $ \Sense_n  =  \exp(n\delta_s\alpha) [\Sens_n]^i$ for the initial condition  $\Sense_0 = e^i$ (the $i$th basis vector in $\Re^d$),  where  
	$ [\Sens_n]^i$ denotes the $i$th column of  $\Sens_n$.
	
	\label{s:SAlinearBivariate}
\end{subequations}

Analysis is greatly simplified because  $A_{n+1} = A(\Phi_{n+1})$  and  $b_{n+1} = b(\Phi_{n+1})$ do not depend upon $\theta_n$, giving   $ \Sense_n  =  \exp(n\delta_s\alpha) A_n\cdots A_1 \Sense_0$. 

\begin{subequations}
	
	Let $\derbarf^*$, $b^*$ denote the  respective means of  $A$, $b$ under $\uppi$,  and  
	$\haA$,  $\hab$ denote the respective zero-mean solutions to Poisson's equation:
	\begin{align}
		\Expect[\haA(\Phi_{n+1} )  \mid \Phi_n =x ]  &  =  \haA(x)  -  [A(x)  - \derbarf^*]
		\label{e:haA}
		\\
		\Expect[\hab(\Phi_{n+1} ) \mid \Phi_n =x ]  &  =  \hab(x)  -  [b(x)  - b^*].
		\label{e:hab}
	\end{align}
	We then have $ \haf(\theta_{n}, \Phi_{n+1}) = \haA_{n+1} \theta_n -\hab_{n+1} $.
	In \Cref{t:noise-decomp-linear}~(i) the disturbance decomposition introduced by \Cref{t:noise-decomp} is simplified,
	and in (ii) it is refined:
	
	\label{e:haAb}
\end{subequations}

\begin{subequations}
	\begin{lemma}
		\label[lemma]{t:noise-decomp-linear}
		For the linear SA recursion,
		under {\em (A2)}:
		\whamrm{(i)}
		The disturbance decomposition   \eqref{e:DeltaDecomp} 
		with $\alpha_n=\alpha$ holds
		in which   
		\begin{equation}
			\label{e:LinDistDecomTerms}
			\begin{aligned}
				\MD_{n+1}   &=   \MDA_{n+1} \theta_n -\MDb_{n+1}   \,,   \qquad \clT_{n+1}  =   \uppsiA_{n+1} \theta_{n+1}  - \uppsib_{n+1}   b_{n+1},
				\\
				\Oops_{n+1}  &=   \uppsiA_{n+1}  [  A_{n+1} \theta_n  -  b_{n+1} ],
			\end{aligned} 
		\end{equation}
		where  $ \MDA_{n+1}  =  \Expect[ \haA_{n+1} \mid \clF_n]  - \haA_{n+1}  $,  $ \MDb_{n+1}  =  \Expect[ \hab_{n+1} \mid \clF_n]  - \hab_{n+1}  $; 
		$ \uppsiA_{n+1}    =  A_{n+1} - \haA_{n+1} $,  $ \uppsib_{n+1}  =   b_{n+1} - \hab_{n+1} $.
		
		\whamrm{(ii)}   We have for each $n$,
		\begin{equation}
			\Oops_{n+1}   =   \barH \theta_n  + \bargamma  +  \MD^H_{n+1} - \clT_{n+1}^H + \clT_{n}^H - \alpha\Oops_{n+1}^H ,
			\label{e:OopsDecomp}
		\end{equation}
		in which  the terms are defined in  the proof, found in \Cref{s:linearSAapp},   and
		$\{  \MD^H_{n+1} \}$ is a martingale difference sequence.
		
		\whamrm{(iii)}   The two covariance matrices \eqref{e:SigmaW} admit the approximation
		\[
		\Sigma_{\MD} =  \Sigma_{\MD^*}  +    
		\alpha  \ZMD         
		+ O(\alpha^2),
		\]
		in which  $ \ZMD     $ is  defined in  the proof, found in \Cref{s:linearSAapp}.		  
		\qed
	\end{lemma}

\end{subequations}


The mean flow for the bivariate SA recursion \eqref{s:SAlinearBivariate} is defined by the pair of ODEs,	  
\[
\begin{aligned}
	\ddt \odestate_t  &=    \derbarf^* (\odestate_t  - \theta^*) \,,
	\, \quad
	\ddt z_t     = \bigl[   M_\delta  +   \derbarf^* \bigr]    z_t.
\end{aligned} 
\]
Hence the Lipschitz conditions required in \Cref{t:sens} are satisfied.    
An approximation for the covariance  matrix $\SigmaPR$ appearing in
\eqref{e:PRCov_Markov}
requires additional notation:  
\begin{equation}
	\begin{aligned} 
		\MD^{H*}_{n+1}  & \eqdef \Oops^*_{n+1}  - \Expect[  \Oops^*_{n+1}   \mid \clF_n]  \,, \ \ \text{   with} \ \   \Oops^*_{n+1}  =  H_{n+1} \theta^* + \gamma_{n+1}   
		\\
		\barOops^*  & \eqdef    \Expect_\uppi[ \Oops^*_{n+1} ]  =   - \Expect_\uppi[\haA_{n+1} (A_{n} \theta^* - b_{n})], 
	\end{aligned} 
	\label{e:MDOops}
\end{equation}
and recall that     $\MD^*_{n+1}   \eqdef  \haf(\theta^* , \Phi_{n+1})  - \Expect[  \haf(\theta^*, \Phi_{n+1})   \mid \clF_{n}]    $. 

\begin{mytheorem}
	\label[mytheorem]{t:sensLinear}
	Suppose that   {\em (A1)-(A5)}  hold for the linear SA recursion 	\eqref{e:LinSAgen},
	and that $\delta_s>0$ in \eqref{e:LinSens2}
	is chosen so that $  M_\delta(\alpha) + \derbarf^*$ is Hurwitz for all $\alpha \in (0, \alpha_0]$,   with $\alpha_0$ defined in   \Cref{t:sens}.   
	Then, the bivariate SA recursion  \eqref{s:SAlinearBivariate} also satisfies  {\em (A1)-(A5)}.
	
	Hence the conclusions of   \Cref{t:sens} hold along with the following refinements: 
	there is $\bdd{t:sensLinear}<\infty$  such that   for $0<\alpha \le \alpha_0$:

	\whamrm{(i)}     With $\clV$   defined in 	\eqref{e:clV} we have that for $ \bdde{t:sensLinear} =2  \delta_s$,  and any  $z=\Psi_0$,
	\begin{equation}
		\Expect[  \|  \theta_n - \theta_n^\infty  \|^2 ]    \le    \bdd{t:sensLinear}  \exp( -n  \bdde{t:sensLinear} \alpha  )      \clV^{1/2}(\Psi_0) \,, \ \ n\ge 0.
		\label{e:L2LyapExp}
	\end{equation}

	\whamrm{(ii)} 	
	The limit \eqref{e:TargetBias} holds  with
	$\displaystyle
	\barOops_\alpha
	=   - \Expect_\upvarpi[\haA_{n+1} (A_{n} \theta_n - b_{n})]
	= \barOops^* + O(\alpha)
	$. The  bias approximation
	\eqref{e:Bias} also holds,   with 
	$\Zthbias  \eqdef 
	[\derbarf^*]^{-1} \barOops^*$.

	\whamrm{(iii)}   
	The limit \eqref{e:PRCov_Markov}  holds  in which   
	\begin{equation}
		\ZthetaPR =    \YthetaPR +  \YthetaPR^\transpose    \ \ \textit{with} \ \    \YthetaPR \eqdef G^*   
		\big(  \half  \ZMD  + \Expect[   \MD^{*}_{n+1}  (  \MD^{H*}_{n+1}  )^\transpose ]   +  \derbarf^* \SigmaTheta^* \barH^\transpose \big)    {G^*}^\transpose.
		\label{e:ZsensLInear}
	\end{equation}
\end{mytheorem}

Similarly to was obtained for nonlinear recursions in \Cref{t:corr_sens}, stronger conclusions are obtained under additional conditions on the noise of the algorithm. In particular, we obtain a generalization of the results in \cite{bacmou13,durmounausam24} for linear SA with $\bfPhi$ i.i.d. in which asymptotic unbiasedness of averaged estimates was established.
\begin{corollary}
	\label[corollary]{t:corr_sens_lin}
	Suppose that the assumptions of \Cref{t:sensLinear} hold. Suppose moreover that either of the following conditions hold:
	{\rm (a)} The noise is additive, as in \eqref{e:WarsawStochasticQSA}, or 
	{\rm (b)}  The process $\bfDelta$ is a martingale difference sequence. 
	
	Then, the conclusions of \Cref{t:sensLinear} hold. Moreover, the right hand side of \eqref{e:TargetBias} and \eqref{e:Bias} is zero.
	\qed
\end{corollary}

%

\section{Theory without convergence}
\label{s:nocon}
In the  following result Assumption (A5) of \Cref{t:sens} is relaxed, so in particular there may be many solutions to $\barf(\theta) =0$.   In particular, 
\Cref{t:optAssumptions}~(iii) shows that, under mild assumptions,  the mean of the empirical target bias \eqref{e:empTargetBias}
vanishes and that its covariance   is small.      In fact, we show in \Cref{s:opt} that we may achieve $\SigmaCLT^{\Probe} =0$ through design, so that the limit in  \eqref{e:NullTargetBiasMSE} is zero.    These strong results are available only for the additive noise model \eqref{e:WarsawStochasticQSA}.

\begin{mytheorem}
	\label[mytheorem]{t:optAssumptions}%
	\begin{subequations}%
		Consider the SA recursion \eqref{e:WarsawStochasticQSA} subject to the following conditions:   
		$\barf  $  is globally Lipschitz continuous, and 
		the ODE@$\infty$ is   asymptotically stable.  Assume  that $\Probe_n = F(\Phi_n)$ for each $n$, in which 
		$\bfPhi$ satisfies (A2ii) and $F\colon\state\to\Re^d$ is a measurable function satisfying $\uppi(F_i) =0$ for each $i$.  
		Moreover, 
		for  an everywhere positive function $p_\tProbe \colon\state\times \Re^d \to (0,\infty)$,  an integer $n_0\ge 1$,   $\bdd{t:optAssumptions}<\infty$,
		and $\bddepsy{t:optAssumptions} >0$,  
		\begin{align}
			\Expect[ \exp(  \bddepsy{t:optAssumptions} \| \Probe_{1} \|^2 ) &\mid \Phi_0=x]  \le     \exp(   W( x)        +   \bdd{t:optAssumptions}) 
			\label{e:ProbeAssumptions}
			\\
			\Prob\{ \Probe_{n_0} \in A & \mid \Phi_0=x \} =  \int_A  p_\tProbe (x, r) \, dr \,, \quad x\in\state\,,   A\in\clB(\Re^d).
			\label{e:ProbeAssumptionsDensity}
		\end{align}%
		\label{e:optAssumptionsPhi}%
	\end{subequations}%
	Then, there is $\alpha_0>0$   such that the following hold for each $\alpha\in (0,\alpha_0]$:\whamrm{(i)}  
	The joint process $\bfPsi$ is a $\psi$-irreducible and aperiodic Markov chain satisfying   {\em (DV3)}:
	\[
	\Expect\bigl[  \exp\bigr(  V_\circ(\Psi_{k+1})      \bigr) \mid \Psi_k= z = (\theta;x) \bigr]  
	\le  \exp\bigr(  V_\circ(z)  - W_\circ(z) +  b_\circ \ind_{C_\circ}(z)\bigl),
	\]
	where $  b_\circ$ is a finite constant,  $C_\circ  \subset C_\tTheta \times C_{\tState}$   is a small set for the joint process   with
	$C_\tTheta\subset\Re^d$ compact and  $ C_{\tState} \in\bx$ small for $\bfPhi$.    
	The  Lyapunov function is  of the form $V_\circ(\theta,x) = \half [\alpha^{-1} \nu(\theta) + V(x)]$,  in which $\nu \colon\Re^d\to \Re_+$ is continuously differentiable  with Lipschitz gradient,  
	and $W_\circ(\theta,x) =  1+ \epsy_\circ [  \| \theta \|^2 + W(x)]$     
	for some $\epsy_\circ>0$.   
	
	\whamrm{(ii)}   
	For any measurable function $g\colon\Re^d\to\Re^m$ with polynomial growth,
	\[
	\begin{aligned}
		& \lim_{n\to\infty} \Expect[  g( \theta_{n} )  \mid \Psi_0 =z ]  = \upvarpi(g)
		\\
		&   \frac{1}{ \sqrt{N} } S_N^g  \darrow  N(0,\Sigma_g) \,, \ \ N\to\infty\,,
	\end{aligned}
	\]   
	where $S_N = \sum_{n=1}^N [ g( \theta_{n} )  - \upvarpi(g) ]$,  $N\ge 1$, and $\Sigma_g\ge0$.

	\whamrm{(iii)}  
	The limits \eqref{e:TargetBiasResults}
	hold,   with $\SigmaCLT^{\Probe}$ denoting  the asymptotic covariance of $ \{\Probe_{n}      \}$.
	
\end{mytheorem}

\Proof     The proofs of parts  (i) and (ii)  are contained in  \Cref{s:optAssumptionsProof}.   
Part (iii) follows from parts (i) and (ii) combined with the following representation, obtained by 
averaging each side of \eqref{e:WarsawStochasticQSA}:   
\begin{equation*}
	\upbeta_N^\barftwo   
	=    \frac{1}{\alpha}\frac{1}{N} [\tiltheta_N -\tiltheta_0]  
	+ \frac{1}{\sqrt{N}}  D_N
	\,,\qquad D_N  =  \frac{1}{\sqrt{N}}  \sum_{n=1}^{N} \Probe_n  
\end{equation*}
The proof is completed on applying the definition $\displaystyle \SigmaCLT^{\Probe}
= \lim_{N\to\infty} \Expect[D_N D_N^\transpose]$  (see  \cite[Thm.~2.3]{glymey96a}).  
\qed
%
%

\subsection{Stochastic gradient descent}
\label{s:SGD}
Stochastic approximation is commonly employed to obtain an approximation of the minimizer $\theta^* \in \argmin \Obj(\theta)$ in which $\Obj \colon\Re^d\to\Re$ is known as the objective function.  This function together with a matrix-valued function $G \colon\Re^d\to\Re^{d\times d}$   will be assumed continuously differentiable.

The stochastic gradient descent (SGD) algorithm considered here is an example of  \eqref{e:WarsawStochasticQSA} 
in which $\barf(\theta) = - G(\theta)  \nabla \Obj\, (\theta)$ and the $d$-dimensional ``probing'' or ``exploration'' signal $\{  \Probe_{n+1} \}$ is a design choice.  Exploration is introduced to avoid local minima in applications for which the objective $\Obj$ is not convex.   
Polyak's higher-order techniques, including the famed introduction of momentum, can be cast as \eqref{e:WarsawStochasticQSA}  for some function     $\barf$
and a re-interpretation of $\theta$  
\cite{pol87,boyvan04a}.

We henceforth restrict  to the standard algorithm in two special cases:   Let $\{W_n \}$ denote an i.i.d.\  sequence on $\Re^d$ with zero mean,  
$\Expect[ \exp(\epsy_0 \|  W_n \|^2) ] <\infty$ for some $\epsy_0>0$, and consider:
\\
\begin{tabular}{lrll}
	\textbf{1.}  &  i.i.d.\ exploration: &  $\Probe_{n+1}= W_{n+1}$, & $\SigmaCLT^{\Probe} =  \Cov(W_n)$;
	\\
	\textbf{2.}  & zig-zag exploration: &  $\Probe_{n+1}= [ W_{n+1} - W_n ]/\sqrt{2}$, \phantom{1cm} & $\SigmaCLT^{\Probe} = 0,$   
\end{tabular} 
\\
where the normalization is imposed to ensure $\Expect_\uppi [  \|  \Probe_{n+1} \|^2 ]$ is identical in each case.

\begin{corollary}
	\label[corollary]{t:optConclusions}
	Suppose  the assumptions of \Cref{t:optAssumptions} hold. Then:
	
	\whamrm{(i)} 
	In case~1, $\Expect[ \| \upbeta_N^\barftwo \|^2 ]  = O(1/{N})$;
	
	\whamrm{(ii)}  
	In case~2, $ \Expect[ \| \upbeta_N^\barftwo \|^2 ]  = O(1/{N^2})$. If in addition the assumptions of  \Cref{t:sens} hold, we have that $\Ztheta=0$ in case~2,  so that \eqref{e:theta2th-moment} is improved to $   \Expect_\upvarpi[  \|    \tiltheta_0  \|^2]    = O(\alpha^2)$.  
	\qed
\end{corollary}


The following result provides sufficient conditions for  \eqref{e:optAssumptionsPhi}.   Denote $ \Lambda(\epsy)\eqdef \log(\Expect[ \exp(\epsy \|  W_n \|^2) ] )$  for $\epsy\in\Re$.
\begin{proposition}
	\label[proposition]{t:TowOptCases}    Suppose that  $\{W_n \}$ is an i.i.d.\  sequence on $\Re^d$ with zero mean,  and $ \Lambda(\epsy_0) <\infty$ for some $\epsy_0>0$.   Assume in addition that its marginal has a density $p_W$ that is continuous and everywhere positive on $\Re^d$.   Then,
	
	\whamrm{(i)} 
	In case~1 the  Markov chain  $\{\Phi_n = W_n \}$ satisfies {\em (A2ii)} and \eqref{e:optAssumptionsPhi} with $V(x) = \epsy_0 \|  x \|^2 $,   $W(x) =  1 +  V(x)$,   $b =\Lambda(\epsy_0)+1$,  
	and $C = \state $.

	\whamrm{(ii)}  
	In case~2 the  Markov chain  $\{\Phi_n = (W_{n-1};W_n) \}$ satisfies 
	{\em (A2ii)} and \eqref{e:optAssumptionsPhi} with 
	$V(w,w') = \epsy_0[ \half \| w \|^2 +\| w' \|^2 ] $,   
	$W(w,w') = 1+  \half \epsy_0[   \| w \|^2 +\| w' \|^2 ] $,   
	$b =\Lambda(\epsy_0)+1$,  
	and $C = \state$. 
	\qed
\end{proposition}
The proof of \Cref{t:TowOptCases} is postponed to   \Cref{s:optAssumptionsProof}.


\section{Examples of applications}
\label{s:ex}


\subsection{Optimization}
\label{s:opt}

Here, we examine the convergence of SGD estimates,
obtained through the recursion in~\eqref{e:WarsawStochasticQSA},
as discussed in detail in \Cref{s:SGD}.
In particular both i.i.d.\ and zig-zag exploration 
are considered, corresponding to cases~1 and~2,
respectively,
as described in \Cref{s:SGD}.

\begin{figure}[ht!]
	\centering
	\includegraphics[width=0.8\hsize]{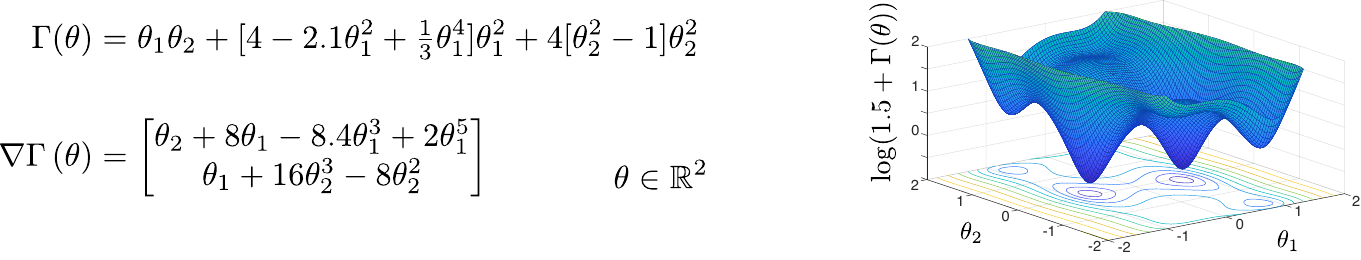}  
\caption{Multiple local minima for the six-hump camel-back function}
\label{f:camel}
\end{figure}

\wham{Multiple local minima.}  
First, we consider the 
objective function $\Obj(\theta)$ given 
by the \textit{six-hump camel-back function}. 
The function and its
gradient are shown in \Cref{f:camel}.
In fact, a plot of the function $\log(1.5 + \Obj(\theta))$ is shown 
on the right hand side, where the logarithm is introduced to better 
exhibit the multiple local minima.

\noindent
Choosing the diagonal matrix
$G(\theta) = \diag( 1/(1 + \theta_1^4)    ; 1/(1 + \theta_2^2) )$, we obtain,
for each $\theta\in\Re^2$,  with $\barf(\theta) = - G(\theta)  \nabla \Obj\, (\theta)$, that
\[
\begin{aligned}
\barf_\infty (\theta)    =  \lim_{c\to\infty}  \frac{1}{c} \barf(c\theta ) 
& =
-
\lim_{c\to\infty}  \frac{1}{c}  \begin{pmatrix} 
	2 (c\theta_1)^5 [1 + (c\theta_1)^4]^{-1}   
	\\
	16 (c\theta_2)^3   [1 + (c\theta_2)^2 ]^{-1}      \end{pmatrix}  =    -   \begin{pmatrix}   2  &0   \\  0  & 16  \end{pmatrix}  \theta  
\\
\lim_{c\to\infty}   \derbarf(c\theta )      &=     -   \begin{pmatrix}   2  &0   \\  0  & 16  \end{pmatrix}  \,, \quad \textit{where $\derbarf(\theta) =  \partial_\theta \barf(\theta)$.} 
\end{aligned}
\]
It follows that 
the ODE@$\infty$ is  asymptotically stable, and that   $\barf   $ is globally Lipschitz continuous.

\wham{Numerical results.} We examine the numerical estimates $\{\theta_n\}$
and $\{\Obj(\theta_n)\}$.
In the results surveyed here, 
the i.i.d.\ process used in cases 1 or 2 was taken Gaussian,
$W_n \sim N(0, \sigma_W^2 I)$.   Large $\sigma_W$ and/or large $\alpha$ is required for efficient exploration. 
Data was collected from $M = 10^5$ independent experiments,  with runlegth  $N=10^4$ in each experiment. We examine two scenarios.

\begin{figure}[ht!]
\centering
\includegraphics[width=0.72\hsize]{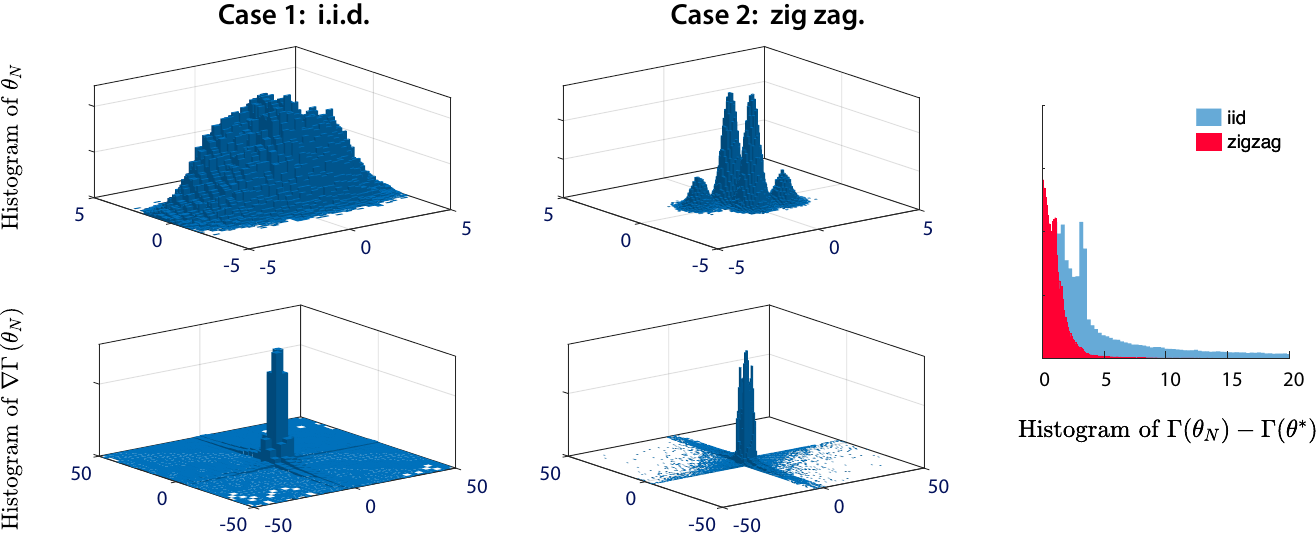}
\caption{SGD for the camel back function.   Large exploration and small gain:    $\sigma^2_W = 400$ and $\alpha=0.02$.   }%
\label{fig:N4M4ssp02sig400}
\end{figure}

\whamit{A. High variance small gain ($\sigma^2_W = 400$ and $\alpha=0.02$).}   \Cref{fig:N4M4ssp02sig400}   shows histograms of $\{ \theta_N^i  : 1\le i\le M \}$,
histograms of the respective gradients,   and a histogram of the observed error in objective $\{ \Obj(\theta_N^i ) - \Obj(\theta^*) : 1\le i\le M \}$.   
The relatively small variability observed in $\{  \nabla \Obj\, ( \theta_N^i ) : 1\le i\le M \}$ might be anticipated by \Cref{t:optAssumptions}~(iii).   
(In case~1, outliers were removed in the histogram for the sake
of visual clarity.)

\begin{figure}[ht!]
\centering
\includegraphics[width=0.72\hsize]{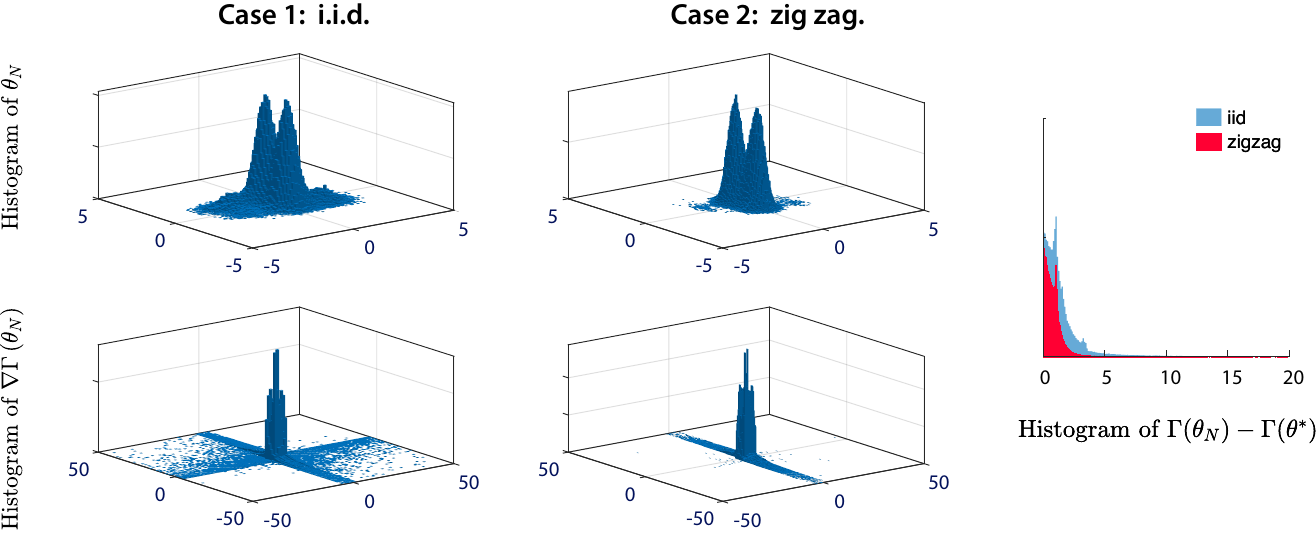}
\caption{SGD for the camel back function. Moderate exploration and large gain:   $\sigma^2_W = 10$ and $\alpha=0.1$.} %
\label{fig:N4M4ssp1sig10}
\end{figure}

\whamit{B. Moderate variance large gain ($\sigma^2_W = 10$ and $\alpha=0.1$).}    With the smaller variance and larger gain the performance in case~1 is greatly improved---see 
\Cref{fig:N4M4ssp1sig10}.  

\wham{An example with global convergence.}  
Suppose that the assumptions of \Cref{t:optAssumptions}
hold and, in addition, that the matrix $G(\theta)$ is invertible for 
each $\theta$  and that  the mean flow is globally asymptotically stable.  
The unique  equilibrium is necessarily equal to the optimizer,  $\theta^* = \argmin \Obj(\theta)$  since $ \nabla \Obj\, (\theta^*) =0$ and there cannot be any local minima or maxima.     This also implies that  $\derbarf^*  \eqdef \partial_\theta\barf\, (\theta^*)  = - G(\theta^*)  \nabla^2 \Obj\, (\theta^*)$, so that  (A5) holds provided  $G(\theta^*)  $ and $\nabla^2 \Obj\, (\theta^*)$ are positive definite matrices.  

An example of this situation is obtained by a modification of  
a well-known example,
the Styblinski-Tang function:
$\Obj(\theta) = \half  [ \theta_1^4 + \theta_2^4  + \xi (\theta_1 + \theta_2) - 16 \|\theta\|^2  ]$.  
This is  the Styblinski-Tang function when $\xi =5$,  
in which case there are five stationary points and one global optimizer
at $\theta_1^*=\theta_2^* \approx -2.9$.      On considering the partial derivatives 
$\partial_{\theta_i} \Obj\, (\theta) =   \half [4 \theta_i^3 +\xi   -32 \theta_i]$,  $i=1,2$, Lipschitz continuity of $\barf$ is obtained using 
$G(\theta) = \diag( 1/(1 + \theta_1^2)    ; 1/(1 + \theta_2^2) )$. Moreover,
with this choice,  $\barf_\infty(\theta) = - 2\theta$, for any value of $\xi$,
and all of the assumptions of  \Cref{t:sens} hold for sufficiently large $\xi$ so that (A5) holds.   

The numerical experiments that follow used  $\xi =50$ in which case there are no stationary points outside of the 
global optimizer at  $\theta_1^*=\theta_2^* \approx - 3.4$.           We have $\nabla^2 \Obj\, (\theta^*)  =  \sigma^2_{\lilObj} I$ with $\sigma^2_{\lilObj} = 
6( \theta^*_i)^2     - 16  >0$,   and $\derbarf^* = -   \sigma^2_{\lilA} I$ with $\sigma^2_{\lilA} = \sigma^2_{\lilObj}/[1+( \theta^*_i)^2  ] \approx 4$,
independent of $i$.
In case~1, the   minimal covariance \eqref{e:SigmaPRopt} is $ \SigmaTheta^* =     \Sigma_{\MD^*}  /\sigma^4_{\lilA} = (\sigma_W/\sigma^2_{\lilA})^2 I$,
and in case~2, $ \SigmaTheta^* =  0$.

\begin{figure}[ht!]
\centering
\includegraphics[width=0.72\hsize]{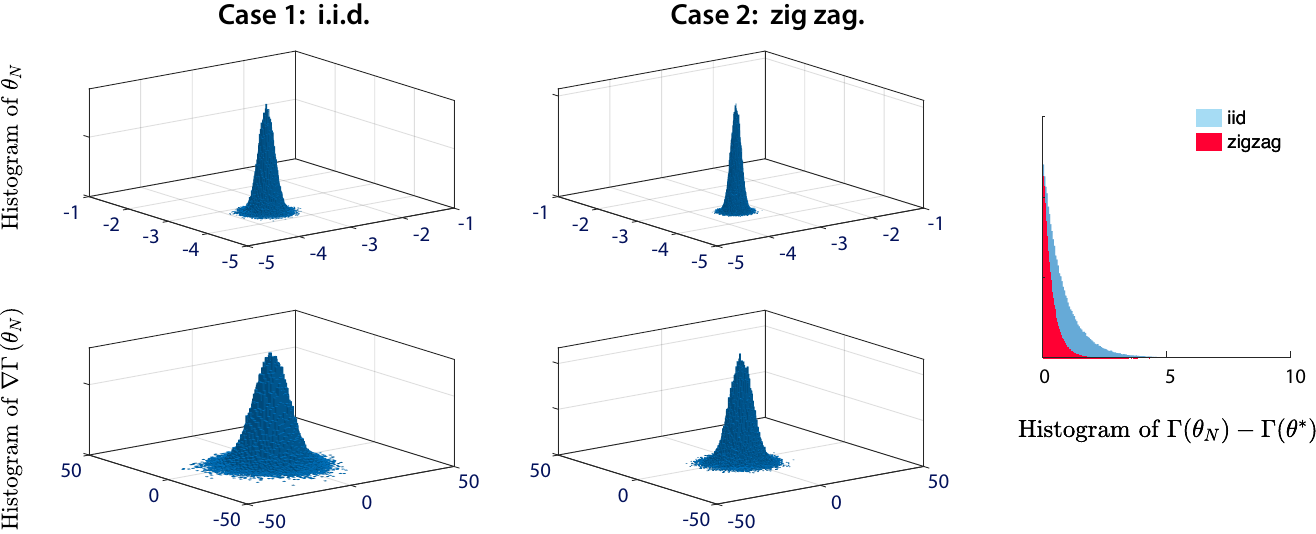}
\caption{SGD for the modified Styblinski-Tang function using moderate exploration and small gain.     }%
\label{fig:Styblinski-Tang}
\end{figure}

There is no reason for concern about exploration for this example.   With  
$\sigma_W =  1$ and $\alpha=0.1$,    performance was similar in each of the two cases:  
$\max_m   \Obj(\theta_N^m)  = \Obj(\theta^*) +  e_m$ with $e_m  = 12$ in case~1 and $e_m =   4$ in case~2,  and the average over $m$ of
$\{ \theta_N^m : 1\le m\le M \}$  was  nearly optimal in each case.     However, the histogram of  $\{ \Obj( \theta_N^m) - \Obj(\theta^*) : 1\le m\le M \}$  shown in	\Cref{fig:Styblinski-Tang} illustrates the benefit of using a probing signal with null asymptotic covariance.


%
%
%
%

\subsection{Reinforcement learning} 
\label{s:TD}

\wham{TD learning.}
In typical algorithms found in reinforcement learning,
the recursion cannot be expressed as the additive noise model \eqref{e:WarsawStochasticQSA}, 
which means that bias is a potential problem.   In particular, the limit \eqref{e:NullTargetBiasMSE} is not possible when $  \barOops_\alpha
\neq 0$  (recall	\eqref{e:TargetBias}).    In the following, the conclusions of  \Cref{t:sensLinear}  are   illustrated with a simple application of TD learning.

%
%
%

Suppose that $\bfmX$ is a Markov chain on a state space $\Spx $, assumed here to be Euclidean.    Given a cost function $c\colon\Spx\to\Re_+$ and 
discount factor $\gamma \in (0,1)$,  the associated value function is denoted  
\begin{equation}
J(x) = \sum_{k = 0}^\infty \gamma^k \Expect[c(X_k)| X_0 = x]   \,,   \quad x\in\Spx .
\label{e:valfunc}
\end{equation}   
The goal of TD-learning is to estimate this value function within a specified function class. 
For a linear function class $\{ J^\theta = \theta^\transpose\psi  :  \theta\in\Re^d \}$ with  basis $\psi: \Spx \to \Re^d$,
the   TD($\lambda$)-learning recursion is
\begin{subequations}
\begin{align}
\theta_{n+1}  & = \theta_n +\alpha_{n+1}     \Tdiff_{n+1} \elig_{n+1}
\label{e:TDnaive}
\\
\Tdiff_{n+1} & =  - J^{\theta_n}(X_n)  + c(X_n)   + \disc   J^{\theta_n} (X_{n+1}) 
\label{e:TDtheta}
\\
\elig_{n+1}  & = \lambda \disc \elig_n + \psi(X_{n+1})\,,  
\quad  
n\ge 0.
\label{e:eligTD}
\end{align}
Recursion  \eqref{e:eligTD} defines the ``eligibility vectors'' $\{\elig_n\}$, in which $\lambda\in [0,1]$.
The recursion  \eqref{e:TDnaive} is an instance of linear SA in which $\Phi_{n+1} = [X_n ; X_{n+1};  \elig_{n+1}]$.

\label{e:TDall}
\end{subequations}

In the example that follows we take $\Spx =\Re$, the cost function quadratic $c(x) = x^2$,   and $\bfmX$ evolves according to linear dynamics  $X_{n+1} = F X_{n} + W_{n+1}$  in which $  |F| <1$,  and $\{ W_n \}$ is an i.i.d.\  second-order process with zero mean.   Then the value function is also quadratic,
\begin{equation}
J(x) = \theta^*_1 + \theta^*_2 x^2 \, ,
\quad  
\theta^* \eqdef   [ 		\theta^*_1  ;  \theta^*_2 ]   \in \Re^2  \,,
\label{e:genform} 
\end{equation}
which motivates the  function class  with $d=2$ and 
$\psi(x) = (x^2;1)$.   In this special case the Markov chain $\Phi_{n+1} = [X_n ; X_{n+1};  \elig_{n+1}]$ satisfies (DV3),   for any choice of $\lambda$,  
in which $V$ and $W$ are quadratic functions on $\Re^3$, provided $\{W_n \}$ satisfies the assumptions of \Cref{t:TowOptCases}.

\wham{Numerical results.}  
Experiments were conducted in the linear Gaussian case,   $W_{n} \sim N(0,\sigma_W^2) $ with $\sigma_W^2=1$,  
$\gamma =0.7$,  $F=0.5$,  $\lambda =0$, and using the 2-dimensional ideal basis  $\psi(x) = (x^2;1)$.  
The exact value function is obtained with $\theta^* =  ( 1.2121...  ; 2.8282...)$.
The results surveyed here compare two  different choices of step-size:
$\alpha_n \equiv \alpha =  10^{-2}$ and $\alpha_{n} = \min\{\alpha, n^{-0.8}\}$.    The time horizon was taken to be  $N= 5\times 10^5$, and
$M=200$ independent runs  were conducted to obtain estimates of bias and variance.

\begin{figure}[ht!]
\centering
\includegraphics[width=0.85\hsize]{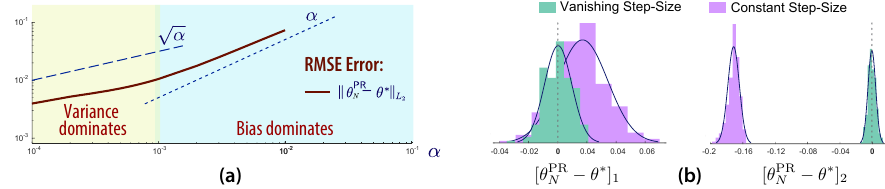}
\caption{Performance of TD($\lambda$) learning with $\lambda=0$.  (a) $L_2$ norm of estimation error for $\thetaPR_N$ with the fixed step-size algorithm as a function of $\alpha$. (b) Histogram of estimation error for vanishing and constant steps-size algorithms for each dimension of $\thetaPR_N$. }%
\label{fig:TD}
\end{figure}

\Cref{fig:TD}~(a) shows the $L_2$ norm of the estimation error obtained from PR averaged estimates as a function of $\alpha$. 
The plot is consistent with the heuristic  \eqref{e:CLTapprox}:
\whamit{Small step-size:} for $\alpha<10^{-3}$ the  $L_2$-error is approximated well by an affine function of $\sqrt{\alpha}$.   This would be anticipated by \eqref{e:CLTapprox}  and the theory in this paper, establishing that the bias $\| \thetaPR_\infty -\theta^*\|$ is of order $\alpha$.

\whamit{Large step-size:}   for $\alpha>10^{-3}$  the $L_2$-error is approximately  affine as a function $\alpha$ rather than $\sqrt{\alpha}$.
It appears that the   bias dominates the variance in this regime, 
at least for this $N$.   

Histograms for the estimation error corresponding to both choices of $\{ \alpha_n \}$ are shown in \Cref{fig:TD}~(b),   for each component of $\thetaPR_N$.   The impact of the bias formula 	\eqref{e:Bias} is evident in the constant gain algorithm.    In particular,     the histogram of error for the second parameter   is centered far from~$0$.  
Further experiments show that bias grows with   $\lambda \in [0,1]$.  This is unfortunate, since  the value $\lambda=1$ is required for unbiased estimates in the actor-critic method \cite{kontsi03a,CSRL}.         

These experiments used a very small discount factor ($\gamma = 0.7$) in order to obtain reliable estimates of bias and variance with moderate runlength.   
The next example is designed so that all of the statistics in   \Cref{t:sensLinear}   are easily computed.

\subsection{Impact of statistical memory}
\label{s:Seanexp}

We consider a scalar linear SA recursion,
\begin{equation}
\begin{aligned}
f(\theta_n, \Phi_{n+1})  &=  A_{n+1} \theta_n    - b  + \clW_{n+1}\\
\text{where }
A_{n+1} &= -1  + \clW_{n+1}  \,,
\quad
\clW_{n+1} =  \beta \clW_n  +  \sqrt{1- \beta^2}  \clW^\circ_{n+1},
\end{aligned}
\label{e:TDlike}
\end{equation}
and $\{\clW^\circ_n\}$ is  i.i.d.\ and Gaussian $N(0,1)$. It follows that  $\derbarf^*=-1$ and $\theta^*=-b$.   
The sequence $\{\clW_n\}$ resembles the eligibility vector appearing in the TD-algorithms of reinforcement learning---recall \eqref{e:eligTD}. 

The scaling by $ \sqrt{1- \beta^2} $ in \eqref{e:TDlike} is introduced to ensure that the steady-state variance of $\clW_n$ is unity,  but the asymptotic variance is large when $\beta\sim 1$:
$\SigmaCLT^{\clW}  = (1+\beta)/(1-\beta)$, from
\[\begin{aligned}
\SigmaCLT^{\clW}  = \sum_{n=-\infty}^\infty \Expect[ \clW_0  \clW_n]  &=  -\Expect[ \clW_0^2]+   2  \sum_{n=0}^\infty \Expect[ \clW_0  \clW_n] .
\end{aligned}
\]

A proof of the following conclusions can be found in \Cref{s:NumProofs}.

\begin{subequations}
\begin{proposition}
\label[proposition]{t:LinExample}
Consider the linear SA recursion in which $\{\clW_n\}$  evolves according to the linear recursion \eqref{e:TDlike},
with   $0\le \beta<1$  and  $\bfmN$ a standard i.i.d.\ Gaussian sequence,  $\clW^\circ_n\sim N(0,1)$ for each $n$.  Then:
\wham{(i)} Provided an invariant probability measure $\varpi \sim (\theta_n,\Phi_n)$ exists, the bias is
\begin{align}
	\Expect_\varpi [ \theta_0 ]  &= \theta^* 
	- \alpha \Expect_\varpi [ \Oops_{2}  ] 
	\label{e:BiasLinExample}
	\\
	\text{with}\quad
	\Expect_\varpi \bigl[ \Oops_2  ]  &=   -  \frac{\beta}{1-\beta}    [ 1+ \theta^*  ]   \\
	& -
	\frac{1}{1-\beta} 
	\Expect_\varpi \bigl[  \clW_{2}  \bigl[   ( {\clW}_{1} -1) [\theta_0-\theta^*  ]   \bigr] .
	\label{e:UpsSSlinExample}
\end{align}

\wham{(ii)}  The optimal asymptotic covariance \eqref{e:SigmaPRopt} is the scalar
\begin{equation}
	\SigmaTheta^* =  [1+\theta^*]^2  \SigmaCLT^{\clW}.
	\label{e:PRlin}
\end{equation}
\end{proposition}

\label{e:AllLinearStats2}
\end{subequations}

The SA recursion \eqref{e:TDlike} was implemented for  $\beta =0.9$, $\theta^\ocp =10$ and several choices of step-size: $\alpha_n \equiv \alpha$ for constant step-size and $\alpha_n = \tfrac{1}{2} n^{-\rho}$ for vanishing step-size. 
Five values of $\alpha$ were tested for the fixed step-size algorithm,  and five values of $\rho$ for the vanishing step-size case:
\begin{equation} 
\begin{gathered}
\alpha   \in 
\{
5\times 10^{-4} ,  2.8\times 10^{-3}, 1.58 \times 10^{-2} , 8.89 \times 10^{-2}  , 0.5
\}
\\
\rho   \in 
\{ 
0.4000  , 0.5375  , 0.6750   , 0.8125  ,  0.9  
\}
\end{gathered}
\label{eq:settings}
\end{equation}

The estimates for the fixed step-size algorithm remained bounded in $n$ for the range of $\alpha$ tested. 
The same was true for the vanishing step-size algorithm, as predicted by theory \cite{borchedevkonmey25}. 
In application of PR-averaging    \eqref{e:PRdefSA}, the value  $N_0=0.2 N$ was chosen in all ten cases.    
With the given numerical values,  applying \eqref{e:PRlin}
gives the approximation for the vanishing gain algorithm,
$
(N-N_0)\Expect[ (\thetaPR_N -\theta^*)^2]   \approx \SigmaTheta^*  \approx  2.3\times 10^3$.

\begin{figure}[ht!]
\centering
\includegraphics[width=0.6\hsize]{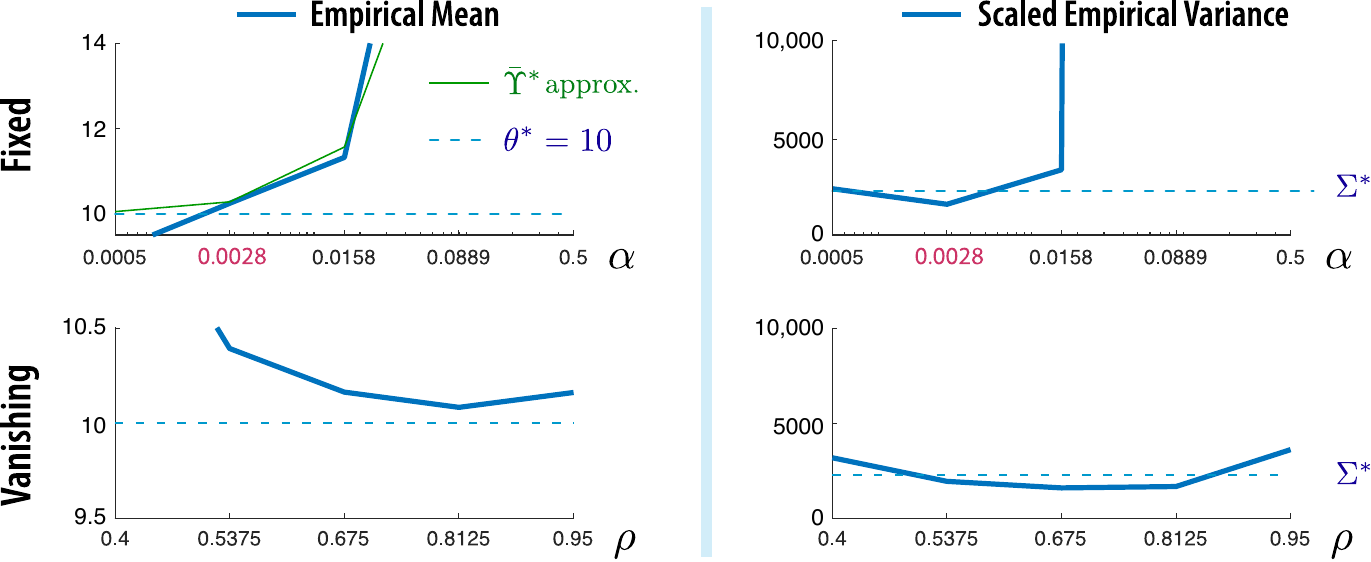}
\caption{Comparison of empirical bias and variance obtained from PR-averaging,
as functions of $\alpha$, for the recursion \eqref{e:TDlike}.}
\label{fig:MeanVarPlots}
\end{figure}

\Cref{fig:MeanVarPlots} shows the estimates of mean and variance obtained in each case.    The plot does not reveal much information for the fixed step-size algorithms because most values of $\alpha$ gave poor results.   The singular winner over all fixed step-size gains was  $\alpha^\star = 2.8\times 10^{-3} $,  resulting in  $ \Expect_\varpi [ \theta_0 ] \approx 
10.29$  and $(N-N_0)\Cov_{\varpi} [ \theta_0 ] \approx 0.7*\SigmaTheta^*$.  
The other four performed far worse.

Each of the experiments using a vanishing gain resulted in  variance of approximately equal to what was obtained using $\alpha^\star$ and with smaller  bias.     
Large bias can be anticipated for the fixed step-size algorithms by consideration of  \eqref{e:UpsSSlinExample}.   
For small $\alpha $ we obtain an approximation by ignoring the second term in this expression:
\[
\Expect_{\varpi} [ \theta_0 ]  = \theta^* 
- \alpha
\Expect_{\varpi} \bigl[ \Oops_2  ] \approx  \theta^*+   \frac{\beta}{1-\beta}    [ 1+ \theta^*  ]  \alpha   =   \theta^*+ 99\alpha.
\]
For $\alpha^\star $, we have $  \theta^*+ 99\alpha^\star  \approx 10.28$, so this approximation nearly matches the approximation  $ \Expect_{\varpi} [ \theta_0 ] \approx 
10.29$ obtained through simulation.
See the plot on the upper left hand side of \Cref{fig:MeanVarPlots} for a comparison of this approximation with the empirical mean.     For  the smallest value of $\alpha$ tested, the parameter estimates are far from steady-state by the end of the run. 
In this case we typically observe \textit{negative} bias.
The cause of the negative bias for   $\alpha = 5\times 10^{-4}$ is explained by the fact that $\theta_0^i$ is drawn from $N(0,25)$,
which has zero mean,  while $\theta^*=10$. 

\begin{figure}[ht!]
\centering
\includegraphics[width= 0.6\hsize]{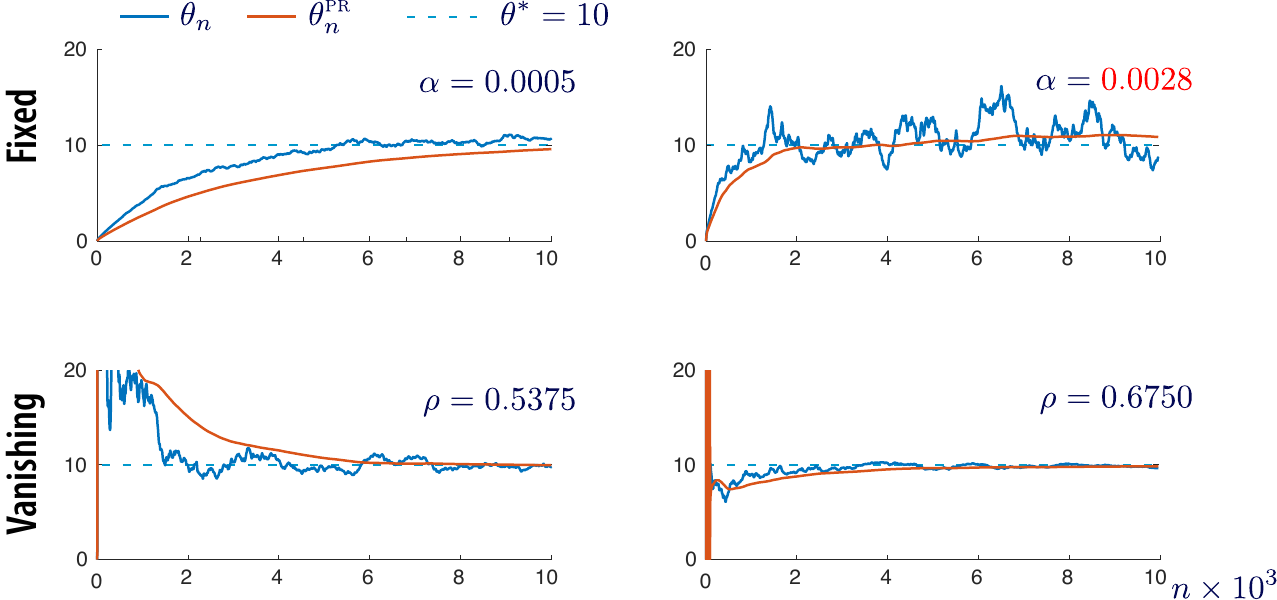}  
\caption{Evolution of estimates with and without PR averaging from four experiments.}
\label{fig:FourPlotsLin} 
\end{figure}

\Cref{fig:FourPlotsLin} shows sample paths with and without averaging for two selected values of fixed step-size, and two values of $\rho$ for   vanishing step-size,  with initialization $\theta_0=0$ in each case. The plots using PR-averaging were obtained via \eqref{e:PRdefSA} with $N_0=\lfloor 0.8 N\rfloor$ for each $N$. 
It is clear that
$\alpha = 5\times 10^{-4}$ fails,  
and  $\alpha = 2.8\times 10^{-3}$ performs much better.

\begin{figure}[ht!]
\centering
\includegraphics[width= .6\hsize]{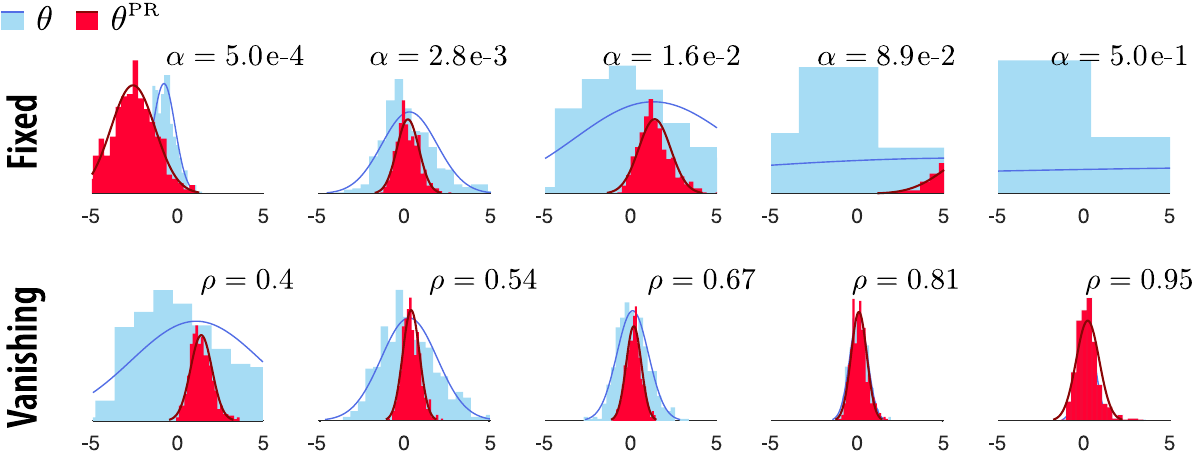}  
\caption{Histograms of the error $\tiltheta_N $ and  $\tilthetaPR_N $ obtained from $M=500$ independent runs in each of 10 experiments.   The top row shows results using the fixed step-size algorithm for various values of $\alpha$, and the bottom row shows results with vanishing step-size.
}
\label{fig:MeanVarHists}
\end{figure}

The ten subplots   in \Cref{fig:MeanVarHists} show histograms of $\{ \theta_N^i  \,,   {\thetaPR_N}^i  :   1\le i\le M \}$ for each of the ten settings
for $\alpha$ and $\rho$ as in \Cref{eq:settings}.   
The results using a vanishing step-size are not sensitive to $\rho$,   
even though the assumptions of the theory 
(which require $\half < \rho  <1$)
are violated for the smallest value  
$ \rho   =     0.4$.

\section{Conclusions and further work}
\label{s:Concl}

This work demonstrates that constant step-size algorithms for learning are appropriate in certain settings,  but may be prone to significant bias in others.   The source of bias is the combination of multiplicative noise and statistical memory,  which is typically unavoidable in applications to reinforcement learning.

There are many open paths for research:

\whamb 
Recall the asymptotic target bias $\barOops_\alpha$ appearing in \eqref{e:TargetBias}.
In prior work, methods were developed to ensure $\barOops_\alpha = \Zero$    for  \textit{quasi-stochastic approximation}  \cite{laumey22e,laumey22b}, in which   $\bfPhi$ is deterministic, and subject to special restrictions.  Extensions of this approach appeared recently in \cite{huochexie22,allgas24} based on variants of Richardson-Romberg extrapolation.     Little is known about the impact on variance.

\whamb   There are open questions in the domain of SGD:

$\circ$   \Cref{t:optAssumptions} tells us that the mean squared empirical target bias $\{ \upbeta_N^\barftwo \}$  decays at rate $O(1/N)$ for i.i.d.\ exploration,  and this can be improved to $O(1/N^2)$ for zig-zag.  With deterministic exploration, theory in  \cite{laumey22e,laumey22b} might be applied to establish the convergence rate $O(1/N^m)$ for any fixed $m<4$.  \textit{What is, then, the fundamental limit for algorithm performance? }

$\circ$    If the mean flow associated with \eqref{e:WarsawStochasticQSA} is globally asymptotically stable and the other assumptions of \Cref{t:sens} hold,   then the convergence rate of  $\{ \thetaPR_n \}$ in mean square will be of order $O(1/N^2)$ for zig-zag exploration.   However, this is only if the gain $\alpha$ is chosen sufficiently small.  The exact same convergence rates would hold for a carefully designed vanishing gain algorithm,    and there is little concern about gain selection since  stability of the SA recursion is automatic.   For the example considered in \Cref{t:optAssumptions}, the choice $\alpha_n = 1/n^\rho$  results in the  optimal  rate of convergence through PR averaging for any choice of $\rho \in (0,1)$
\cite{laumey24a}.
\textit{What is, then, the motivation for a constant gain algorithm?}

$\circ$   If the objective has many local minima, as in the example shown in \Cref{f:camel},  then there is ample motivation for a non-vanishing gain algorithm.   However,  in cases for which $\barf$ has multiple roots, it makes no sense to employ PR averaging,   and the performance  metrics emphasized in this paper are not always reliable predictors of success.     It is our conjecture that the techniques in this paper can be used to estimate a global minimum much more quickly than using  methods such as simulated annealing.

%
%
%
%
%
%
%
%
%
%
%
%

\bibliographystyle{abbrv}


\def\cprime{$'$}\def\cprime{$'$}

\clearpage
\appendix

\section{Appendices}



\subsection{Mean flow stability theory}

The mean flow is exponentially asymptotically stable (EAS)  if  for positive constants $\bddMF$ and $\rhoMF$,  and any   $\odestate_0\in\Re^d$,
\begin{equation}
\| \tilodestate_t \|   \le  \bddMF \| \tilodestate_0 \|    \exp(-\rhoMF t)  \,,\qquad t\ge 0,
\label{e:EAS}
\end{equation}
where $\tilodestate_t = \odestate_t -\theta^*$.       
\Cref{t:qsaprobeExpAS} below tells us that EAS holds under (A5).  
Part (i) of \Cref{t:qsaprobeExpAS} is elementary,  and the bound \eqref{e:EAS}
in (ii)    is taken from
\cite[Prop.~A.11]{laumey25a}.    The bound \eqref{e:EASsens}
follows from  \eqref{e:EAS}  and the Hurwitz assumption on $\derbarf^*$.

\begin{proposition}
\label[proposition]{t:qsaprobeExpAS}
The following hold for the mean flow 	\eqref{e:meanflow} under \AthreeV:

\whamrm{(i)}   If $\barf$ is globally Lipschitz continuous then there is    $ \varrho_e>0$ and $ b_e<\infty$ such that for any  initial condition satisfying $\| \odestate_0\| \ge 1/ \varrho_e$,   
\begin{equation}
\| \odestate_t \| \le    b_e \| \odestate_0\| e^{- \varrho_e t} \,,\qquad 
0\le t\le  T(\odestate_0)
\eqdef  \min \{  t\ge 0 : \| \odestate_t \|  \le 1/ \varrho_e \} \, .
\label{e:ExpStableAlmost}
\end{equation}		

\whamrm{(ii)} 
Suppose in addition that (A5) holds.  Then the mean flow is EAS.  
Moreover, the sensitivity process \eqref{e:sensODE}
satisfies for   positive constants $\bddMF$, $\rhoMF$,  and $\kappa_0\ge 1$,
\begin{equation}
\| \sens_t \|   \le  \bddMF \| \tilodestate_0 \|^{\kappa_0}   \exp(-\rhoMF t)  \,,\qquad t\ge 0.
\label{e:EASsens}
\end{equation}
\end{proposition}

\subsection{Markov chain bounds}
\label{s:mc}

We collect together general results for the Markov chain $\bfPhi$.   

Under (A2),
$\bfPhi$   is  $v$-uniformly ergodic with Lyapunov function $v=e^V$, from which we obtain the following  results   \cite{MT}.   The matrix
$\SigmaCLT^g$ is known as the asymptotic variance of $\{g(\Phi_k) \}$ \cite[Sec.~17.5]{MT}. 

\begin{proposition}
\label[proposition]{t:Vuni}
Suppose that {\em (A2ii)} holds.   Then,  there exists $b_v< \infty$ and $\varrho_v\in (0,1)$  such that 
for any measurable function   $g: \Re^d \times \state \to \Re^d$ and any $x\in\state$:

\whamrm{(i)}  If there is $b_g <\infty$ such that
$\|g(x)\| \le  b_g v(x)$ for each $x$,   then
\[
\|  \Expect [ \tilg(\Phi_n) \mid \Phi_0 =x]   \|  \le b_v  b_g  v(x)   \varrho_v^n \,, \quad n\ge 1\,, 
\]
where $\tilg(x) = g(x) - \uppi(g)$ for all $x$.

\whamrm{(i)} 
If the bound is strengthened to $\|g(x)\|^2 \le b_g^2 v(x)$ for each $x$ we have,  as $N\to\infty$, 
\begin{equation}
\begin{aligned}
	\frac{1}{N}    \Cov\Big(\sum_{i=1}^N  g(\Phi_i) \Big)   \to  \SigmaCLT^{g}    
	\quad
	\textit{and}
	\quad
	\frac{1}{\sqrt{N} }   \sum_{i=1}^N  \tilg(\Phi_i) \darrow  N(0, \SigmaCLT^{g}   )   
\end{aligned}
\label{e:hollandCLT}
\end{equation}
where    the second limit is in distribution.
\qed
\end{proposition}

The next result establishes that  the function $V_p (x) =   [ 2 + V(x)]^{p+1}  $ satisfies the drift condition,
\begin{equation}
\Expect\bigl[     V_p(\Phi_{k+1}) -  V_p(\Phi_{k})     \mid \Phi_k=x \bigr]  
\le     -  \delta_p  W(x) [ 1+V(x) ]^p  +   b_p^v \ind_{C_p}(x)  \,, \ \ x\in\state\,, 
\label{e:Vp3}
\end{equation}
for finite constants $\delta_p, b_p^v$,  and a  set $C_p$ that is small.

\begin{proposition}
\label[proposition]{t:bounds-Hgen}	
Suppose that {\em (A2)} holds.  Then, for each integer $p\ge 1$  the bound \eqref{e:Vp3} holds with  $\delta_p>0$, $ b_p^v<\infty$, and   $C_p$  a small set.   Moreover, there exists  $\bdd{t:bounds-Hgen}_p<\infty$ such that the following hold:
If $g\colon \state\to\Re^d$ is a measurable function satisfying  
$ \| g( x)   \|  \leq  b_g  W(x) [1+V(x) ]^p$   for some   $  b_g <\infty$ and all     $x\in\state $,
then, there exists $\hag \colon   \state\to\Re^d$   solving 	Poisson's equation,
$\Expect[\hag(\Phi_{n+1}) - \hag(\Phi_{n}) \mid \Phi_n = x] 
= 
- g( x) +  \uppi(g)$.
Moreover, it can be chosen so that $ \uppi(\hag)  \eqdef  \Expect_\uppi[\hag( \Phi_n)] = 0$, and  satisfying the bound
\[
\|  \hag(x) \|  \le   b_g  \bdd{t:bounds-Hgen}_p   [1+V(x) ]^{p+1}  \,, \quad x\in\state \, .
\]
\end{proposition}

\Proof   
The second derivative test shows that the function 
$G_p(x) = [ 2 + \log x]^{p+1}$  is concave on $[1,\infty)$.    
We also have  $V_p (x) \eqdef  [ 2 + V(x)]^{p+1} = G_p(v(x))   $ for $x\in\state$.    Hence,
by Jensen's inequality and (DV3),
\[
\Expect\bigl[     V_p(\Phi_{k+1})     \mid \Phi_k=x \bigr]  
\le   \bigr( 2 +  V(x)  - W(x) +  b \ind_C(x)  \bigl)^p .
\]
And under (A2ii) we can find $\delta_p>0$,  $b_p^v<\infty$ such that \eqref{e:Vp3} holds,  and for $r\ge 1$   $C_p =  S_W(r)  := \{ x :  W(x)\le r \}$ is   small under (A2).     This is condition (V3) of \cite{MT},
and the desired conclusions for $\hag$ follow from \cite[Thm.~2.3]{glymey96a}.
\qed

\smallskip

\Cref{t:DV3multBdd} represents a first step in the proof of  the coupling result \eqref{e:Coupling} in \Cref{t:sens}.   
The bound  \eqref{e:DV3multBdd+}  appears as Prop.~3 
in \cite{borchedevkonmey25}.  
We then obtain \eqref{e:DV3multBdd+L} as a corollary when $L$ satisfies (A4).   

\begin{subequations}
\begin{proposition}
\label[proposition]{t:DV3multBdd}
The following holds under {\em (A2)}  and  {\em (A4)}:  there exists a non-decreasing function $\bdd{t:DV3multBdd} \colon \Re_+\to\Re_+$ such that for any  $N\ge1$ and any non-negative sequence $\{\delta_k : 0\le k\le N-1 \}$,   on denoting $r = \sum \delta_k$,    
\begin{align} 
	\Expect_x\Bigl[\exp\Bigl( V(\Phi_N) +\sum_{k=0}^{N-1}   \delta_k W(\Phi_k)  \Bigr)  \Bigr]
	\le  \bdd{t:DV3multBdd} (r) v(x) \, ,  \quad &&    \text{provided} \    0<r\le 1;
	\label{e:DV3multBdd+}
	\\
	\Expect_x\Bigl[\exp\Bigl( V(\Phi_N) +\sum_{k=0}^{N-1}   \delta_k L(\Phi_k)  \Bigr)  \Bigr]
	\le  \bdd{t:DV3multBdd} (r)   v(x) \, ,  \quad &&    \text{for all} \   r>0.
	\label{e:DV3multBdd+L}
\end{align}  
\end{proposition}
\end{subequations}

A version of the following was first established in \cite{fostwe98} for a geometrically ergodic chain.   While the bound \eqref{e:couplePhi}
was only established with $v_+$ a bounded function on $\state$, the extension below follows from $v$-geometric regularity of $\bfPhi$ (see \cite[Sec.17.4]{MT}),
and the fact that $v_+$ is bounded by a constant times $v$ under (DV3) and the definition \eqref{e:clV}.

\begin{proposition}
\label[proposition]{t:PhiCouples}   
If {\em (DV3)} holds then for each $x \in \state$,  there is a joint realization of the Markov chain $\bfPhi$, with $\Phi_0=x$,  and
a stationary version of the Markov chain $\{ \Phi^\infty_k : -\infty < k <\infty \}$.      The  trajectories couple:  	there is a stopping time $\taucpl$ for the joint process such that  $ \Phi_k = \Phi_k^\infty$ for $ k\ge \taucpl$.  Moreover,   there are constants  $\bcpl$, $ \rhocpl<1$ independent of $x$, such that   for $k\ge 1$,
\begin{equation}
\Expect[ v_+(\Phi_k) \ind \{ k < \taucpl \} ]  \le   \bcpl \rhocpl^k  v(x).
\label{e:couplePhi}
\end{equation}
\end{proposition}

\subsection{Moment bounds and bias}
\label{s:momentbdds+bias}

We require the following moment bounds  as a first step in the proofs of \Cref{t:sens,t:optAssumptions,t:sensLinear}:

\begin{mytheorem}
\label[mytheorem]{t:BigBounds}
Suppose that   {\em (A1)-(A4)}  hold.      
Then there exists    $\bdd{t:BigBounds}<\infty$     and  $\alphaTmp{t:BigBounds}>0$ such that for $0<\alpha \le \alphaTmp{t:BigBounds}$:

\whamrm{(i)}
$\displaystyle
\sup_{k  , z}  \frac{1}{\clV(z)} \Expect\bigl[ \clV(\Psi_k) \mid \Psi_0 = z \bigr] \le \bdd{t:BigBounds} . 
$

\whamrm{(ii)}  
If in addition {\em (A5)} holds then,
$\displaystyle
\limsup_{n\to\infty}   \Expect\bigl[  \| \tiltheta_n\|^4  \bigr] \le \bdd{t:BigBounds}  \alpha^2$,
for any $ \Psi_0 \in\Re^d\times\state$.
\end{mytheorem}

The first step in the proof of \Cref{t:BigBounds} is to define ``sampling times'' in an ODE analysis: for each $k$, let $\SAtime_k = \alpha k$. For a given $T>0$, let  $T_0=0$ and $T_{n+1} = \min \{\SAtime_k: \SAtime_k \geq T_n + T\}$. Consequently, the sequence $\{T_n\}$ satisfies $\displaystyle T\leq T_{n+1} - T_{n} \leq T+  {\alpha} $ for each $n$. Let $m_0\eqdef 0$ and $m_n$ the integer satisfying $\SAtime_{m_n} = T_n$ for each $n\geq1 $.   Whenever Assumption (A5) is imposed, so that the
mean flow is EAS, it  will be assumed   that the value of $T>0$   is chosen sufficiently large so that 
\begin{equation}
\| \tilodestate_t\|   \le  \half  \| \tilodestate_0\|   \ \  \textit{ for       $\odestate_0\in\Re^d$,  \quad and}   \quad \| e^{\derbarf^*t} \|_F \le \tfrac{1}{4}  \qquad \textit{for all  $t\ge T$}.
\label{e:odeContractionT}
\end{equation}

Let  $\{ \odestate^{(n)}_{t}: t \geq T_n\} $ denote the solution to the mean flow \eqref{e:meanflow},   initialized  with $\odestate^{(n)}_{T_n} = \theta_{m_n}$. 
In the following this is compared with the \textit{interpolated parameter process}:
\begin{equation}
\text{$ \ODEstate_t  = \theta_k$   when $t=\SAtime_k$, for each $k\ge 0$,    }
\label{e:ODEstate}
\end{equation}
and defined for all $t$ through piecewise linear interpolation.

While \cite{borchedevkonmey25} concerns SA with vanishing step-size,  
no assumptions 
preclude $\{\alpha_k\}$ from being constant on any of the finite intervals 
$\{ k : m_n\le k < m_{n+1 }\}$.  
Consequently,  any of the finite-interval bounds in  \cite{borchedevkonmey25} are valid here,  on choosing $\alpha_k \equiv \alpha$ on any such interval.  
The following is one example.  		Recall the function $\clV$ is defined just above \eqref{e:clV}.

\begin{subequations}

\begin{proposition}
\label[proposition]{t:results_CLT}
Under  {\em (A1)-(A4)}, there is $\alpha_0>0$,  $\bdde{t:results_CLT}<1$,  and $\bdd{t:results_CLT}<\infty$,  			
such that the following bound holds for $0<\alpha\le \alpha_0$,   and any initial condition $\Psi_0$: 
\begin{align}
	\Expect[ \clV(\Psi_{k+ j })   \mid \clF_{k  }   ]  & \le  \bdd{t:results_CLT}  \clV(\Psi_{k})      \, ,   && k\ge  0 \,, \  j\ge 1
	\label{e:PreLcontraction}
	\\
	\Expect[ \clV(\Psi_{k+m_1 })   \mid \clF_{k  }   ]   &\le  \bdde{t:results_CLT}   \clV(\Psi_{k})    + \bdd{t:results_CLT}  \, ,    && k\ge  0.
	\label{e:Lcontraction}
\end{align}
\end{proposition}

\end{subequations}

\Proof
The bound \eqref{e:PreLcontraction} is  \cite[Prop.~12]{borchedevkonmey25},  and   \eqref{e:Lcontraction}
is established in the proof of   \cite[Thm.~2]{borchedevkonmey25}.
\qed

Denote
\begin{equation}
\begin{aligned}
y_k^{(n)} \eqdef  \theta_k - \odestate_{\SAtime_k}^{(n)}     \,,
\qquad    z_k^{(n)} \eqdef      \frac{1}{\sqrt{\alpha}} y_k^{(n)}   \,,    \ \ \qquad k\ge m_n. 
\end{aligned}
\label{e:scaled-error-def}
\end{equation}
Lemma 13 of  \cite{borchedevkonmey25} establishes the approximation,
\begin{align}
z_{k+1}^{(n)}  
& =   z_k^{(n)} + \alpha    \derbarf_k  z_k^{(n) }  +   \sqrt{\alpha} [  \Delta_{k+1} + \clE_k^T + \clE_k^D  ]  \,, 
\label{e:SA-z}
\end{align}
in which   $\derbarf_k \eqdef \derbarf(\odestate^{(n)}_{\SAtime_k} )$,   
$\clE_k^T  \eqdef  \barf(\theta_k)-
\barf(\odestate^{(n)}_{\SAtime_k})  -A  (\odestate^{(n)}_{\SAtime_k})   (\theta_k -\odestate^{(n)}_{\SAtime_k}    ) $,
and
$ \clE_k^D$  denotes  the error in replacing   the integral $\int_{\SAtime_k}^{\SAtime_{k+1}} \barf(\odestate^{(n)}_t)\, dt$  by
$\alpha    \barf(\odestate^{(n)}_{\SAtime_k})  $.   This is used in the proof of the following.

\begin{proposition}
\label[proposition]{t:tight}
The following holds under    {\em (A1)-(A5)}:  There exists $\alphaTmp{t:tight} >0$ such that, 
for any $T>0$,
there exists a deterministic constant $\bdd{t:tight}<\infty$ independent of $0<\alpha\le \alpha_0$ such that
\[ 
\Expect\bigl[\|z_k^{(n)}\|^4  \mid \clF_{m_n} \bigr]  \leq \bdd{t:tight} \clV(\Psi_{m_n})  \,,  \quad   m_n < k   \le  m_{n+1} \, ,   \  n\ge 0\, .
\]
\end{proposition}

\Proof
This bound is similar to those in  \cite[Prop.~15]{borchedevkonmey25}. 
However, in that prior result the bounds involved a full expectation.   We show here how the proof can be adapted to obtain the desired result.    

Let $b^{(n)}_k =     \Expect[ \|z_k^{(n)} \|^4 \mid \clF_{m_n}] ^{1/4}$.   We show that 
for a constant $\bdd{e:15BG}$ depending only on $T$,     
\begin{equation}
b^{(n)}_k \le   \bdd{e:15BG} \Big[ \clV(\Psi_{m_n})    +   \alpha    \sum_{j=m_n+1}^{k} b^{(n)}_j   \Big]\,, \qquad  m_n< k\le m_{n+1}\, .
\label{e:15BG}
\end{equation}
This combined with the discrete discrete Bellman-Gronwall lemma implies the desired result with $\bdd{t:tight} =  \bdd{e:15BG} \exp( \bdd{e:15BG} T)$.

The proof of \eqref{e:15BG} begins  by expanding the representation \eqref{e:SA-z}:  for $m_n\le l < k \le m_{n+1}$,  
\begin{equation}
\label{e:zn-acc-sum}
\begin{aligned}
z_k^{(n)} = \ScaledST_{l, k-1}  z_{l}^{(n)}  +  \sqrt{\alpha} \sum_{j=l}^{k-1} & \ScaledST_{j+1, k-1} [\clE_j^T  + \clE_j^D +   \Delta_{j+1}], 
\\
\text{where } \ScaledST_{l, k-1}  \eqdef &
\begin{cases}
	\prod_{i=l}^{k-1}[I+\alpha \derbarf(\odestate^{(n)}_{\SAtime_{i}})]  & l < k  \\
	I, &   l=k \, .
\end{cases} 
\end{aligned}
\end{equation} 
Using the fact that $\| \derbarf(\varble)\|_F$ is uniformly bounded over $\Re^d$,  it follows that there is $r_0>0$ such that   
\[
\|\ScaledST_{l, k-1} -I\| \leq e^{r_0 T}  (\SAtime_k - \SAtime_{l})    \,, \qquad m_n\le l\le k \le m_{n+1}.
\]

\begin{subequations}

Bounds on the fourth moment of each of the terms in \eqref{e:zn-acc-sum} may be found in   \cite[Lemma 18]{borchedevkonmey25}.   
The corresponding conditional expectations are as follows:  for a deterministic constant $ b_{18}$ depending only on $T$,
\begin{align}
&  \sqrt{\alpha}  \Expect\big[   \big \|  \sum_{j=l}^{k-1}  \ScaledST_{j+1, k-1} \clE_j^T \big \|^4  \, \big | \,  \clF_{m_n} \big]^{1/4}  
\le   \alpha  b_{18}  \Expect\big[    \sum_{j=l}^{k-1}    \|   z^{(n)}_j      \|^4 \, \big | \,  \clF_{m_n} \big] ^{1/4}
\\
&  \sqrt{\alpha}  \Expect\big[   \big \|  \sum_{j=l}^{k-1}  \ScaledST_{j+1, k-1} [\clE_j^D - \alpha \Oops_{j+1}] \big \|^4  \, \big | \,  \clF_{m_n} \big] ^{1/4}   \le
\alpha b_{18} \clV^{1/4}(\Psi_{m_n} )\sqrt{\SAtime_k - \SAtime_l}
\label{e:z-D-bdd}
\\
&  \sqrt{\alpha}  \Expect\big[   \big \|  \sum_{j=l}^{k-1}  \ScaledST_{j+1, k-1}  [\clT_{j+1}-\clT_{j}] \big \|^4  \, \big | \,  \clF_{m_n} \big] ^{1/4}   \le \alpha  b_{18}   \clV^{1/4}(\Psi_{m_n} )\sqrt{\SAtime_k - \SAtime_l}
\label{e:z-tele-bdd}
\\
&  \sqrt{\alpha}  \Expect\big[   \big \|  \sum_{j=l}^{k-1}  \ScaledST_{j+1, k-1} \MD_{j+1} \big \|^4  \, \big | \,  \clF_{m_n} \big] ^{1/4}   \le    b_{18}   \clV^{1/4}(\Psi_{m_n} )\sqrt{\SAtime_k - \SAtime_l}.
\label{e:z-Mart-bdd}
\end{align}
The first three inequalities require \Cref{t:results_CLT} and the following bounds:   
\[
\begin{aligned}
\|\clE_k^T \|  &=  \sqrt{\alpha} O(  \|z_k^{(n)}\|)    
\,, 
&
\| \clE_k^D\| & = O(\alpha \| \theta_k\| ),
\\
\| \clT_{k} \|&  \le b_{\haf}  [1+V(\Phi_{k} )]   \|   \theta_k \| \,,
\qquad
& \|  \Oops_{k+1} \|& \le      b_{\haf}   [1+V(\Phi_{k+1} )]    \|  f( \theta_k ,\Phi_{k+1}) \|.
\end{aligned}
\]
The bound \eqref{e:z-tele-bdd} also requires summation by parts.  
The final bound \eqref{e:z-Mart-bdd} uses 
Burkholder's inequality (see \cite[Lem.~6, Ch.~3, Part~II]{benmetpri12})
combined with    
\[
\Expect[ 
\| \MD_{n+1}   \|^2  \mid \clF_{n} ]   \le b_f ^2  [ 1 + \| \theta_n\|]^2   [1+V(\Phi_{n})]^2,
\]
which follows from \Cref{t:bounds-H}.   See   \cite[Lemma~8]{borchedevkonmey25} for details on this step.  
We now set $l=m_n$.
Recalling that  $z_{m_n}^{(n)} =0$, the bounds above combined with 	\eqref{e:zn-acc-sum} imply the desired recursive bound 
\eqref{e:15BG}.  
\label{e:AllTheConditionalBdds}
\end{subequations}
\qed

\begin{subequations}

\begin{corollary}
\label[corollary]{t:ConditionalMS}
Under the assumptions of \Cref{t:tight},  if $T>0$ satisfies  \eqref{e:odeContractionT}, then  there is $\bdd{t:ConditionalMS}$ independent of $\alpha$ such that the following hold a.s.\ 
for $p=2$ and $p=4$:  
\begin{align}
	\Expect\bigl[   \| \tiltheta_{ k} \|^p  \mid \clF_{m_n} \bigr]  & \leq  \bdd{t:ConditionalMS} \| \tiltheta_{m_{n} } \|^p   +  p  \bigl[ \alpha^2   \bdd{t:tight}  \clV(\Psi_{m_n}) \bigr]^{p/4},       &&    m_n \le k\le m_{n+1},
	\label{e:ConditionalMSa}
	\\
	\Expect\bigl[   \| \tiltheta_{m_{n+1} } \|^p  \mid \clF_{m_n} \bigr]  & \leq   \half \| \tiltheta_{m_{n} } \|^p  +  p  \bigl[ \alpha^2   \bdd{t:tight}  \clV(\Psi_{m_n}) \bigr]^{p/4},       &&   n\ge 0\, .
	\label{e:ConditionalMS}
\end{align}
\end{corollary}

\end{subequations}

\Proof
For simplicity we restrict to $p=2$.
The proof combines three bounds:

\wham{1.}  
$
\Expect\bigl[\|z_{k} ^{(n)}\|^2  \mid \clF_{m_n} \bigr]  \leq  \sqrt{\bdd{t:tight} } \clV^{1/2}(\Psi_{m_n})  
$ follows from
\Cref{t:tight} and Jensen's inequality.

\wham{2.}    
$  \| \tiltheta_{ k } \|  \le   \|   \odestate_{\SAtime_k }^{(n)}  - \theta^*\|  +   \sqrt{\alpha} \|z_{k} ^{(n)}\| $ follows from   the definitions in \eqref{e:scaled-error-def}.
Moreover, the Lipschitz assumption on $\barf$ implies that there is
$\ell>0$ such that $ 
\|   \odestate_{\SAtime_k }^{(n)}  - \theta^*\|  \le    
e^{\ell [ \SAtime_k -T_n] }    \|  \tiltheta_{m_n}  \|$.    

\wham{3.}    
$ \|   \odestate_{T_{n+1} }^{(n)}  - \theta^*\| \le   \half \| \tiltheta_{m_{n} } \|^2$  holds in view of
\eqref{e:odeContractionT}.

The bounds in  2 and 3 apply 
the initialization    $ \odestate_{T_{n}}^{(n)}  = \theta_{m_n}$.

\smallskip

Combining 1 and 2   implies  \eqref{e:ConditionalMSa}  with  $ \bdd{t:ConditionalMS}  = 2 e^{\ell [ T_{n+1} -T_n] } $  (which is independent of $n$).

Next, we combine 2 and 3   using $k=m_{n+1}$ to conclude that 
$
\| \tiltheta_{m_{n+1} } \|^2  \le     \half \| \tiltheta_{m_{n} } \|^2 +    2 \alpha \|z_{m_{n+1}} ^{(n)}\|^2 
$, 
which when combined with 1 gives \eqref{e:ConditionalMS}.   
\qed

\whamit{Proof of \Cref{t:BigBounds}.}
Part (i) follows from \Cref{t:results_CLT} (i).  
Part~(ii) easily follows from   \Cref{t:ConditionalMS} with $p=4$:
\[
\begin{aligned}
\limsup_{k\to\infty} 
\Expect\big[   \| \tiltheta_{ k} \|^4   \bigr]  & \leq  \bdd{t:ConditionalMS}   \limsup_{n\to\infty} \Expect\bigl[    \| \tiltheta_{m_{n} } \|^4   \bigr] 
+  4    \alpha^2   \bdd{t:tight} b_\clV
\\
\textit{where} \ \
\Expect\bigl[   \| \tiltheta_{m_{n+1} } \|^4  \bigr]  & \leq   \half     \Expect\bigl[    \| \tiltheta_{m_{n} } \|^4   \bigr] 
+  4    \alpha^2   \bdd{t:tight}  \Expect\big[   \clV(\Psi_{m_n}) \bigr]  \,,  \ \ n\ge 0\,,
\end{aligned}
\]
and $b_\clV =
\limsup_{n\to\infty} \Expect [   \clV(\Psi_{m_n}) ] $.   
\Cref{t:results_CLT} implies that $b_\clV \le   \bdd{t:results_CLT} / (1 -   \bdde{t:results_CLT}  )$.
The recursive bound above gives
\[
\limsup_{n\to\infty}  \Expect\bigl[   \| \tiltheta_{m_{n+1} } \|^4  \bigr]    \le 
8  \alpha^2   \bdd{t:tight}    b_\clV  .
\]
This establishes  \Cref{t:BigBounds}~(ii)   with $ \bdd{t:BigBounds} =   4 \bdd{t:tight} [
1+ 2 \bdd{t:ConditionalMS}   ]     b_\clV$.
\qed

\wham{Sample path moment bounds.}   The conditional moment bounds in 
\Cref{t:results_CLT} and 
\Cref{t:tight} lead to sample path bounds required in analysis of the sensitivity process.

\begin{subequations}

\begin{proposition}
\label[proposition]{t:BigBddLLN}
Under  {\em (A1)-(A5)} we have, with probability one, from each initial condition,
\begin{align}
	\limsup_{N\to\infty} \frac{1}{N} \sum_{n=1}^N  \clV^{1/2}(\Psi_{m_n })   
	& \le   \bdd{t:BigBddLLN}_\nu
	\eqdef    \frac{ \bdd{t:results_CLT} }{ 2\bdde{t:results_CLT}} \frac{1}{1- \bdde{t:results_CLT}^{1/2} }  
	\label{e:BigBddLLNsqrtV}
	\\
	\limsup_{N\to\infty} \frac{1}{N} \sum_{n=1}^N  \| \tiltheta_{m_n} \|^2  
	& \le   \alpha   \bdd{t:BigBddLLN}_\tTheta \eqdef  \alpha  \big[ 4  \bdd{t:BigBddLLN}_\nu   \sqrt{ \bdd{t:tight}  } \big] ,
	\label{e:BigBddLLNtheta}
\end{align}
where $ \bdd{t:results_CLT} $,  $ \bdde{t:results_CLT} $ are given in \Cref{t:results_CLT},
and $\bdd{t:tight}$ is given in \Cref{t:tight}.
\end{proposition}

\end{subequations}

\Proof
The proof  appeals to the strong law of large numbers for martingales.   In the case of \eqref{e:BigBddLLNsqrtV} we apply the following:
\begin{equation}
\lim_{N\to\infty} \frac{1}{N} \sum_{n=1}^N  \clD_{n+1} =0\,, \quad \textit{with} \ \  \clD_{n+1} \eqdef   \clV^{1/2}  (\Psi_{m_{n+1} })   -\Expect[ \clV^{1/2}  (\Psi_{m_{n+1} }) \mid \clF_{m_n }   ].
\label{e:SLLNmart}
\end{equation}
The limit holds due to     \Cref{t:results_CLT}, which  implies that $\{ \clD_{n+1} \}$  is a square-integrable martingale difference sequence.  

The inequality 	\eqref{e:Lcontraction} combined with Jensen's inequality gives for each $n$,
\[
\begin{aligned}
\Expect[ \clV^{1/2}  (\Psi_{m_{n+1} }) \mid \clF_{m_n }   ] 
&  \le \Bigl(  \bdde{t:results_CLT}   \clV(\Psi_{m_n})    + \bdd{t:results_CLT}  \Bigr)^{1/2}
\le  \bdde{t:results_CLT}^{1/2}     \clV^{1/2}  (\Psi_{m_n})    +         \frac{1}{ 2\bdde{t:results_CLT}  }\bdd{t:results_CLT},
\end{aligned}
\]
where the final inequality uses concavity of the square root, and the fact that $\clV(z)\ge 1$ for all $z$.  
From the definitions this gives the lower bound $ \clD_{n+1}   \ge  (1 -   \bdde{t:results_CLT}^{1/2}   )   \clV^{1/2}  (\Psi_{m_n})    -    \bdd{t:results_CLT} /(2   \bdde{t:results_CLT}  )$,  and \eqref{e:BigBddLLNsqrtV} follows from \eqref{e:SLLNmart}.

The proof of  \eqref{e:SLLNmart}  is based on 
$
\{ \clD^{\tTheta}_{n+1} \eqdef   \| \tiltheta_{m_{n+1} } \|^2   -\Expect[  \| \tiltheta_{m_{n+1} } \|^2  \mid \clF_{m_n }   ]
:  n\ge 0 \}$, which is also a square-integrable martingale difference sequence.     The upper bound \eqref{e:ConditionalMS}
is equivalent to   the lower bound  $  \clD^{\tTheta}_{n+1} \ge  - 2 \alpha  \sqrt{\bdd{t:tight} } \clV^{1/2}(\Psi_{m_n})  + \half \| \tiltheta_{m_{n} } \|^2 $.    The remainder of the proof of \eqref{e:SLLNmart} follows again from the strong law of large numbers for square integrable martingales. 
\qed

\subsection{Coupling and ergodicity for nonlinear SA}
\label{s:CouplingNonlinear}

We begin with a few simple observations to motivate the sensitivity process.   
Our goal is to show that parameter sequences $\{ \theta_n ,\theta_n' : n\ge 0\} $   from distinct initial conditions couple, in the sense that $\lim_{n \to \infty} \| \theta_n - \theta_n' \| =0$, provided they are driven by the same Markov chain.    
\Cref{t:PhiCouples}  is used to establish coupling when  $ \{ \theta_n' \}$ is a stationary version of the parameter process:
Let $\theta^r_0 = \theta_0 + r[ \theta_0^\infty - \theta_0]$ for $0\le r\le 1$,     and let  $\{\theta^r_n  : n\ge 0$ denote the parameter process with this initial condition,  and   $\Psi^r_n = (\theta^r_n ;  \Phi^\infty_n)$.   Finally, denote by 
$\{ \Sens_n(\theta_0^r) : n\ge 0\}$   the sensitivity process with parameter process initialized at $\theta_0^r $.   
We then obtain by the fundamental theorem of calculus,
\begin{equation}
\theta_n - \theta_n^\infty =   \clA_n  [ \theta_0 - \theta_0^\infty ]  \,, \quad \clA_n \eqdef   \int_0^1 \Sens_n(\theta_0^t) \, dt .
\label{e:SensDelta-theta}
\end{equation}

Bounds on the matrix process $\{ \clA_n \}$ involve 
a \textit{partial sensitivity process}:
for each $n$ and $k\ge m_n$ denote $ \Sens_k^{(n)}  = \partial_{\theta_{m_n}}  \theta_k$ for $k\ge m_n$.   Similar to 	\eqref{e:Sens1}, this evolves according to the random linear system,
\begin{equation}
\Sens_{k+1}^{(n)} -\Sens_k^{(n)}   =      \alpha A_{k+1}   \Sens_k^{(n)}    \,,  \qquad k\ge m_n\,, \ \    \Sens_{m_n}^{(n)}  =I.
\label{e:Sens-k}
\end{equation}
The original sensitivity process 
can be recovered:    by the chain rule,  for $n\ge1$ and $k\ge m_n $,
\begin{equation}
\Sens_k =    \Sens_k^{(n)}  \Sens_{m_{n}}^{(n-1)}   \cdots \Sens_{m_2}^{(1)}      \Sens_{m_1}^{(0)} .
\label{e:SensFromSens-k}
\end{equation}

\begin{subequations}

Much of the work in this section involves the proof of the following.

\begin{proposition}
\label[proposition]{t:PartialSensBdd} 
Suppose that   {\em (A1)-(A5)}  hold.   Then, there is $\alphaTmp{t:PartialSensBdd} >0$ such that 
for any $\alpha \in (0,\alphaTmp{t:PartialSensBdd} ]$,

\whamrm{(i)}   
There is a   constant $ \bdd{t:PartialSensBdd}<\infty$ such that  for any initial condition $\Psi_0$ and any $n\ge 0$,   
\begin{align}
	\Expect\bigl[   \|\Sens^{(n)}_{m_{n+1}}\|_F^2   \mid    \clF_{m_n}   \bigr]  
	&  \le      \half  \bigl( 1 
	+   \bdd{t:PartialSensBdd}   \big [   \| \tiltheta_n\|^2 +  \alpha  \clV^{1/2} (  \Psi_{m_n} )   \big] \big),
	\label{e:PartialSensBdd}
\end{align} 
and consequently $   \Expect\bigl[   \|\Sens^{(n)}_{m_{n+1}}\|^4_F   \mid    \clF_{0}   \bigr]      \le    \bdd{t:PartialSensBdd}  v(\Phi_0)$.

\whamrm{(ii)}     Consider $\bfPhi = \bfPhi^\infty$.    On  letting $\{  \Sens^{(n)}_k (\theta_0^r) \;;\; k\ge m_n \,, \ n\ge 0\}$  denote the partial sensitivity process with initial condition $\theta_0^r$,
\begin{equation}
	\int_0^1    \Expect\bigl[   \|\Sens^{(n)}_{m_{n+1}}(\theta_0^r) \|_F^2   \mid    \clF_{m_n}   \bigr]   \,  dr   
	\le      \half    \int_0^1   \bigl( 1 
	+   \bdd{t:PartialSensBdd}   \big [   \| \tiltheta_n^r\|^2 +  \alpha  \clV^{1/2} (  \Psi_{m_n}^r )   \big] \big) \,  dr .
	\label{e:PartialSensBddConvexIC}
\end{equation}
\end{proposition}

\label{e:PartialSensBddAll} 
\end{subequations}

Motivation for the proposition comes from the following:

\begin{proposition}
\label[proposition]{t:SensLLN}
Under  {\em (A1)-(A5)},
there is  $\bdde{t:SensLLN}>0$  and 
$\alphaTmp{t:SensLLN}  \in (0,
\alphaTmp{t:PartialSensBdd} ]$  
such that the following hold 
for any initial condition 
$\Psi_0 = (\theta_0; x)$,  and any step-size $\alpha\in (0,\alphaTmp{t:SensLLN} ]$.

\whamrm{(i)}
If   \eqref{e:PartialSensBdd}  holds then  
$\displaystyle 
\limsup_{n\to\infty} \frac{1}{n} \log  \|\Sens_n\|_F  \le - \bdde{t:SensLLN} \alpha     $  \   a.s.

\whamrm{(ii)}
If   \eqref{e:PartialSensBddConvexIC}  holds  along with {\em (A1)--(A5)}, 
then \Cref{t:sens}  (iii) holds:   there is a unique invariant measure $\upvarpi$ for $\bfPsi$,  and  for each initial condition $\Psi_0$  there is a a stationary realization  $\bfPsi^\infty = \{ \Psi^\infty_n  : -\infty <n < \infty\}$ on the same probability space such that 
\begin{equation}
\limsup_{n\to\infty}  \frac{1}{n} \log  \| \theta_n - \theta_n^\infty \|  
\le
\limsup_{n\to\infty} \frac{1}{n} \log  \|\clA_n\|_F  
\le   -\bdde{t:SensLLN}\alpha.
\label{e:t:sensiii}
\end{equation}
\end{proposition}

\Proof
\Cref{t:PhiCouples}  allows us to construct, 
for each initial condition $\Psi_0 =z$,
a joint realization
$\{\Psi_n,  \Psi^\infty_n  :  0 <n < \infty\}$ in which $\bfPhi^\infty$ is stationary,  and the two processes couple at time $\taucpl$.  
Since our goal is to establish almost-sure results,   we may restrict attention to $n\ge \taucpl$.   To simplify notation, in the following we replace $\bfPhi$ by $\bfPhi^\infty$.   In particular,  $\Psi_n = (\theta_n;\Phi_n^\infty)$ for all $n$.    

Applying 		
\eqref{e:SensFromSens-k},
for $ m_n \le k  < m_{n+1}$,
\[
\frac{1}{k}
\log \| \Sens_k \|_F  
\le  \frac{1}{k}  \log \|  \Sens_k^{(n)} \|_F
+ 
\frac{\alpha }{T}   \frac{1}{n}  \sum_{j=1}^n   \clZ_{m_j}  \,,  \quad \clZ_{m_j} \eqdef  \half   \log \|  \Sens_{m_{j}}^{(j-1)}  \|_F^2 ,
\]
where we have used the fact that $  k  = Tn /\alpha +O(1)  $ with $O(1)$ is bounded as $k\to\infty$ for the range of $k$ considered.

If  \eqref{e:PartialSensBdd} holds, then by Jensen's inequality and the bound $\log(1+x) \le x$,
\[
\Expect\bigl[  \clZ_{m_{j+1}}    \mid    \clF_{m_j}   \bigr]   
\le 
- \half \log(2)   
+\half   \bdd{t:PartialSensBdd}   \big [   \| \tiltheta_{m_j} \|^2 +  \alpha  \clV^{1/2} (  \Psi_{m_j} )   \big] .
\]
Appealing to the strong law of large numbers for martingales,   as in the proof of  \Cref{t:BigBddLLN},
\[
\limsup_{n\to\infty}
\frac{1}{k}
\log \| \Sens_k \|_F  
\le   
\frac{\alpha}{T}
\limsup_{n\to\infty}
\frac{1}{n}  \sum_{j=1}^n     \Expect\bigl[  \clZ_{m_{j+1}}    \mid    \clF_{m_j}   \bigr]     
\le
-
\frac{1}{2}
\frac{\alpha}{T} \Bigl (
\log(2)   
-   \alpha \bdd{t:PartialSensBdd}   \big [   \bdd{t:BigBddLLN}_\tTheta  +      \bdd{t:BigBddLLN}_\nu  \big] 
\Bigr),
\]
where the second inequality follows from \Cref{t:BigBddLLN}.     This establishes (i),    with $\bdde{t:SensLLN} <   \log(2)  /(2T)$,   
and $\alphaTmp{t:SensLLN}   <  
\log(2)   \big[  \bdd{t:PartialSensBdd}     [   \bdd{t:BigBddLLN}_\tTheta  +      \bdd{t:BigBddLLN}_\nu  ] \big]^{-1}$.

We now establish (ii).   
The proof of  the second inequality in \eqref{e:t:sensiii}
is identical to the proof of (i):  
If \eqref{e:PartialSensBddConvexIC} holds then,
\[
\begin{aligned}
\limsup_{n\to\infty}  
\frac{1}{n} \log \| \clA_n\|_F &  \le \limsup_{n\to\infty}  \frac{1}{n} \sum_{k=1}^n  \log   \int_0^1    \Expect\bigl[   \|\Sens^{(n)}_{m_{n+1}}(\theta_0^r) \|_F^2   \mid    \clF_{m_n}   \bigr]   \,  dr   
\\
&
\le     - \half \log(2)  + \half    \bdd{t:PartialSensBdd}  
\limsup_{n\to\infty}  \frac{1}{n} \sum_{j=1}^n \int_0^1\big [   \| \tiltheta_{m_j}^r\|^2 +  \alpha  \clV^{1/2} (  \Psi_{m_j}^r )   \big]   \, dr    
\\
&      \le    -\bdde{t:SensLLN}\alpha,
\end{aligned}
\]
where the final inequality is a straightforward generalization of \Cref{t:BigBddLLN},  
with  $\bdde{t:sens}$   identical to the value in part (i). 

The existence  of  $\bfPsi^\infty  $  follows from \Cref{t:BigBounds},
and then \eqref{e:t:sensiii} follows from  the representation
\eqref{e:SensDelta-theta}.
\qed

\smallskip

\Cref{t:DV3multBdd} will be applied as a first step in the proof,  using $\delta_k \equiv p\alpha$ for $0\le k <N$  with $p\ge 1$,  
and $N = m_{n+1} -m_n =  \lceil T/\alpha \rceil$ is independent of $n$.

\begin{lemma}
\label[lemma]{t:bddsSens}
Suppose that   {\em (A1)-(A5)}  hold.   
Then, there is $\bdd{t:bddsSens}<\infty$ such that  for $0<\alpha \le \alphaTmp{t:SensLLN}$ and each~$n$, 
\[
\sup_{\theta,x}   \frac{1}{v(x)}  \max_{m_n\le k \le m_{n+1} }   \Expect\Bigl[   \|\Sens^{(n)}_k\|^4  v(\Phi_{k+1}) \mid  \theta_{m_n} = \theta \, ,  \Phi_{m_n} =x\Bigr]  
\le \bdd{t:bddsSens}.
\]
\end{lemma}

\Proof 
Recall $\Sens^{(n)}_{m_n} = I$.    
From \eqref{e:Sens-k} we obtain for $m_n < k \le m_{n+1} $,
\[
\| \Sens_k^{(n)}   \|^4 =    \|  [ I+\alpha A_k]  [ I+\alpha A_{k-1}]\cdots [ I+\alpha A_{m_n}]  \| ^4
\le \exp\Bigl( 4 \alpha \sum_{i=m_n+1}^{k }   \| A_i\| \Bigr).
\]
Next apply  (A2i) to obtain  $ \| A_i\| \le L(\Phi_i)$ for each $i$,  so that on applying \eqref{e:DV3multBdd+L}  of 
\Cref{t:DV3multBdd},
\[
\Expect\Bigl[   \|\Sens^{(n)}_k\|^4  v(\Phi_{k+1}) \mid  \clF_{m_n + 1}  \Bigr]
\le
\Expect\Bigl[  \exp\Bigl(  4 \alpha \sum_{i=m_n+1}^{k}   L(\Phi_i)    + V(\Phi_{k+1})  \Bigr)  \mid  \clF_{m_n + 1}  \Bigr]
\le
\bdd{t:DV3multBdd} (r)   v(\Phi_{m_n + 1} ),
\]
where $r= 4[T+\alpha]$.     Under (DV3) the desired conclusion follows:
\[
\Expect\Bigl[   \|\Sens^{(n)}_k\|^4  v(\Phi_{k+1}) \mid  \Psi_{m_n} = (\theta; x)\Bigr] 
\le 
\bdd{t:DV3multBdd} (r)  \Expect[    v( \Phi_{m_n+1} )   \mid \Phi_{m_n} =x   ]  
\le   \bdd{t:DV3multBdd} (r)  v(x).
\] 
\qed

\wham{A disturbance decomposition for the sensitivity process.}

Denote  $\haderf(\theta,x) = \partial_\theta \haf(\theta,x)  $  for $ \theta\in\Re^d$,  $x\in\state$.   
It is easily shown that $\haderf(\theta,\varble)  $   solves Poisson's equation with forcing function $A(\theta,\varble)$ for any $ \theta\in\Re^d$.  The lemma that follows is established exactly as in \Cref{t:noise-decomp}. The bounds in \eqref{e:AdecompBdds}
hold under the Lipschitz bound in (A2), exactly as in the proof of \Cref{t:bddsSens}.

Denote   $\uppsiA (\theta,x) \eqdef  A(\theta,x) - 
\haderf(\theta,x)$  for $\theta\in\Re^d$,  $x\in\state$.

\begin{subequations}

\begin{lemma}
\label[lemma]{t:Adecomp}
Suppose that   {\em (A1)-(A5)}  hold.   
Then,
$
A_{k+1}   =    \derbarf(\theta_k) +  \MDA_{k+1} - \clTA_{k+1} + \clTA_{k} - \alpha\OopsA_{k+1}
$ for $k\ge 0$,   in which   
\begin{align}
	\MDA_{k+1} &= \haderf(\theta_k, \Phi_{k+1}) - \Expect[\haderf(\theta_k, \Phi_{k+1})  \mid \clF_{k}],  
	\label{e:MDSens}
	\\
	\clTA_{k+1} &=  \uppsiA (\theta_{k+1}, \Phi_{k+1}),
	\label{e:telescopeSens}
	\\
	\OopsA_{k+1} &=  \frac{1}{\alpha} \bigl[ \uppsiA (\theta_{k+1}, \Phi_{k+1}) - \uppsiA (\theta_{k},  \Phi_{k+1}) \bigr].
	\label{e:UpsilonSens}
\end{align}
Moreover, the terms admit bounds independent of $\{\theta_k\}$: for a 
fixed constant $\bdd{t:Adecomp}$,      
\begin{equation}
	\begin{aligned}
		&
		\| \MDA_{k+1} \|_F  \le    \bdd{t:Adecomp} \bigl[1+ V(\Phi_{k+1})+ \Expect[V(\Phi_{k+1}) \mid \clF_k ]  \bigr]
		\\
		&  \|  \OopsA_{k+1}\|_F +  \|\clTA_{k+1}\|_F \le    \bdd{t:Adecomp} [1+V(\Phi_{k+1}) ] .
	\end{aligned} 
	\label{e:AdecompBdds}
\end{equation}
\end{lemma}

\label{e:noise-decompSens}
\end{subequations}

Application of this disturbance decomposition is very different from the ODE analysis used to establish convergence of the parameter process.   In \Cref{t:bound-taylor-reminder} we justify the representation, 
\begin{equation}
\Sens_{k+1}^{(n)}    =   \DetST \Sens_k^{(n)}  +   \alpha  \DistSens_{k+1}   \,,\quad \DetST=[I +      \alpha  \derbarf^*   ] \,, 
\label{e:LinearSensitivity}
\end{equation}
in which $\{  \DistSens_{k+1}  \}$ satisfies useful moment bounds.

\begin{lemma}
\label[lemma]{t:bound-taylor-reminder}
Suppose that   {\em (A1)-(A5)}  hold.   The representation 
\eqref{e:LinearSensitivity}
then follows with
$  \DistSens_{k+1}   
=  \DistSensND_{k+1}   
+  \DistSensTS_{k+1}$, in which  $  \DistSensTS_{k+1}   =  [\derbarf(\theta_k)  -  \derbarf^* ]\Sens_k^{(n)}  $  and
\begin{equation}
\begin{aligned}
	\DistSensND_{k+1}    =   \MDSens_{k+1}    -  [\clTSens_{k+1} &- \clTSens_{k} ]  - \alpha \OopsSens_{k+1}  
	\\
	\text{with}  \quad    \MDSens_{k+1}   &  \eqdef \MDA_{k+1}  \Sens_{k}^{(n)}   \,, \qquad   \clTSens_{k+1}    \eqdef    \clTA_{k+1}  \Sens_{k+1}^{(n)}  \,,
	\\
	\OopsSens_{k+1}   &  \eqdef   \OopsA_{k+1}    \Sens_k^{(n)}  +   \clTA_{k+1}  A_{k+1} \Sens_{k}^{(n)}.
\end{aligned}
\label{e:DeltaSensDecomp}
\end{equation}
\end{lemma}

\Proof 
We have \eqref{e:LinearSensitivity} with $  \DistSens_{k+1}   
=  \DistSensND_{k+1}   
+  \DistSensTS_{k+1}$,   and applying \Cref{t:Adecomp} gives
$ \DistSensND_{k+1}     =  [  - \clTA_{k+1} + \clTA_{k} - \alpha\OopsA_{k+1} ] \Sens_k^{(n)} $.  
The proof is completed on adding and subtracting $ \clTA_{k+1}  \Sens_{k+1}^{(n)} $ and employing \eqref{e:Sens-k}.  
\qed

\smallskip

The representation \eqref{e:LinearSensitivity} is similar to \eqref{e:zn-acc-sum}, but the analysis that follows is far less complex than our treatment of $\{ z_k^{(n)} \}$,  largely because we    
easily obtain strong moment bounds on  $ \{  \DistSens_{k}  \;;\; k > m_n  \} $:
\Cref{t:ConditionalMS}   is used to obtain bounds on  $ \{  \DistSensTS_{k} \;;\; k > m_n  \} $,   and
\Cref{t:Adecomp}  is used in consideration of  $ \{ \DistSensND_{k} \;;\; k > m_n  \} $.

The next result implies that  the dominant term in $ \{  \DistSens_{k}  \;;\; k> m_n  \} $,  in terms of 
its impact on the variance of $\{ \Sens_k^{(n)}  \}$, 
is the martingale difference sequence $ \{     \MDSens_{k} \;;\; k> m_n \}$.

\begin{lemma}
\label[lemma]{t:SensRealization}
The representation
$\Sens_{k+1}^{(n)} = \clX_{k+1}^{(n)}  +  \clY_{k+1}^{(n)}  +   \alpha  \clU_{k+1}^{(n)}  $ 
holds for each $n\ge 0$ and $k\ge m_n$, 
where    
$ \clU_{k+1}^{(n)}   \eqdef \DetST^{-1}  \clTSens_{k+1} -  \DetST^{k-m_n}  \clTSens_{m_n}  $,  		
with $\DetST$ defined in \eqref{e:LinearSensitivity},
and
\begin{equation}
\begin{aligned}
	\clX_{k+1}^{(n)}  &  =   \DetST\clX_k^{(n)}  +   \alpha      \MDSens_{k+1}      \,,   \quad
	&&    \clX_{m_n}^{(n)}   =  I
	\\
	\clY_{k+1}^{(n)}   & =   \DetST\clY_k^{(n)}  +      \alpha    \DistSensY_{k+1}       \,,   \quad
	&&    \clY_{m_n}^{(n)}   =  0
	\\
	&&&  
	\DistSensY_{k+1}    =    \DistSensTS_{k+1}  -  \alpha\OopsSens_{k+1}  
	- \alpha \derbarf^* \DetST^{-1}\clTSens_{k+1},
\end{aligned} 
\label{e:LinearSensitivityNicer}
\end{equation}
in which $\DetST=[I +      \alpha  \derbarf^*   ] $.
\end{lemma}
\Proof 
From \eqref{e:LinearSensitivity} and the initialization $ \Sens_{m_n}^{(n)}  =I$,
we obtain for $k > m_n$,
\begin{equation}
\begin{aligned}
\Sens_{k+1}^{(n)}   & =  \DetST^{k+1-m_n}   +   \alpha \sum_{j=m_n+1}^{k+1}    \DetST^{k+1-j}   \DistSens_{j}   
=     \clX_{k+1}^{(n)}   +    \alpha \sum_{j=m_n+1}^{k+1}    \DetST^{k+1-j}    \DistSensYa_j,
\end{aligned}
\label{e:UnravelS}
\end{equation}
where $\DistSensYa_j = 
\DistSens_{j}   - \MDSens_{j} = -[
\clTSens_{j} - \clTSens_{j-1}   +\alpha \OopsSens_{j} ]$.
To establish the representation \eqref{e:LinearSensitivityNicer}, apply summation by parts:
\[
\begin{aligned}
- \sum_{j=m_n+1}^{k+1}    \DetST^{k+1-j} [  \clTSens_{j} - \clTSens_{j-1} ]
&= - \DetST^{-1} \clTSens_{k+1}  +  \DetST^{k-m_n} \clTSens_{m_n}
+ \sum_{j=m_n+1}^{k+1} (\DetST^{k-j}  - \DetST^{k+1-j}  ) \clTSens_{j}
\\
&= \clU_{k+1}^{(n)}  
-  \alpha  \sum_{j=m_n+1}^{k+1} \DetST^{k+1-j} \derbarf^* \DetST^{-1}   \clTSens_{j},
\end{aligned}
\]
so that $ \alpha \sum_{j=m_n+1}^{k+1}    \DetST^{k+1-j}    \DistSensYa_j 
=  \clY_{k+1}^{(n)}  +  \alpha  \clU_{k+1}^{(n)} $.
This combined with \eqref{e:UnravelS} completes the proof.
\qed

\begin{lemma}
\label[lemma]{t:DiscretizeMatrixExponential}
Suppose that $\derbarf^*$ is Hurwitz.  Then,  there exists $\alphaTmp{t:DiscretizeMatrixExponential}>0$,
$\bdd{t:DiscretizeMatrixExponential}<\infty$,
and 
$\bddepsy{t:DiscretizeMatrixExponential} >0$ such that
\[
\| [I+\alpha \derbarf^*]^n  -   e^{n\alpha \derbarf^*} \|_F  \le \alpha^2 n \bdd{t:DiscretizeMatrixExponential}
\exp(-\bddepsy{t:DiscretizeMatrixExponential} \alpha n) 
\le  \frac{  \bdd{t:DiscretizeMatrixExponential} } { \bddepsy{t:DiscretizeMatrixExponential} }   \alpha \,, 
\quad n\ge 0.
\]
\end{lemma}

\Proof 
Write  $  [I+\alpha \derbarf^*]^n  =   e^{n\alpha A^\alpha} $   where $  A^\alpha  \eqdef \alpha^{-1} \log(I+ \alpha \derbarf^* )$.
Let  $R^\alpha =   [ A^\alpha - \derbarf^*]/\alpha$, giving  $ [I+\alpha \derbarf^*]^n    -  e^{n\alpha \derbarf^*}    =      e^{n\alpha \derbarf^*} \big[     e^{n\alpha^2 R^\alpha}  -  I  \big]$.     In analogy with the scalar bound $e^x -1  \le x e^x$ for $x\ge 0$, we have
\[
\| e^{n\alpha^2 R^\alpha}  -  I \|_F   \le \sum_{k=1}^\infty  \frac{1}{k!} \big[ n\alpha^2 \| R^\alpha \|_F]^k  
\le   \alpha^2 \big[  n  \| R^\alpha \|_F 
\exp( n\alpha^2 \| R^\alpha \|_F) \big].
\]
Consequently,  
\[
\begin{aligned}
\|  [I+\alpha \derbarf^*]^n    -   e^{n\alpha \derbarf^*}  \|_F
&    \le \|     \exp(n\alpha \derbarf^*)   \|_F   \|  e^{n\alpha^2 R^\alpha}  -  I   \|_F
\\
&  \le  
\alpha^2  n    \| R^\alpha \|_F   \big[    \|   \exp(n\alpha \derbarf^*)  \|_F     
\exp( n\alpha^2 \| R^\alpha \|_F)  \big].
\end{aligned} 
\]
From the Taylor series approximation of the logarithm we have   $\lim_{\alpha\downarrow 0} R^\alpha = -\half (\derbarf^*)^2$.     
The proof follows from these bounds and  the Hurwitz assumption  on $\derbarf^*$.
\qed

\smallskip

%

\whamit{Proof of \Cref{t:PartialSensBdd}.}    
The proof of  \eqref{e:PartialSensBdd} begins with an application of \Cref{t:SensRealization} to obtain for each $n$,
\begin{equation}
\Expect\bigl[   \|\Sens^{(n)}_{m_{n+1}}\|_F^2   \mid    \clF_{m_n}   \bigr]^{1/2} 
\le  B_{\clX}(n) + B_{\clY}(n) +  \alpha B_{\clU}(n),
\label{e:PartialStep1}
\end{equation} 
where  $  B_{\clX}(n) \eqdef
\Expect\bigl[   \|\clX^{(n)}_{m_{n+1}}\|_F^2   \mid    \clF_{m_n}   \bigr]^{1/2}$,  and the definitions of the other two terms are analogous.  

\begin{subequations}
From \eqref{e:LinearSensitivityNicer} we obtain
$
\clX_{k+1}^{(n)}   = 
\DetST^{k+1-m_n}   +   \alpha \sum_{j=m_n+1}^{k+1}    \DetST^{k+1-j}   \MDSens_{j}    
$,
in which $\{ \MDSens_{j}    \}$  is a matrix-valued martingale difference sequence.   Consequently,
with $  \Sigma^{Ws}_{j} \eqdef \Expect[  \MDSens_{j}   ( \MDSens_{j} )^\transpose  ]$,
\[
\Expect\bigl[   \| \clX^{(n)}_{m_{n+1}}\|^2_F   \mid    \clF_{m_n}   \bigr]  
=    \|  \DetST^{m_{n+1}-m_n-1 } \|_F^2   +   \alpha^2 \trace  \Big( \sum_{j=m_n}^{m_{n+1} -1 }    \DetST^{m_{n+1}-j}  \Sigma^{Ws}_{j}   ( \DetST^{m_{n+1}-j}  )^\transpose  \Big).
\]
In view of \Cref{t:DiscretizeMatrixExponential} we have $  \|  \DetST^{m_{n+1}-m_n-1 } \|_F^2 =  \|  \DetST^{m_1 } \|_F^2  
\le   2  \| e^{\derbarf^*T} \|_F^2$ for all $\alpha>0$ sufficiently small,
so that by  \eqref{e:odeContractionT} we have $ \|  \DetST^{m_{n+1}-m_n-1 } \|_F^2 \le 1/8$.   
An application of 
\Cref{t:bddsSens} gives for some constant $b_X$,
\begin{equation}
\begin{aligned}
	B_{\clX}^2(n) =
	\Expect\bigl[   \|\clX^{(n)}_{m_{n+1}}\|^2_F   \mid    \clF_{m_n}   \bigr]  
	&   \le   \tfrac{1}{8}   +   b_X \alpha^2    [ m_{n+1} - m_n]   v^{1/2}(\Phi_{m_n})   
	\\
	&    \le
	\tfrac{1}{8}   +   b_X \alpha     [ T+ 1 ]   \clV^{1/2}(\Psi_{m_n}),
\end{aligned} 
\label{e:BXbdd}
\end{equation}
where the second inequality uses $m_{n+1} - m_n \le (T+1)/\alpha$.

The bound on $ B_{\clY}(n)$ is more complex.    As above we obtain from    \eqref{e:LinearSensitivityNicer},
\[
\clY_{k+1}^{(n)}   =    \alpha \sum_{j=m_n+1}^{k+1}    \DetST^{k+1-j}   \DistSensY_{j}    \,, \qquad k\ge m_n
\]
so that on applying the triangle inequality, 
\[   
\begin{aligned}
B_{\clY}(n) & =
\Expect\big[   \|\clY^{(n)}_{m_{n+1}}\|^2_F   \mid    \clF_{m_n}   \big] ^{1/2} 
\le    \alpha      \sum_{j=m_n+1}^{m_{n+1}}   \|   \DetST^{m_{n+1}+1-j}\|_F   \Expect \big [\|  \DistSensY_{j}  \|_F^2  \mid    \clF_{m_n}   \big] ^{1/2} 
\\
& \le  b_{\DetST}    \max_{m_n<j\le m_{n+1}}   \Expect \big [\|  \DistSensY_{j}  \|_F^2  \mid    \clF_{m_n}   \big] ^{1/2} 
\,,  \qquad  \text{with} \ \  b_{\DetST}   =  \alpha \sum_{k=0}^\infty  \|  \DetST^k\|_F .
\end{aligned}
\]
A second application of the triangle inequality gives 
\[
\begin{aligned}
B_{\clY}(n)  &\le   B_{\lilTS}(n) + \alpha  B_{\OopsSens}(n) 
+ \alpha \| \derbarf^* \DetST^{-1} \|_F  B_{\clTSens}(n) 
\\
& \text{with} \ \ 
B_{\lilTS}(n)   \eqdef
\max_{m_n<j\le m_{n+1}}   \Expect \big [  \|  \DistSensTS_{j}  \|_F^2  \mid    \clF_{m_n}   \big] ^{1/2} ,
\end{aligned}
\]
and the other terms defined similarly  (recall the definition  of $ \DistSensY_{k+1}$ in   \eqref{e:LinearSensitivityNicer}).   

Applying the definition  \eqref{e:DeltaSensDecomp} and \eqref{e:AdecompBdds}  gives for each $j$,    
\[
\begin{aligned}
\|\OopsSens_{j}    \|_F   & \le  \bigl[  \|  \OopsA_{j}  \|_F      +  \| \clTA_{j}  \|_F  \| A_{j}\|_F \big] \|  \Sens_{{j-1}}^{(n)} \|_F
\le     \bdd{t:Adecomp}     \bigl[   1      +   V(\Phi_j)   \big]^2 \|  \Sens_{{j-1}}^{(n)} \|_F
\\
\| \clTSens_{j}    \|_F   & \le
\|  \clTA_{j} \|_F  \|  \Sens_{j}^{(n)} \|_F   \le    \bdd{t:Adecomp}  V(\Phi_{j}) \|  \Sens_{j}^{(n)} \|_F.  
\end{aligned}
\]
We also have $[1+V(x)]^2 \le 4 v^{1/4}(x) $ for   all $x$, so by  \Cref{t:bddsSens} and Jensen's inequality, 
\begin{equation}
B_{\OopsSens}(n) + B_{\clTSens}(n)   \le   b_Y^a  \clV^{1/4}(\Psi_{m_n}) \ \  a.s.,  \quad \text{  for some $b_Y^a<\infty$.}
\label{e:bYa}
\end{equation}     

Letting  $\ell_A$ denote the  Lipschitz bound on $\derbarf$ we have   $\|  \DistSensTS_{j}  \|_F  \le \ell_A \| \tiltheta_{j+1} \| $.
Hence by \Cref{t:ConditionalMS},
$
B_{\lilTS}^2(n) 
\leq  \ell_A^2 \big[ \bdd{t:ConditionalMS} \| \tiltheta_{m_{n} } \|^2  + 2 \alpha  \sqrt{   \bdd{t:tight}  } \clV^{1/2} (\Psi_{m_n}) \big]  
$.
Combining this and \eqref{e:bYa} gives 
\begin{equation}
B_{\clY}(n)   \le  b_Y^a  \clV^{1/4}(\Psi_{m_n})  +  \ell_A \Big[ \bdd{t:ConditionalMS} \| \tiltheta_{m_{n} } \|^2  + 2 \alpha  \sqrt{   \bdd{t:tight}  } \clV^{1/2} (\Psi_{m_n}) \Big]^{1/2}.  
\label{e:BYbdd}
\end{equation}

Finally, precisely as in the proof of \eqref{e:bYa} we have
$
B_{\clU}(n)    \le  b_U  \clV^{1/4}(\Psi_{m_n}) $  $a.s.$,  
for some constant $b_U$.  This bound combined with  \eqref{e:BXbdd}, \eqref{e:BYbdd}  and the identity   \eqref{e:PartialStep1}   completes the proof of  \eqref{e:PartialSensBdd}.
\end{subequations}%
\qed

\smallskip

Armed with   \Cref{t:BigBounds} and the conclusions of \Cref{t:SensLLN}~(ii),   we establish a sequence of lemmas that will provide the proof of \Cref{t:sens} and parts of \Cref{t:sensLinear}.    In particular, we set $\alpha_0 = \min(\alphaTmp{t:BigBounds},\alphaTmp{t:PartialSensBdd})$ in  \Cref{t:sens}.
We begin with a crucial result related to part (ii) of  \Cref{t:sens}.

\begin{lemma}
\label[lemma]{t:AsyOops+Bias}
Under the assumptions of  \Cref{t:sens},
on setting $\alpha_0 = \min(\alphaTmp{t:BigBounds},\alphaTmp{t:SensLLN})$, the following hold for  $0<\alpha \le  \alpha_0$:  
\whamrm{(i)}   
If $g\colon\Re^d\times\state\to \Re$ satisfies $| g(z)|  \le  \clV(z) /[1 + \log(\clV(z))]$ for each $z$,    
and $g(\varble,x)$ is continuous on $\Re^d$ for each $x$, then for each initial condition
\[
\lim_{n\to\infty}
\Expect[  g(\Psi_n) ] 
=
\Expect_\upvarpi[  g(\Psi_0) ] . 
\]

\whamrm{(ii)}   
The limit \eqref{e:TargetBias} holds   with 
$\displaystyle
\barOops_\alpha \eqdef       \Expect_\upvarpi  [ \Oops_{n+1}]$.

\whamrm{(iii)}     $\|   \thbias \|     \le    \bdd{t:AsyOops+Bias} \alpha $,  with  $  \bdd{t:AsyOops+Bias}$ independent of $\alpha$.  
\end{lemma}

\Proof
The ergodic theorem in 
(i) follows from  \Cref{t:PhiCouples} and  \eqref{e:t:sensiii} if $g$ is bounded.  The extension to unbounded functions follows from 
\Cref{t:BigBounds} (implying uniform integrability) under the growth condition on~$| g(z)| $.

Combine \eqref{e:noisyEuler} 
and the decomposition \eqref{e:DeltaDecomp},
with $\alpha_n=\alpha$,
to obtain, for any $n$, 
\[
0 =    \Expect_{\upvarpi}[\theta_{n+1} - \theta_n ]   =   \alpha  \Expect_{\upvarpi}[ \barf(\theta_n) +  \MD_{n+1} - \clT_{n+1} + \clT_{n} - \alpha\Oops_{n+1} ].
\]
Invariance and the martingale difference property for $\{  \MD_{n+1} \}$ implies that $ \Expect_{\upvarpi}[ \barf(\theta_n)  ]=    \alpha \barOops_\alpha$       
with $\displaystyle
\barOops_\alpha \eqdef       \Expect_\upvarpi  [ \Oops_{n+1}]$,  which gives (ii).
Next, apply the fundamental theorem of calculus to write 
\begin{equation}
\barf(\theta_n)  =  \barf(\theta^*) +  \derbarf(\theta^*)  [\theta_n -\theta^*] 
+
\Bigl[ \int_0^1[  \derbarf(\theta^r_n) -\derbarf(\theta^*) ] \, dr\Bigr]  [\theta_n -\theta^*] ,
\label{e:FTC}
\end{equation}
where $\theta^r_n = (1-r) \theta^* + r\theta_n$, so that  $\theta^r_n - \theta^*  =  r [\theta_n -\theta^*] $.
Recalling $\barf(\theta^*) =0$,  the notation  $G^*=-[\derbarf^*]^{-1}$,
and letting $\ell_A$ be a Lipschitz constant for $A$, we obtain (iii)   with
\[ 
\bdd{t:AsyOops+Bias} =   \| G^*\|_F  \sup_{0<\alpha\le \alpha_0}
\big[   \barOops_\alpha   +    \half  \tfrac{1}{\alpha}   \,  \ell_A   \Expect_{\upvarpi}[   \|\theta_n -\theta^* \|^2 ].
\]
\Cref{t:BigBounds} implies that the supremum is finite.  
\qed
\begin{subequations}

\smallskip

We next  consider approximation of steady state expectations of 
$\Expect_\upvarpi[g(\Psi_0)]$ by  $\Expect_\upvarpi[\barg(\theta_0)]$
and by $\barg(\theta^*)$.  
\begin{lemma}
\label[lemma]{t:PolyMeans}
Under the assumptions of  \Cref{t:sens}, set $\alpha_0 = \min(\alphaTmp{t:BigBounds},\alphaTmp{t:SensLLN})$.   For   any $p\ge 1$  there exists $ 
\bdd{t:PolyMeans}= \bdd{t:PolyMeans}(p) < \infty$ such that the following holds:
if  $g\colon \Re^d \times\state\to\Re$ satisfies for some 
$b_g <\infty$,   and all $\theta, \theta'  \in  \Re^d$,  $ x\in \state$,
\[
| g (\theta, x) |   \le   b_g [1+\|\theta\|]^3 [1+V(x)]^p    \quad \textit{and} \quad
| g (\theta, x)  - g (\theta', x) |   \le   b_g [1+\|\theta\|  +\|\theta'\|  ]^2 [1+V(x)]^p  \|  \theta - \theta'\| \,  ,
\]
then for   $\alpha\in (0,\alpha_0]$  and  $k\ge 0$,
\begin{equation}
	\big |  \Expect_\upvarpi[g(\theta_0 , \Phi_k)]  -  \Expect_\upvarpi[\barg(\theta_0)]  \big |   \le  \alpha b_g \bdd{t:PolyMeans} .
	\label{e:GenMeanApprox} 
\end{equation}
Suppose in addition that 
\[
\| \derbarf_g (\theta)  - \derbarf_g (\theta') \|   \le   b_g [1+\|\theta\|  +\|\theta'\|  ]^2  \|  \theta - \theta'\|,
\]
for all $\theta,\theta'$,   where
$\derbarf_g \eqdef\partial_\theta \barg$.   Then,
\begin{equation}
	\big |    \Expect_\upvarpi[ \barg(\theta_0)]    -   \barg(\theta^*)  \big |   \le  \alpha b_g \bdd{t:PolyMeans} .
	\label{e:GenMeanApproxB}
\end{equation}
\end{lemma}

\end{subequations}

\Proof
Recall that for a function $g$ with domain $\Re^d \times\state$ we say that   $\hag$ is a solution to Poisson's equation 
if $\Expect[\hag(\theta, \Phi_{n+1}) - \hag(\theta, \Phi_{n}) \mid \Phi_n = x] 
= 
- g(\theta, x) + \barg(\theta)$,   where $\barg(\theta) \eqdef \int g(\theta,x)  \uppi(dx)$.    In view of \Cref{t:bounds-Hgen} a solution exists,   satisfying    $ \int \hag(\theta,x)  \uppi(dx) = 0$,  and for every $\theta,\theta',x$ we have $| \hag(\theta, x) | \le   b_g  \bdd{t:bounds-Hgen}_p   [1+\|\theta\|]^3[1+V(x) ]^{p+1} $ along with the Lipschitz bound, 
\[
| \hag (\theta, x)  -  \hag (\theta', x) |   \le   b_g \bdd{t:bounds-Hgen}_p   [1+ \|\theta\| + \|\theta'\|]^2 [1+V(x)]^{p+1}  \|  \theta - \theta'\|.
\]

Poisson's equation gives $\Expect[\hag(\theta_n, \Phi_{n+k+1}) - \hag(\theta_n, \Phi_{n+k}) \mid \clF_{n+k}] 
= 
- g(\theta_n, \Phi_{n+k}) + \barg(\theta_n)$, and then by the smoothing property of conditional expectation and invariance,
\[
\Expect_\upvarpi \big[  g(\theta_0, \Phi_{k}) -  \barg(\theta_0)    \big]
=
\Expect_\upvarpi[ \hag(\theta_0, \Phi_{k}) - \hag(\theta_0, \Phi_{k+1})  ] 
= \Expect_\upvarpi[ \hag(\theta_{1}, \Phi_{k+1}) - \hag(\theta_0, \Phi_{k+1})  ] .
\]
Applying the Lipschitz bound,  
\[
\big |  \Expect_\upvarpi[  g(\theta_0, \Phi_{k}) ]  -  \Expect_\upvarpi[\barg(\theta_0)]  \big |  \le    
\alpha 
b_g \bdd{t:bounds-Hgen}_p   \Expect_\upvarpi  \bigl [ [1+\|\theta_0\|+\|\theta_{1}\|]^2 [1+V(\Phi_{k+1})]^{p+1}  \|  f(\theta_0, \Phi_{1} )\|  \big].
\]

\Cref{t:BigBounds}~(ii) implies that there is  a constant $b_\nu$  depending only on $p$ such that 
\[
\Expect_\upvarpi  \bigl [ [1+\|\theta_0\|+\|\theta_{1}\|]^2 [1+V(\Phi_{k+1})]^{p+1}  \|  f(\theta_0, \Phi_{1} )\|   \mid \clF_0 \big]  
\le    b_\nu  \bdd{t:BigBounds}   \clV(\Psi_{0}).
\]
Combining these bounds gives \eqref{e:GenMeanApprox}
with $ \bdd{t:PolyMeans}  = b_\nu  \bdd{t:BigBounds} \bdd{t:bounds-Hgen}_p  $.

For \eqref{e:GenMeanApproxB} we apply the fundamental theorem of calculus to obtain \eqref{e:FTC} with $\barf$ replaced by 
$\barg$,   and $A$ replaced by  $\derbarf_g$.           
Taking expectations  and applying the Lipschitz bound for $\derbarf_g$, we obtain via the Cauchy-Schwarz inequality
\[
\begin{aligned}
\big |    \Expect_\upvarpi[ \barg(\theta_0)]    -   \barg(\theta^*)    -  \derbarf_g (\theta^*)   \thbias  \big |  
& \le  
b_g  \Expect_\upvarpi \big[   [1+\|\theta_0 \|  +\|\theta^*\|  ]^2  \|  \theta_0 - \theta^*\|^2 \bigr]
\\
& \le  
b_g   \sqrt{    \Expect_\upvarpi \big[   [1+\|\theta_0 \|  +\|\theta^*\|  ]^4    \Expect_\upvarpi \big[   \|  \theta_0 - \theta^*\|^4 \bigr]    }.
\end{aligned}
\]
The right hand side is $O(\alpha)$ by  \Cref{t:BigBounds}  and  $\| \thbias \| =O(\alpha)$ by \Cref{t:AsyOops+Bias}, from which we obtain  \eqref{e:GenMeanApproxB}.
\qed

\begin{lemma}
\label[lemma]{t:SigmaDeltaDeltaStar}
Under the assumptions of  \Cref{t:sens},
the covariance matrices \eqref{e:SigmaW} satisfy 
$
\Sigma_{\MD}  =      \Sigma_{\MD^*}  +  O(\alpha)
$.
\end{lemma}

\Proof 
Fix indices $i,j$ and write  $g(\theta,x)  =  \Expect[ \MD_1^i \MD_1^j \mid \theta_0=\theta\,, \ \Phi_0=x]$.   
The function $g$ satisfies all of the assumptions  required in  \Cref{t:PolyMeans},   and by the definitions in \eqref{e:SigmaW} we have $
\Sigma_{\MD}^{i,j} = 
\Expect_\upvarpi[   g(\theta_0, \Phi_0) ]  $  and $   \Sigma_{\MD^*}^{i,j} =   \barg(\theta^*)  $.  
This combined with 
\Cref{t:PolyMeans} and  \Cref{t:BigBounds} implies the desired result:
$
\big|    \Sigma_{\MD}^{i,j}  -   \Sigma_{\MD^*}^{i,j}  \big|  \le    |  \Expect_\upvarpi[   \barg(\theta_0) -   \barg(\theta^*) ]  | + O(\alpha)  
= O(\alpha)  $.
\qed

\smallskip

\whamit{Proof of \Cref{t:sens}.}   
\Cref{t:PartialSensBdd} and \Cref{t:SensLLN}~(iii) establish \Cref{t:sens}~(iii).  
\Cref{t:BigBounds} gives
\eqref{e:theta4th-moment},   and   \eqref{e:theta2th-moment}  follows from
\Cref{t:AsyOops+Bias}. 
To complete the proof of (i)   we   establish the approximation  $\ZthetaAlpha= \Ztheta    +O(\alpha)$,  with 
$\ZthetaAlpha\eqdef \alpha^{-1}  
\Expect_{\upvarpi} [  \tiltheta_n \tiltheta^\transpose_n ] $ and   $\Ztheta\ge 0$ is the solution to  \eqref{e:Sigma0LyapEqn}.   For this write,
\[
\tiltheta_{n+1} =   [I +  \alpha \derbarf^*   ]\tiltheta_{n}+ \alpha \big[    \Delta_{n+1}    +  \clE^{\,\lilTS}_n
\big],
\]
where $ \clE^{\,\lilTS}_n = \barf(\theta_n) -    \derbarf^* \tiltheta_{n}$ satisfies, for some constant  $b^{\lilTS} $,   
\[
\| \clE^{\,\lilTS}_n \| \le  b^{\lilTS}  \| \tiltheta_{n} \|^2  \qquad
\| \clE^{\,\lilTS}_{n+1} -  \clE^{\,\lilTS}_n  \| \le \alpha \,    b^{\lilTS}  \| \tiltheta_{n} \|^2 .
\]
From the disturbance 
decomposition \eqref{e:DeltaDecomp} with $\alpha_n=\alpha$ we also have
\begin{equation}
\begin{aligned} 
\|   \Expect_{\upvarpi} [  \clE^{\,\lilTS}_n\Delta_{n+1}^\transpose ]   \|_F 
& = \|   \Expect_{\upvarpi} [  \clE^{\,\lilTS}_n  (- \clT_{n+1} + \clT_{n} - \alpha\Oops_{n+1})^\transpose ]   \|_F  
\\
& = \|   \Expect_{\upvarpi} [  (  \clE^{\,\lilTS}_{n+1} - \clE^{\,\lilTS}_n)  \clT_{n+1}^\transpose  
- \alpha \clE^{\,\lilTS}_n \Oops_{n+1})^\transpose ]   \|_F =  O(\alpha^2) ,
\\[.5em]
\||    \Expect_{\upvarpi}  [   \tiltheta_n  \Delta_{n+1}^\transpose    ] \|_F   
& = \|   \Expect_{\upvarpi} [   \tiltheta_n    (- \clT_{n+1} + \clT_{n} - \alpha\Oops_{n+1})^\transpose ]   \|_F  
\\
& = \|   \Expect_{\upvarpi} [  ( \tiltheta_{n+1} -  \tiltheta_n  ) \clT_{n+1}^\transpose  
- \alpha  \tiltheta_n  \Oops_{n+1})^\transpose ]   \|_F =  O(\alpha) .
\end{aligned}
\label{e:CovProof0}
\end{equation}

Consequently,  by invariance and the bound $  \Expect_{\upvarpi} [ \|  \tiltheta_n \|^2 ] = O(\alpha)$, 
\begin{equation}
\Expect_{\upvarpi} [  \tiltheta_{n+1} \tiltheta^\transpose_{n+1}] 
=
[I +  \alpha \derbarf^*   ]   \Expect_{\upvarpi} [  \tiltheta_n \tiltheta^\transpose_n ]   [I +  \alpha \derbarf^*   ]^\transpose  
+ \alpha   [Y + Y^\transpose]
+ \alpha^2    \Sigma_{\Delta}      +  O(\alpha^3),
\label{e:CovProof1}
\end{equation}
where $Y =    [I +  \alpha \derbarf^*   ]  \Expect_{\upvarpi}  [   \tiltheta_n  \Delta_{n+1}^\transpose    ]  
=    \Expect_{\upvarpi}  [   \tiltheta_n  \Delta_{n+1}^\transpose    ] +O(\alpha^2)$. 
The expectation may be expressed
\[
\begin{aligned}
\Expect_{\upvarpi}  [   \tiltheta_n  \Delta_{n+1}^\transpose    ]    
& 
=   \Expect_{\upvarpi}  [   \tiltheta_n  \Delta_{n+1}^{*\transpose }   ]  + O(\alpha^2)
\\  
& 
=   \Expect_{\upvarpi}  \big[   \tiltheta_n  \big[     \haDelta_{n+1}^{*\transpose }  -   \haDelta_{n+2}^{*\transpose }    \big]  \big]  + O(\alpha^2)
\\  
& 
=   \Expect_{\upvarpi}  \big[   \big[   \tiltheta_{n+1} - \tiltheta_n  \big]       \haDelta_{n+2}^{*\transpose }  \big]    + O(\alpha^2)  
\\  
& 
=   \alpha   \Expect_{\upvarpi}  \big[  \big(\barf(\theta_n) +  \Delta_{n+1}  \big)   \haDelta_{n+2}^{*\transpose }  \big]    + O(\alpha^2) ,
\end{aligned}
\]
where the first identity follows  \Cref{t:PolyMeans},   and the second is a consequence of Poisson's equation and the definition $  \haDelta^*_{n+1} = 
\haf(\theta^*, \Phi_{n+1})$:
\[
\Expect [   \haDelta^*_{n+1} \mid \clF_n ]  =  \haDelta^*_n -  \Delta_n^*.
\]
The third identity follows from invariance, similar to \eqref{e:CovProof0}:   $  \Expect_{\upvarpi} [  \tiltheta_n    \haDelta_{n+1}^{*\transpose } ] =   \Expect_{\upvarpi} [  \tiltheta_{n+1} 
\haDelta_{n+2}^{*\transpose } ] $.  Following the steps below \eqref{e:CovProof1}, applying  \Cref{t:PolyMeans} with $g(\theta_n,\Phi_{n+1}) = \barf(\theta_n) \haDelta_{n+2}^{*\transpose } $ and Poisson's equation,
gives 
\[
Y = 
\Expect_{\upvarpi}  [   \tiltheta_n  \Delta_{n+1}^\transpose    ]   + O(\alpha^2) 
=     \Expect_{\upvarpi}  \big[   \Delta^*_{n}     \haDelta_{n+1}^{*\transpose }  \big]    + O(\alpha^2)  
=     \Expect_{\upvarpi}  \big[   \Delta^*_{n}     \haDelta_{n}^{*\transpose }  \big]       -   \Expect_{\upvarpi}  \big[   \Delta^*_{n}     \Delta_{n}^{*\transpose }  \big]          + O(\alpha^2) . 
\]
This and \eqref{e:CovProof1}
implies the desired approximation $\ZthetaAlpha= \Ztheta +O(\alpha)$.


\Cref{t:AsyOops+Bias} combined with \Cref{t:PolyMeans} gives  \eqref{e:TargetBias} with
\[
\barOops_\alpha \eqdef       \Expect_\upvarpi  [ g(  \theta_n,\Phi_{n+1})  ]  =   \barg(\theta^*) +O(\alpha) \,,
\]
where $g(\theta_n,\Phi_{n+1} )  =  \Oops_{n+1}$. 
From the definition \eqref{e:Upsilon},    for any $\theta,x,\alpha$,
on denoting  $ \theta^+_{x,\alpha} = \theta + \alpha f(\theta,x)$,
\[
g(\theta,x) =   \frac{1}{\alpha}  \big[   \uppsi(   \theta^+_{x,\alpha}, x  )   -   \uppsi(\theta, x)     \big] .
\]
We have $ \partial_\theta   \uppsi = 
A(\theta,x) -\haderf(\theta,x)  $,  
where    $ A(\theta, x) =  \partial_\theta f(\theta,x) $,   $ \haderf(\theta, x) =  \partial_\theta \haf(\theta,x) $.    Consequently, 
under the assumptions of the theorem,     there is a fixed constant  $\bdd{t:AsyOops+Bias}$ such that
\[
g(\theta,x) =   [ \partial_\theta   \uppsi\, (\theta,x) ] f(\theta,x) +   \alpha B(\theta, x,\alpha)    V(x)   ,
\]
where  $  \|  B(\theta, x,\alpha) \|  \le  \bdd{t:AsyOops+Bias}   $  for  $\theta\in\Re^d$,   $x\in\state$,   and $\alpha\in (0,1]$.  
This establishes   \eqref{e:OopsApprox}.   
\qed

\subsection{Coupling and ergodicity for linear SA}
\label{s:linearSAapp}

The linear recursion \eqref{e:LinSAgen} may be expressed
\[
\theta_n = \Sens_{n} \theta_0 - \alpha \sum_{k=1}^n   \clL_{n,k} b_k  \,,  \ \ \textit{with  $\Sens_n =\clL_{n,0}$, }\quad n\ge 0\,,
\] 
in which the  random \textit{state transition matrices}  $\{ \clL_{n,k} :   k>n \ge 0 \}$ are defined by induction:   
\[
\clL_{ n, n} = I \,, \ \ \textit{and} \quad
\clL_{ n, k} =    [I+\alpha A_{k}]   \clL_{ n, k-1}\,, \quad k> n \ge 0\, .
\]
It follows that   $ \theta_n - \theta_n^\infty  =  \Sens_n [\theta_0 - \theta_0^\infty ]$, from which the following result is not surprising:

\begin{proposition}
\label[proposition]{t:SAlinErgodicity}
The following hold under the assumptions of \Cref{t:sensLinear}.
there is  $\bdd{t:SAlinErgodicity} <\infty$  such that
for each $0<\alpha  \le   \alpha_0 \eqdef \alpha_0 = \min(\alphaTmp{t:BigBounds},\alphaTmp{t:SensLLN})$:

\whamrm{(i)}
$ \Expect[ \|    \Sens_n\|_F^4 ]  \le \bdd{t:SAlinErgodicity}   \exp(-n \bdde{t:sensLinear} \alpha)   v(x) $
for each $\Psi_0 = z = (\theta;x)$ and 
$n\ge 0$,  
with $ \bdde{t:sensLinear} =2  \delta_s$.  

\whamrm{(ii)}
The uniform bound \eqref{e:L2LyapExp} holds  for each $\Psi_0$.

\whamrm{(iii)}  Suppose that  $G\colon \Re^d \times \state\to \Re$ is a measurable function satisfying  the following bounds: 
for some  $b_G<\infty$,   all   $  z= (x,u)\in \Re^d \times \state $, 
and    all $\theta,\theta' \in\Re^d$,
\[  
\begin{aligned}  &| G(\theta, x) |  \le   b_G  \clV^{1/4}(\theta, x) \, 
\\
&| G(\theta,x) - G(\theta', x)| 
\leq b_G  v^{1/4}(x)    \|\theta - \theta'\| .
\end{aligned}
\] 
Then,  $\Expect_\upvarpi[| G(\Psi_n) |^2]  \leq  \bdd{t:SAlinErgodicity}   b_G $, and
with  $\tilG = G - \upvarpi(G)$ we have each initial condition  $\Psi_0 = z $,
\begin{equation}
|\Expect[\tilG(\Psi_n) \mid \Psi_0=z]  |   \le b_G   \bdd{t:SAlinErgodicity}   \exp(-n \bdde{t:sensLinear} \alpha)      \clV^{1/2}(z)\,, \ \ n\ge 0.
\label{e:linSAergodic}
\end{equation}
\end{proposition}

\Proof    
\Cref{t:BigBounds}  applied to \eqref{s:SAlinearBivariate}  implies that  $ \Expect[ \|    \Sens_n\|_F^4 ]  \le \bdd{t:BigBounds} \clV(\Psi_0)  \exp( - 4 n\delta_s\alpha) $.   However, since $\{\Sens_n \}$ does not depend upon $\theta_0$,   we can set $\theta_0 \equiv 0$ to replace $ \clV(\Psi_0) $ by 
$v_+(\Phi_0)$.     Next,  by the Markov property,  followed by the smoothing property of conditional expectation,  
\[
\begin{aligned}
\Expect[ \|    \Sens_n\|_F^4  \mid   \Phi_1 ]  &\le   \bdd{t:BigBounds} v_+(\Phi_1)  \exp( -4 [n-1] \delta_s\alpha),
\\
\text{giving} \ \
\Expect[ \|    \Sens_n\|_F^4   ]  &\le   \bdd{t:BigBounds}  e^b  v(x)  \exp( -4 [n-1]\delta_s\alpha) \,, \ \ n\ge 1\, ,
\end{aligned}
\]
where $b>0$ is the constant in DV3.    This establishes (i).

By the Cauchy-Schwarz inequality,  for each  $\Psi_0 = z = (\theta;x)$,
\[
\Expect[ \|  \theta_n - \theta_n^\infty \|^2]  \le  \Expect[ \|    \Sens_n\|_F^2 \| \theta_0 - \theta_0^\infty\|^2 ]  
\le  \sqrt{ \Expect[ \|    \Sens_n\|_F^4 ]    \Expect[ \|     \theta_0 - \theta_0^\infty\|^4 ]  }.
\]
We have $   \Expect[ \|     \theta_0 - \theta_0^\infty\|^4 ]  \le b_1  (1 +  \beta_0 \|\theta\|^4)  $ for some $b_1$, where $\beta_0$ appears in the definition of $\clV$ in 	\eqref{e:clV}.   Combining this with part (i) gives  \eqref{e:L2LyapExp}, establishing (ii).

To establish the final result, note that
for any initial condition $\Psi_0=z$,
\[
\Expect[\tilG(\Psi_n) ]  =   \Expect[ G(\theta_n,\Phi_n) - G(\theta_n^\infty, \Phi_n^\infty)  ] 
\]
Under the assumptions of (iii) we have  $| G(\theta_n,\Phi_n) - G(\theta_n^\infty, \Phi_n^\infty)| 
\leq b_G  v^{1/4}(\Phi_n^\infty)    \|\theta_n - \theta_n^\infty\| $ for $n\ge   \taucpl$.     
Hence the proof of (iii) is obtained by combining \Cref{t:PhiCouples}   and \Cref{t:SAlinErgodicity}~(ii). 
\qed

\smallskip	

The remainder of the proof of \Cref{t:sensLinear} requires that we take a closer look at bias and variance.  
Approximations of the steady-state mean of a random vector $g( \Psi_k)$ are obtained via solutions to Poisson's equation.  
We require approximations only of functions that   are linear or quadratic as a function of $\theta_k$.

For a pair of measurable functions $L,M\colon\state\to\Re^{d\times d}$, 
denote $L_m =  L(\Phi_m) $,   $M_m =  M(\Phi_m) $ for $m\ge 0$, and 
consider the vector valued stochastic process $ \{  L_k \theta_k  :   k\ge 0 \} $  and the scalar valued stochastic process 
$ \{  \theta_k^\transpose  M_k \theta_k  :   k\ge 0  \} $.   
Under the assumptions imposed in \Cref{t:BiasLinear} there are solutions to Poisson's equation $\haL, \haM$,  normalized to have zero mean, and solve for $k\ge0$,
\begin{equation}
\Expect[    \haL (\Phi_{k+1})    \mid \clF_k] =    \haL (\Phi_{k})    -  [  L_k   - \barL ] \,,  \ \ 
\Expect[    \haM (\Phi_{k+1})    \mid \clF_k] =    \haM (\Phi_{k})    -  [  M_k   - \barM ] \,,  
\label{e:PoissonForBiasFormulae}
\end{equation}
where $\barL$,  $\barM$ denote respective steady-state means.
We write  $\haL_k =  \haL(\Phi_{k}) $,   $\haM_k =   \haM(\Phi_{k}) $.

\begin{proposition}
\label[proposition]{t:BiasLinear}  
Suppose that the assumptions of \Cref{t:sensLinear} hold, and  suppose the measurable functions $L$ and $M$ satisfy the bounds, 
\begin{equation}
\| L( x)  \|^2_F  +   \| M( x)   \|^2_F  \leq  b_0    [1+V(x) ]^p \,, \quad x\in\state \,,
\label{e:LMbdds}
\end{equation}
for some $p\ge 1$ and  $  b_0 <\infty$.
Then, there is a constant $\bdd{t:noise-decomp-gen} <\infty$ such that  for each $\alpha \in (0, \alpha_0]$,  
\whamrm{(i)}
Zero-mean solutions to \eqref{e:PoissonForBiasFormulae} exist satisfying  for a constant $b_p$, 
\[
\| \haL( x)  \|^2_F  +   \| \haM( x)   \|^2_F  \leq   b_p  b_0    [1+V(x) ]^{p+1} \,, \quad x\in\state.
\]

\whamrm{(ii)}      $\Expect_\upvarpi[ L_k \theta_k]  = \barL \theta^*   + \alpha   [  \Uppi^{L*}  +  \barL \Zthbias ] + O(\alpha^2) $,   
where   $ \Uppi^{L*}  =   \Expect_\upvarpi[ \haL_{k+1} f_{k+1}(\theta_k) ]$,   with $f_{k+1}(\theta_k)  = A_{k+1}\theta_k - b_{k+1}$.

The special case $L_k\equiv I$ gives  
$ \thbias \eqdef \Expect_\upvarpi[ \tiltheta_k ]   =  \alpha \Zthbias   + O(\alpha^2)$,  
with  $\Zthbias  \eqdef 
[\derbarf^*]^{-1} \barOops^*$.

\whamrm{(iii)}  $\Expect_\upvarpi[ \theta_k^\transpose M_k  \theta_k]  =  {\theta^*}^\transpose \barM\theta^* + 
\alpha \big[   \Uppi^{M*}     +  \trace (\barM \Ztheta) \big]
+ O(\alpha^2)$,      with $\Ztheta\ge 0$ the solution to \eqref{e:Sigma0LyapEqn},   and 
$  \Uppi^{M*}  =   \Expect_\upvarpi \big[ { \theta^*}^\transpose  (  \haM_{k+1} +   \haM_{k+1} ^\transpose)  f_{k+1}(\theta^*)      \big] $.
\end{proposition}

\Proof  Part (i) is immediate from \Cref{t:bounds-Hgen}.
For (ii) apply Poisson's equation and invariance to obtain 
\[
\begin{aligned}
\Expect_\upvarpi[ L_k \theta_k]  &=
\Expect_\upvarpi[ (   \haL_k - \haL_{k+1}  +\barL)  \theta_k]  
\\
&=     \barL \Expect_\upvarpi[   \theta_k]  +
\Expect_\upvarpi[  \haL_{k+1}  (  \theta_{k+1} -\theta_k  ) ]           
\\
&=     \barL \theta^*  +   \barL\thbias  +  \alpha \Uppi^L \,, 
\end{aligned}
\]
where  $ \Uppi^L  =    \Expect_\upvarpi[ \haL_{k+1} f_{k+1}(\theta_k) ]$.    In particular,   $ \|  \Expect_\upvarpi[ L_k \theta_k]    -   \barL \theta^* \| =   O(\alpha)$  since by \Cref{t:sens} we have  $\|\thbias\| = O(\alpha)$.

For the linear SA recursion we may write $ \Expect[ \haL_{k+1} f_{k+1}(\theta_k)  \mid \clF_k]  =   L_k^1 \theta_k + g_k$ in which $\{ L_k^1, g_k\}$ satisfy bounds of the form \eqref{e:LMbdds} for a possibly larger value of $p$.    Consequently,
\[
\Uppi^L =   \Expect_\upvarpi[  L_k^1 \theta_k + g_k ]   =  (  \barL^1 \theta^* + \barg ) +O(\alpha)  =  \Uppi^{L*} +O(\alpha) .
\]
Combining these  identities gives 
\begin{equation}
\Expect_\upvarpi[ L_k \theta_k]  = \barL \theta^*  +   \barL\thbias +   \alpha  \Uppi^{L*} +O(\alpha).
\label{e:linBias1}
\end{equation}
It remains to justify the approximation for $\thbias$. 

Taking expectations of each side of \eqref{e:LinSAgen} and applying invariance gives
$ 0= \Expect_\upvarpi[ A_{n+1} \theta_n - b_{n+1}]   = \Expect_\upvarpi [A_{n+1} \tiltheta_n ],$
where the second inequality follows because     $0=\barf(\theta^*) = \Expect_\upvarpi [A_{n+1} \theta^* - b_{n+1}] $.   
Recalling the definition of $\haA$ in \eqref{e:haA}, we repeat the arguments above:
\[
\begin{aligned}
0 = 
\Expect_\upvarpi [A_{n+1} \tiltheta_n ]   & =  \Expect_\upvarpi [  ( \haA_{n+1}- \haA_{n+2} +\derbarf)  \tiltheta_n ] 
\\
&= \derbarf \thbias  +    \Expect_\upvarpi [   \haA_{n+2}  ( \tiltheta_{n+1} - \tiltheta_n  )].
\\
\end{aligned}
\]
From the foregoing we have 
\[
\barOops^* \eqdef  - \Expect_\upvarpi [   \haA_{n+2}  (A_{n+1} \theta^* - b_{n+1}   )] 
=  
- \Expect_\upvarpi [   \haA_{n+2}  (A_{n+1} \theta_n - b_{n+1}   )]      + O(\alpha)  .
\]
This completes the proof that $ \thbias   =  \alpha \Zthbias   + O(\alpha^2)$, which when combined with  \eqref{e:linBias1} completes the proof of (ii).

Part (iii) is similar:  recalling that   $f_{k+1}(\theta_k)  = A_{k+1}\theta_k - b_{k+1}$,
\[
\begin{aligned}
\Expect_\upvarpi[ \theta_k^\transpose M_k  \theta_k]   & =  
\Expect_\upvarpi[ \theta_k^\transpose  (  \haM_k - \haM_{k+1}  +\barM )  \theta_k]     
\\
&=  \Expect_\upvarpi[ \theta_k^\transpose  \barM   \theta_k]   
+\Expect_\upvarpi[ \theta_{k+1}^\transpose   \haM_{k+1}    \theta_{k+1}  ]     
- \Expect_\upvarpi[ \theta_k^\transpose   \haM_{k+1}     \theta_k]     
\\
&=  \Expect_\upvarpi[ \theta_k^\transpose  \barM   \theta_k]   
+  \alpha  \Expect_\upvarpi \big[  \theta_{k}^\transpose   \haM_{k+1}  f_{k+1}(\theta_k)     +   f_{k+1}(\theta_k) ^\transpose   \haM_{k+1} \theta_{k}   \big]
+ O(\alpha^2)
\\
&=   {\theta^*}^\transpose \barM\theta^* +  
\alpha \Big(  \trace (\barM \Ztheta)  +  \Expect_\upvarpi \big[  \theta_{k}^\transpose  ( \haM_{k+1} + \haM_{k+1} ^\transpose) f_{k+1}(\theta_k)       \big]  \Big),
+ O(\alpha^2)
\end{aligned} 
\]
where the final approximation follows from   \eqref{e:theta2th-moment}.   Thus  we have shown that  $ \Expect_\upvarpi[ \theta_k^\transpose M_k  \theta_k]   =  {\theta^*}^\transpose \barM\theta^* + O(\alpha)$ under these conditions on $\{ M_k \}$.    As in the proof of (ii), we may write,
\[
\Expect_\upvarpi[ \theta_{k}^\transpose   \haM_{k+1}  f_{k+1}(\theta_k)  \mid \clF_k]      = \theta_k^\transpose M_k^1  \theta_k  +  d_k^\transpose      \theta_k ,
\]
in which $\{ M_k^1, d_k\}$ satisfy bounds of the form \eqref{e:LMbdds} for a possibly larger value of $p$.    
Consequently,
\[
\begin{aligned}
\Expect_\upvarpi \big[  \theta_{k}^\transpose  ( \haM_{k+1} + \haM_{k+1} ^\transpose) f_{k+1}(\theta_k)       \big]   
&=   \Expect_\upvarpi \big[  \theta_k^\transpose ( M_k^1 + {M_k^1}^\transpose)  \theta_k  +  2 d_k^\transpose      \theta_k \big]
\\
&=   \Expect_\upvarpi \big[ {\theta^*}^\transpose  ( M_k^1 + {M_k^1}^\transpose)  \theta_k  +  2 d_k^\transpose \theta^* \big]   +  O(\alpha)
\\
&=   \Uppi^{M*}  +  O(\alpha), 
\end{aligned}
\]
which completes the proof of (iii).
\qed

\smallskip 

We next consider bounds on the asymptotic covariance.  
A warning is required here:   the matrix $\SigmaPR$ as defined in  \eqref{e:SigmaPR} is not finite if $\thbias\neq 0$ in \eqref{e:Bias},  since for each $N$,
\[
\Expect[ (\thetaPR_N -\theta^*)(\thetaPR_N -\theta^*)^\transpose]  =  \Cov\,(\thetaPR_N )  +    \Expect[ \thetaPR_N ]  \Expect[  \thetaPR_N ] ^\transpose \, .
\]
In view of this, we consider instead the limit of the scaled covariance,  as in \Cref{e:PRCovrep}  below.

Recall the CLT for $\bfPhi$ in \Cref{t:Vuni}~(ii).  Our interest here is asymptotic statistics for functionals of the joint process $\bfPsi$, which is more challenging because it is not necessarily  ergodic in the sense of \cite[Ch.~17]{MT}, because $\psi$-irreducibility may fail.

\begin{subequations}%

Much of the work here focuses on the ``disturbance process'' $\bfDelta$, whose asymptotic covariance is defined similarly:
\begin{equation}
\SigmaCLT^\Delta = \lim_{N\to\infty}
\frac{1}{N}    \Cov\Big(\sum_{i=1}^N \Delta_i    \Big).
\label{e:SigmaDelta}
\end{equation}	
One step in the lemma that follows is obtained by averaging each side  
of \eqref{e:noisyEuler}, with $\alpha_n=\alpha$,
using $\barf(\theta) = \derbarf^* [\theta - \theta^*]$:
\begin{equation}
\thetaPR_N=   \theta^*  - [\derbarf^*]^{-1} \frac{1}{N-N_0}     [ 	S^{\Delta}_N  -  S^\tau_N],
\label{e:PRrep}
\end{equation}
where $S^\tau_N  = 
\frac{1}{\alpha} (\theta_{N} - \theta_{N_0+1}) 
$ and $
S^{\Delta}_N 
=
\sum_{k=N_0+1}^N \Delta_{n} $.

\end{subequations}%

\begin{lemma}
\label[lemma]{t:covtilPR}
Under the assumptions of \Cref{t:sensLinear},  the limit \eqref{e:SigmaDelta} exists,  is finite, and independent of the initial condition $\Psi_0$.   Moreover,  for each initial condition,
\begin{equation}
\lim_{N \to \infty } N\Cov(\thetaPR_N )   
=  
G^* \SigmaCLT^\Delta {G^* }^\transpose.
\label{e:PRCovrep}
\end{equation}
\end{lemma}

\Proof
The proof of the existence of the limit  \eqref{e:SigmaDelta} begins with the representation
\[
\frac{1}{N}    \Cov\Big(\sum_{i=1}^N \Delta_i   \Big)  	
=
\frac{1}{N}  \sum_{i,j=1}^N  R^\Delta_{i,j},
\]
where 
$R^\Delta_{i,j}  =  \Expect[\tilde  \Delta_i {\tilde  \Delta_j}^\transpose]$,  
using $\tilde  \Delta_i  = \Delta_i - \Expect[   \Delta_i ]$.  
Let  $R^{\Delta^\infty}_{i,j}   $ denote the auto-covariance sequence for the stationary process $\{\Delta_{k+1}^\infty = f(\theta_{k}^\infty,\Phi_{k+1}^\infty) - \barf(\theta_{k}^\infty) \}$.    Applying \eqref{e:linSAergodic},   we can find $b_d<\infty $  such that  for all $i,j$,
\[
\begin{aligned}
\| R^\Delta_{i,j} \|_F  & \le  b_d \clV(\Psi_0)   \exp(- |i-j| \bdde{t:sensLinear} \alpha)  
\\
\| R^\Delta_{i,j} - R^{\Delta^\infty}_{i-j}  \|_F  &  \le    b_d \clV(\Psi_0) 
\exp(- \min(i,j) \bdde{t:sensLinear} \alpha).
\end{aligned}  
\]
This establishes both finiteness and independence of the initial condition in  \eqref{e:SigmaDelta},
and the standard alternative representation
$\SigmaCLT^\Delta  = \sum_{n=-\infty}^\infty R^{\Delta^\infty}_n$.

Taking covariances of both sides of \eqref{e:PRrep} and rearranging terms, we obtain
\begin{equation*}
(N - N_0)\Cov(\thetaPR_N )   
=  
\frac{1}{N-N_0} [\derbarf^*]^{-1}  \Cov(S^{\Delta}_N  - S^\tau_N) [{\derbarf^*}^\transpose]^{-1}.
\end{equation*}
The proof is completed on observing that
\[
\lim_{N\to\infty}	\frac{1}{N-N_0} 
\Cov(S^{\Delta}_N) =     \SigmaCLT^\Delta  \,,
\qquad 
\lim_{N\to\infty}	\frac{1}{N-N_0} 
\Cov( S^\tau_N) =  0
\]
\qed

\smallskip

The remainder of this subsection is devoted to approximation of $\SigmaCLT^\Delta$.    In particular, while \Cref{t:covtilPR} asserts that this matrix is finite, an application of  \eqref{e:linSAergodic}
would give the unacceptable bound $\| \SigmaCLT^\Delta \|_F \le O(1/\alpha)$.    \Cref{t:SigmaDeltaDeltaStar} allows to establish a bound that remains bounded for arbitrarily small $\alpha>0$.

As in  \Cref{t:BiasLinear}  we present approximations for a general vector valued stochastic process.   We require approximations of   the asymptotic covariances of several processes, such as   $ \Oops_{k+1}     =   \uppsiA_{k+1}  [  A_{k+1} \theta_{k}  -  b_{k+1} ]$  (recall   \eqref{e:LinDistDecomTerms}).

We consider a vector valued function similar to what is considered in \Cref{t:BiasLinear}:  for   a given pair of measurable functions $L\colon\state\to\Re^{d\times d}$,   $g\colon\state\to\Re^{d }$ we define 
\begin{equation}
\clG_{k+1}  =  L_{k+1} \theta_k + g_{k+1} \,, \quad k\ge0\,,
\label{e:clGvariance}
\end{equation}
where  $L_k =  L(\Phi_{k}) $,   $g_k =   g(\Phi_{k}) $.   
Let $\barL$,  $\barg$ denote the respective steady state means of $L$, $g$,     denote $\barG(\theta) = \barL \theta +\barg$,
and $\barG_k = \barG(\theta_k)$.

It will be assumed that $L,g$ satisfy bounds of the form \eqref{e:LMbdds}, so that we are assured of   solutions to Poisson's equation:
\[
\Expect[    \haL (\Phi_{k+1})    \mid \clF_k] =    \haL (\Phi_{k})    -  [  L_k   - \barL ] \,,  \ \ 
\Expect[    \hag (\Phi_{k+1})    \mid \clF_k] =    \hag (\Phi_{k})    -  [  g_k   - \barg ] \,.
\]
As always, we assume the solutions are normalized to have zero steady-state mean.
Applying \Cref{t:bounds-H}, as in \Cref{t:noise-decomp},  we obtain the decomposition,
\begin{equation}
\clG_{k+1}  =  \barG_k  + \MD^G_{k+1} -  \clT^G_{k+1} +  \clT^G_{k} - \alpha  \Oops^G_{k+1}   
\label{e:Gdecomp}
\end{equation}
We are interested in approximations of the steady-state mean and covariance of $\clG_{n+1} $,  and bounds on asymptotic covariances   
for this sequence and also   $\{ \MD^G_{n+1} \}$, denoted,
\[
\SigmaCLT^{\clG}  = \sum_{n=-\infty}^\infty   \Expect_\upvarpi[ \tilclG_0  \tilclG_n^\transpose]   \,, 
\quad
\SigmaCLT^{\MD^G}  =    \Expect_\upvarpi[\MD^G_0 (\MD^G_n )^\transpose] .
\]
The asymptotic covariance $\SigmaCLT^{\MD^G} $ is simplified due to the martingale difference property.

\begin{subequations}

\begin{proposition}
\label[proposition]{t:noise-decomp-gen}
Suppose that the assumptions of \Cref{t:sensLinear} hold, and consider the stochastic process \eqref{e:clGvariance} subject to the bounds 
\begin{equation*}
	\| L( x)  \|^2_F  +   \| g( x)   \|^2  \leq  b_\clG    [1+V(x) ]^p \,, \quad x\in\state \,,
\end{equation*}
for some $p\ge 1$ and  $  b_\clG <\infty$.
Then, there is a constant $\bdd{t:noise-decomp-gen} <\infty$ such that  for each $\alpha \in (0, \alpha_0]$:
\whamrm{(i)}
The representation \eqref{e:Gdecomp}
holds:  denoting $\haclG_{k+1} =  \haL_{k+1} \theta_k   +  \hag_{k+1}$,
\begin{align}
	\MD^G_{n+1} &  \eqdef \haclG_{n+1}  - \Expect[ \haclG_{n+1}   \mid \clF_n],    
	\label{e:MDG}
	\\[0.5em]
	\clT^G_{n+1} &\eqdef    [   L_{n+1}  -\haL_{n+1}] \theta_{n}   +   g_{n+1}  -\hag_{n+1},
	\label{e:telescopeG}
	\\
	\Oops_{n+1}^G &\eqdef   [   L_{n+1}  -\haL_{n+1}]   [   L_{n+1} \theta_n + g_{n+1} ].
	\label{e:UpsilonG}
\end{align}
The final term may be expressed $ \Oops_{n+1}^G  =  H^G(\Phi_{n+1})  \theta_k + \gamma^G (\Phi_{n+1}) $  in which  
\[
\| H^G( x)  \|^2_F  +   \| \gamma^G( x)   \|^2  \leq   \bdd{t:noise-decomp-gen}  b_\clG    [1+V(x) ]^q\, \ \ \textit{$x\in\state$,  for some $q\ge p$.}
\]

\whamrm{(ii)}   	The steady state mean admits the approximation  
\begin{equation}
	\Expect_\varpi[\clG_n ] =   \barG(\theta^*)     + \alpha\big[  \barL \Zthbias -  (  \barH^G   \theta^* + \bargamma^G ) \big],
	+ O(\alpha^2)     
	\label{e:biasLinearGen}
\end{equation}
where $  \barH^G \eqdef \Expect_{\uppi}[  H^G(\Phi_{n+1}) ]   $,      $ \bargamma^G   \eqdef  \Expect_{\uppi}[   \gamma^G (\Phi_{n+1}) ] $. 

\whamrm{(iii)}      The asymptotic covariances satisfy $ \| \SigmaCLT^{\clG}  \|_F^2     +  \|   \SigmaCLT^{\MD^G} \|_F^2  \le \bdd{t:noise-decomp-gen}   b_\clG$.

\whamrm{(iv)}  
$ \| \SigmaCLT^{\Delta}  \|_F    =    \|    \SigmaPR \|_F   \le \bdd{t:noise-decomp-gen}  $,  and the asymptotic covariance for $\{\barG_n \}$  is
\begin{equation}
	\SigmaCLT^{\MD^{\barG}}  =  \barL  \SigmaPR \barL^\transpose=  \barL  G^* \SigmaCLT^\Delta   (   \barL  G^*  )^\transpose   
	\label{e:barGacov}
\end{equation}
\end{proposition}

\label{e:noise-decomp-gen}
\end{subequations}

\Proof
The proof of (i) is identical to the proof of \Cref{t:noise-decomp}, with bounds obtained using \Cref{t:bounds-Hgen},
and (ii) is immediate from (i).

For (iii), first note that  the existence of $\bdd{t:noise-decomp-gen} $ such that the bound  $  \|   \SigmaCLT^{\MD^G} \|_F^2   \le \bdd{t:noise-decomp-gen}   b_\clG$  holds is immediate from (i) and the martingale difference property (so that $\SigmaCLT^{\MD^G}  = \Cov(\MD_n^G)$, the steady-state covariance).

A similar bound for $ \| \SigmaCLT^{\clG}  \|_F $ is obtained by following the proof of \Cref{t:covtilPR}.
Letting $b_p$ denote a constant for which $ [1+V(x) ]^p \le b_p   v^{1/4}(x)$ for all $x$,  we obtain from   \eqref{e:linSAergodic} of
\Cref{t:SAlinErgodicity},  
\[ 
\| \SigmaCLT^{\clG} \|_F   = \sum_{n=-\infty}^\infty   \| \Expect_\upvarpi[ \tilclG_0  \tilclG_n^\transpose] \|_F   \le  
2  b_p  \sqrt{b_\clG} \bdd{t:SAlinErgodicity}  \sum_{n=0}^\infty   \Expect_\upvarpi[   \| \tilclG_0  \|     \clV^{1/2}(\Psi_0) ]  \exp(-n \bdde{t:sensLinear} \alpha)     .
\]
Under the given assumptions we have 
\[
\Expect_\upvarpi[   \| \tilclG_0  \|     \clV^{1/2}(\Psi_0) ]   \le    \sqrt{  
\Expect_\upvarpi[   \| \clG_0  \|^2 ]     \Expect_\upvarpi[   \clV(\Psi_0) ]   }  
\le  
b_p
\sqrt{  b_\clG  }  \Expect_\upvarpi[   \clV(\Psi_0) ] ,
\] 
giving
$\| \SigmaCLT^{\clG} \|_F   \le 2 b_p b_\clG    \Expect_\upvarpi[   \clV(\Psi_0) ] \bdd{t:SAlinErgodicity} / [1- \exp(-\bdde{t:sensLinear} \alpha) ]$,  which is finite and no larger than  $O(1/\alpha)$.

To obtain a better bound, first  consider the special case $\clG_n = \Delta_n$, for which $\barG_n=0$.     The CLT covariance of the telescoping sequence $\{ \clT_{n+1} -  \clT_{n}\}$ is zero.  Hence,
with $\alpha_n=\alpha$, the  representation \eqref{e:DeltaDecomp} gives 
$
\| \SigmaCLT^\Delta \|_F \le 
\| \SigmaCLT^\MD \|_F  +  \alpha
\| \SigmaCLT^\Oops \|_F 
$.
\Cref{t:SigmaDeltaDeltaStar}  implies that  $\| \SigmaCLT^\MD \|_F $ is uniformly bounded in $\alpha$.  
Also, we can apply the foregoing bounds to obtain $\| \SigmaCLT^\Oops \|_F \le  
b_\Oops b_\clG /\alpha$ with $b_\Oops$ independent of $\alpha$.  This establishes (iv).

Using (iv) we can complete the proof of (iii):
exactly as in the preceding bound on $\| \SigmaCLT^\Delta \|_F $,  
\eqref{e:Gdecomp} gives,
\[
\| \SigmaCLT^\clG \|_F \le 
\| \SigmaCLT^{\MD^{\barG}} \|_F  + 
\| \SigmaCLT^{\MD^G} \|_F  +  \alpha \| \SigmaCLT^{\Oops^G} \|_F.
\]
The final two terms 
admit uniform bounds as in the preceding arguments.
The representation for $ \SigmaCLT^{\MD^{\barG}}  $ in 
\eqref{e:barGacov} follows from   \Cref{t:covtilPR},
whose norm is bounded by a constant times $b_\clG$ on applying (iv). 
\qed

\smallskip

\whamit{Proof of \Cref{t:noise-decomp-linear}.}  Part (i) is immediate from \Cref{t:noise-decomp}.

For (ii)  first observe that (i) implies  $ \Oops_{n+1}  =  H_{n+1} \theta_n + \gamma_{n+1}$,
in which     $ H_{k+1} =  \uppsiA_{k+1}   A_{k+1} $,   $\gamma_{k+1} = -  \uppsiA_{k+1}   
b_{k+1} $, and Poisson's equation implies that  
$\barH \eqdef    \Expect_\uppi  [  \uppsiA_{k+1}   A_{k+1}  ] =   - \Expect_\uppi  [\haA_{k+2}   A_{k+1} ]$.
An application of   \Cref{t:noise-decomp-gen} implies 
the representation \eqref{e:OopsDecomp} holds in which
$\barH = \Expect_\uppi[ \uppsiA_{n+1}   A_{n+1}  ]$,   $\bargamma = -  \Expect_\uppi[ \uppsiA_{n+1}   b_{n+1} ]$,   
$\MD^H_{n+1}   \eqdef  \haOops_{n+1}  - \Expect[  \haOops_{n+1}   \mid \clF_n]   $,
\begin{align*}
\clT^H_{n+1}   &  \eqdef    [   H_{n+1}  -\haH_{n+1}] \theta_{n}   +   \gamma_{n+1}  - \hagamma_{n+1}  
\ \ \textit{and} \ \ 
\Oops_{n+1}^H   \eqdef   [   H_{n+1}  -\haH_{n+1}]   [  A_{n+1} \theta_{n}  -  b_{n+1} ].
\end{align*}

For (iii) express  the $i,j$ entry of the covariance matrix as $ \Sigma_{\MD}^{i,j} =  \Expect[ \theta_k^\transpose M_k  \theta_k+   L_k \theta_k   ] $ in which   $\theta_k^\transpose M_k  \theta_k  +   L_k \theta_k  =  \Expect [ \MD_{k+1}^i\MD_{k+1}^j \mid \clF_k]$.      
We have $L_m =  L^{i,j}(\Phi_m) $,   $M_m =  M^{i,j}(\Phi_m) $ for $m\ge 0$, in which the  matrix valued functions $L^{i,j}$,  $M^{i,j}$
satisfy the assumptions of \Cref{t:BiasLinear}.   Consequently, on applying  \Cref{t:BiasLinear},  
\[
\Sigma_{\MD}^{i,j} = 
\Sigma_{\MD^*}^{i,j} +   \alpha   [   \Uppi^{L*}   +\barL \Zthbias  ] + 
\alpha \big[   \Uppi^{M*}     +  \trace (\barM \Ztheta) \big]
+ O(\alpha^2) .
\]
The $i,j$ entry of   $\ZMD  $  is obtained from the definitions. 
\qed

\smallskip

We are now in a position to prove  \Cref{t:sensLinear}.

\whamit{Proof of \Cref{t:sensLinear}.}   Part (i)  is a consequence of \Cref{t:SAlinErgodicity},  and (ii)  follows from \Cref{t:noise-decomp-gen}~(ii).

The challenge is (iii), which requires the refinement of the noise decomposition \eqref{e:DeltaDecomp} obtained in  \Cref{t:noise-decomp-linear}. 
Applying  \Cref{t:noise-decomp-linear}~(ii),  we revisit \eqref{e:noisyEuler} 
with $\alpha_n=\alpha$ to obtain
\[
\theta_{n+1} = \theta_{n}+ \alpha \big(  
\derbarf^* [\theta_n - \theta^*] +     \MD^\sbullet_{n+1} - \clT^\sbullet_{n+1} + \clT^\sbullet_{n} - \alpha [\barH \theta_n   + \bargamma  ]   +\alpha^2 \Oops_{n+1}^H  
\big),
\]
where $\MD^\sbullet_{n+1} =\MD_{n+1} - \alpha \MD^H_{n+1}  $ and 
$\clT^\sbullet_{n+1} =\clT_{n+1} - \alpha \clT^H_{n+1}  $.
Returning to \eqref{e:PRrep},   
\[
\begin{aligned}
\derbarf^* [	\thetaPR_N -  & \theta^* ]   =     - \frac{1}{N-  N_0}  \Big[S^\tau_N  +   \sum_{k=N_0+1}^N \Delta_k  \Big]
\\
&=     -  \frac{1}{N-  N_0}     \sum_{k=N_0+1}^N  \{  \MD^\sbullet_k   -   \alpha  	   [ \barH \theta_k  + \bargamma   ]  +\alpha^2 \Oops_{k}^H    \}   +   \frac{1}{N-  N_0}     [S^\tau_N+\clT^\sbullet_{N+1}  -   \clT^\sbullet_{N_0+1}   ].
\end{aligned} 
\]
On rearranging terms,
\[
[\derbarf^* -  \alpha \barH ] [	\thetaPR_N -   \theta^* ]  =     -  \frac{1}{N-  N_0}       \sum_{k=1}^N \Delta_k 
=     -  \frac{1}{N-  N_0}     \sum_{k=1}^N  \{  \MD_k^\sbullet   -   \alpha  	   [      \barH  \theta^*+ \bargamma  +\clE_k ]   \},
\]
where $\{ \clE_n \}$ has asymptotic variance of order  $O(\alpha^2)$;  in particular,   \Cref{t:noise-decomp-gen} implies that the 
asymptotic variance of $\{
\Oops_{k}^H \}$ is bounded, uniformly in $\alpha\in (0,\alpha_0]$.    Pre-multiplying by $G_\alpha = -  [\derbarf^* -  \alpha \barH ]^{-1} $   gives
\[
(N - N_0)\Cov(\thetaPR_N )    = G_\alpha \Big[    \frac{1}{N-  N_0}   \sum_{k=1}^N   \Cov( \MD_k^\sbullet  )   \Big] G_\alpha ^\transpose  +    D_N,
\]
where   $\| D_N \|_F \le O(\alpha^2)$.   
Substitution of the approximation   $G_\alpha =   G^*  -  \alpha G^*  \barH G^* + O(\alpha^2)$ then gives
\begin{equation}
\begin{aligned} 
\lim_{N\to\infty}
N &\Cov(\thetaPR_N )  =     \SigmaTheta^*   +   \alpha   Z_\alpha  
+ O(\alpha^2)   ,
\\
&\text{where} \ \   Z_\alpha \eqdef    \tfrac{1}{\alpha} \big(  
[ G^*  -  \alpha G^*  \barH G^* ]  \Sigma_{\MD^\sbullet} [ G^*  -  \alpha G^*  \barH G^* ]  ^\transpose       -  \SigmaTheta^* \big),
\end{aligned}
\label{e:AlmostZ}
\end{equation}
and $ \Sigma_{\MD^\sbullet}   \eqdef  \Cov(    \MD_k^\sbullet )$  is obtained for the stationary version of $\bfPsi$.

To complete the proof of (iii) requires the approximation $Z_\alpha = \ZthetaPR +O(\alpha)$,  and that we identify  the matrix $ \ZthetaPR$.      First, from the definitions,
\[
\Sigma_{\MD^\sbullet}  =       \Sigma_{\MD}   
- \alpha      \big[  \Sigma_\otimes + \Sigma_\otimes^\transpose   \big]  + O(\alpha^{2}) \,, \quad \textit{where $\Sigma_\otimes = \Expect[  \MD_k {\MD^H_k}^\transpose ]$.}
\]
Consequently, writing $\SigmaTheta^{\alpha *} =   G^*     \Sigma_{\MD}  [ G^* ]^\transpose  $  and recalling that $[ G^* ]^{-1} = -\derbarf^*$,   
\[
\begin{aligned}
\alpha  Z_\alpha &=   \SigmaTheta^{\alpha *}   -  \SigmaTheta^*        
+   \alpha G^* \bigl[    \Sigma_\otimes + \Sigma_\otimes^\transpose    \bigr]  [ G^* ]^\transpose 
-\alpha  \bigl[  G^*   \barH \SigmaTheta^{\alpha *}    -  \SigmaTheta^{\alpha *} [ G^* \barH ]^\transpose  \bigr]   +O(\alpha^2)   
\\
&=   \SigmaTheta^{\alpha *}   -  \SigmaTheta^*   + \alpha G^* \bigl[   \Sigma_\otimes^* + {\Sigma_\otimes^*}^\transpose    \bigr]  [ G^* ]^\transpose 
+ \alpha G^* \bigl[  \derbarf^* \SigmaTheta^* \barH^\transpose  + \barH   \SigmaTheta^*[\derbarf^*]^\transpose      \bigr]  [ G^* ]^\transpose   +O(\alpha^2)    
\\
&=    \alpha G^* \bigl[  \ZMD  +  \Sigma_\otimes^* + {\Sigma_\otimes^*}^\transpose    \bigr]  [ G^* ]^\transpose 
+ \alpha G^* \bigl[  \derbarf^* \SigmaTheta^* \barH^\transpose  + \barH   \SigmaTheta^*[\derbarf^*]^\transpose      \bigr]  [ G^* ]^\transpose   +O(\alpha^2) ,
\end{aligned}
\]
where in the second equation  we introduced the notation $\Sigma_\otimes^* =  \Expect[   \MD^{*}_{n+1}  (  \MD^{H*}_{n+1}  )^\transpose ] $,
and used $\SigmaTheta^{\alpha *}   = \SigmaTheta^ *  +O(\alpha)$.   In the final expression we 
applied \Cref{t:noise-decomp-linear}. 
\qed

\subsection{Theory for SGD}
\label{s:optAssumptionsProof}

\Cref{t:optAssumptions} is largely a corollary to \Cref{t:BigBounds}:

\whamit{Proof of \Cref{t:optAssumptions}.}  
To establish (i) we first note that for the function
$V_\circ(\theta,x) = \half [\alpha^{-1} \nu(\theta) + V(x)]$,
regardless of the definition of $\nu$,
\begin{equation}
\begin{aligned}
\Expect\big[ \exp\big(  &V_\circ(\Psi_{n+1})    \mid \Psi_n=(\theta;x)  \big]   
\\
&  \le   
\sqrt{     \Expect\big[ \exp\big(  \alpha^{-1} \nu(\theta_{n+1})     \big) \mid \Psi_n=(\theta;z)  \big]      \Expect\big[ \exp\big(  V(\Phi_{n+1})   \big) \mid \Phi_n=x \big]   }
\\
&  \le   
\sqrt{     \Expect\big[ \exp\big(  \alpha^{-1} \nu(\theta_{n+1})     \big) \mid \Psi_n=(\theta;z)  \big]    }   
\exp\big( \half  [V(x)  -   W(x) +  b   ]   \big),
\end{aligned}
\label{e:optAss1}
\end{equation}
where the second inequality follows from \eqref{e:DV3}.   

We next apply a drift condition for the mean flow, which will inform the definition of $\nu$.  
Prop.~4.34 of \cite{CSRL} implies that there is a Lipschitz continuous function $\bigupgamma_0\colon\Re^d \to\Re_+$ and  constant $ B_0 < \infty$  such that  
$\nabla \bigupgamma_0 \, (\theta)  \cdot \barf (\theta)  \le  -  \|\theta\| + B_0$ for all $\theta$;    
It can be assumed smooth by convolution with a Gaussian kernel,  and a glance at the proof shows that it has linear growth in $\| \theta \|$.   Consequently,  for sufficiently large $b_0 \ge 1$, and a possibly larger constant $B_0$,
the function $\bigupgamma\eqdef  b_0 \bigupgamma_0^2$ satisfies   $\nabla \bigupgamma \, (\theta)  \cdot \barf (\theta) \le  -  \|\theta\|^2 + B_0$.     We define $\nu = \delta \bigupgamma$ for a constant $\delta>0$ that is specified in the following.

Repeating the arguments leading to \eqref{e:SensDelta-theta} we obtain,
\[
\bigupgamma(\theta_{n+1} )  = 
\bigupgamma(\theta_{n} )  +  A_{n+1}^\upgamma [  \theta_{n+1}  -\theta_{n} ]\,, \ \   \clA_{n+1}^\upgamma \eqdef   \int_0^1  \partial_\theta \bigupgamma(\theta_{n+1}^t) \, dt ,
\]
where $\theta_{n+1}^t = (1-t)\theta_n + t \theta_{n+1}$ for each $t$.    Letting $\ell_\upgamma$ denote a Lipschitz constant for $\nabla \bigupgamma$ and
noting that $\|\theta_{n+1}^t - \theta_n \|  = t  \|\theta_{n+1}  -\theta_{n} \|$ we obtain for any  $\beta>0$,
\[
\begin{aligned}
\bigupgamma(\theta_{n+1} ) & \le    \bigupgamma(\theta_{n} )  +    \alpha \nabla \bigupgamma\, (\theta_n) \big[   \barf(\theta_n) +  \Probe_{n+1}    \big]   +   \half  \alpha^2  \ell_\upgamma  \big\| \barf(\theta_n) +  \Probe_{n+1}    \big\|^2 
\\
&  \le     \bigupgamma(\theta_{n} )   -  \alpha   \| \theta_n\|^2  +\alpha   \big[  B_0+ \tfrac{1}{\beta}   \| \nabla \bigupgamma\, (\theta_n) \|^2 +  {\beta}\| \Probe_{n+1} \|^2 \big]
+ \alpha^2 \ell_\upgamma  \big[  \| \barf(\theta_n) \|^2 +  \|  \Probe_{n+1}    \|^2\big] . 
\end{aligned}
\]
We have for some $b_2>0$ and any $\theta$ the bounds $  \| \nabla \bigupgamma\, (\theta) \|^2 +  \| \barf(\theta) \|^2  \le b_2 [1 +  \| \theta\|^2 ]$,  giving  
\[
\tfrac{1}{\alpha} \bigupgamma(\theta_{n+1} )  \le \tfrac{1}{\alpha}     \bigupgamma(\theta_{n} )   -   [1 - b_2(  \tfrac{1}{\beta}  +  \alpha \ell_\upgamma)   ] [1+  \| \theta_n\|^2   ] +1
+ B_0 +  [\alpha\ell_\upgamma + \beta]  \| \Probe_{n+1} \|^2.
\]
Assuming that $\alpha_0 \ell_\upgamma <1/2$ we may set $\beta = 2 b_2/(1 - 2\alpha_0\ell_\upgamma)$ to obtain $1 - b_2 [ \tfrac{1}{\beta}  +  \alpha \ell_\upgamma ]  \ge  \half$ for $0<\alpha\le \alpha_0$, and hence for this range of $\alpha$ and  any $\delta>0$,
\[
\Expect\big[ \exp\big(  
\tfrac{\delta}{\alpha} \bigupgamma(\theta_{n+1} )   \big) \mid \clF_n \big]     
\le   \exp\big( \tfrac{\delta}{\alpha}     \bigupgamma(\theta_{n} )   -    \half  \delta    \| \theta_n\|^2   
+ \delta [1+  B_0]     \big)       \Expect\big[ \exp\big(   \delta[ \alpha_0\ell_\upgamma +\beta ]    \| \Probe_{n+1} \|^2 \big)   \mid \clF_n \big] .
\]
Applying \eqref{e:ProbeAssumptions} we have, whenever $ \delta[ \alpha_0\ell_\upgamma +\beta ]   \le \bddepsy{t:optAssumptions}$,
\[
\Expect\big[ \exp\big( \delta [ \alpha_0\ell_\upgamma +\beta ]    \| \Probe_{n+1} \|^2 \big)   \mid \Phi_n=x  \big]    \le     \exp( \delta [ \alpha_0\ell_\upgamma +\beta ] /\bddepsy{t:optAssumptions}]   [ W( x)        +   \bdd{t:optAssumptions} ] ).
\]
Choose  $\delta>0$ so that  $\delta [ \alpha_0\ell_\upgamma +\beta ] /\bddepsy{t:optAssumptions}  \le 1/2$.  
Recalling that we define   $\nu = \delta \bigupgamma$,
\[
\Expect\big[ \exp\big(  \tfrac{1}{\alpha} \nu(\theta_{n+1})     \big) \mid \Psi_n=(\theta;z)  \big]  
\le   \exp\big(  \tfrac{1}{\alpha} \nu (\theta  )   -    \half  \delta   \| \theta\|^2  
+ \delta [1+  B_0]      + \half  [\bdd{t:optAssumptions}  +       W( x)  ]        \big)  .  
\]
Combining this bound with \eqref{e:optAss1}  gives, for a constant $B_1$,
\[
\Expect\big[ \exp\big(  V_\circ(\Psi_{n+1})    \mid \Psi_n=(\theta;x)  \big]     \le
\exp\big(   V_\circ(\Psi_{n})     -    [ 1+     (\delta  \| \theta\|^2   +  W(x) )/4 ]   + B_1 \big).
\]
This establishes (i), where the terms are specified as follows:
$b_\circ = 2 B_1$,
and since  $W_\circ(\theta,x) =  1+ \epsy_\circ [  \| \theta \|^2 + W(x)]$,   
we fix
$\epsy_\circ   <  \max(1, \delta)/4$ and define 
$ C_{\tState} = \{ (\theta;x) :  (\delta  \| \theta\|^2   +  W(x) )/4     \le   \epsy_\circ (  \| \theta \|^2 + W(x)  )   + B_1\}$.  
By construction we have 
$C_\circ  \subset C_\tTheta \times C_{\tState}$, in which 
$ C_\tTheta= \{ \theta  :   ( \delta/4 -\epsy_\circ )    \| \theta \|^2   \le  B_1 \}$ and 
$ C_{\tState} = \{ x  :  (1/4 - \epsy_\circ)   W(x)   \le  B_1 \}$.

Part (ii) follows from 
\Cref{t:BigBounds}~(i). 

\qed

\smallskip

\whamit{Proof of \Cref{t:TowOptCases}.}
We present the proof of (ii) only since the proof of (i) is similar.  
The existence of the density \eqref{e:ProbeAssumptionsDensity} holds for $n_0=2$. 
For  $x = (w;w') \in\Re^d\times\Re^d$,  the notation $ \Phi_n= x$ means that $W_{n-1}=w$ and $W_n=w'$.   Consequently, 
\[
\Expect[ \exp\big(\epsy \| \Probe_{1} \|^2 \big ) \mid  \Phi_0= x]   \le 
\Expect[ \exp\big( 2\epsy [  \| W_{1} \|^2+\| W_{0} \|^2 ]  \big ) \mid  W_0 = w']    \le \exp\big(  \Lambda (2\epsy)   +  2\epsy \|w'\|^2\big).
\]
We obtain \eqref{e:ProbeAssumptions}  on taking 
$\bddepsy{t:optAssumptions} = \epsy$ with  $\epsy =  \min(1/2, \epsy_0/8)$, and  $\bdd{t:optAssumptions} = \Lambda (2\epsy)$.

Similarly,  with $V(w,w') = \epsy_0[ \half \| w \|^2 +\| w' \|^2 ] $,   
\[
\Expect[ \exp\big( V(\Phi_{n+1} )  \big)\mid  \Phi_n= x]   =     \exp\big(   \half \epsy_0 \|w'\|^2  + \Lambda(\epsy_0)  \big  )   =    \exp\big( V(x)- \half \epsy_0 [    \|w\|^2 + \|w'\|^2  ] + \Lambda(\epsy_0)   \big).
\]
This establishes \eqref{e:DV3}  and the remaining parts of (A2ii) easily follow.
\qed

\subsection{Numerical experiments: impact of memory}
\label{s:NumProofs}

\whamit{Proof of \Cref{t:LinExample}.}	
Applying \eqref{e:Bias} and substituting  $\derbarf^* = -1$ gives the bias formula \eqref{e:BiasLinExample}.    

It remains to establish the formula for $\Expect_\upvarpi  [ \Oops_{0}  ] $ in  \eqref{e:UpsSSlinExample}.
It is convenient to consider the stationary realization
\[
\MD_{n}  =   \sqrt{1+\beta^2}   \sum_{k=0}^\infty \beta^k  \clW^\circ_{n+1 - k},
\]
which is also an $N(0,1)$ random variable for each $n$.   A solution to Poisson's equation is easily obtained:
\[
\widehat\MD_{n}\eqdef  \frac{1}{1-\beta}  \MD_{n}.
\]
This solves $\Expect[\widehat\MD_{n+1} \mid \clF_n] = \widehat\MD_n  - \MD_n$ for each $n$.  More important for understanding $\Oops$ is that $\haA_{n+1} \eqdef \widehat\MD_{n}$ solves Poisson's equation \eqref{e:haA} with $\derbarf^*=-1$.
The representation \eqref{e:Upsilon}  becomes
\[
\begin{aligned} 
\Oops_{n+1} 
&=    -  \frac{1}{1-\beta}  \MD_{n+1}  \bigl[   ( \MD_{n} -1) \theta_n + \theta^* +     \MD_{n} \bigr]
\\
&=    -  \frac{1}{1-\beta}  \MD_{n+1} \MD_{n} \bigl[ 1+ \theta^*    \bigr] 
-  \frac{1}{1-\beta}  \MD_{n+1}      ( \MD_{n} -1) [\theta_n-\theta^*  ]  .
\end{aligned} 
\]
Using $\Expect[\MD_{n+1}    \MD_{n} ] = \beta$ then gives the steady-state formula  \eqref{e:UpsSSlinExample}.

We next establish (ii).
The formula \eqref{e:SigmaPRopt} must be modified when the noise is not white:  $ \Sigma_{\MD^*}$ is replaced by $\Sigma_{\Delta^*}$,  using the steady-state disturbance
\[
\Delta_{n+1}^* \eqdef  f(\theta^*, \Phi_{n+1}) - \barf(\theta^*)  =   f(\theta^*, \Phi_{n+1}) .
\]
In the present example $\Delta_{n+1}^* = (1+\theta^*) \MD_{n}$. Recalling  $\derbarf^*=-1$ then gives $\SigmaTheta^*  = \Sigma_{\Delta^*} $.    

However,   $\Sigma_{\Delta^*}  =    (1+\theta^*)^2\Sigma_\clW $,   with 
\[
\Sigma_\clW 
= 
\sum_{k=-\infty}^\infty   \Expect_\upvarpi[   \MD_0 \MD_k ]  
=
\Expect_\upvarpi[  2 \widehat\MD_0 \MD_0  -  \MD_0^2].
\]
Putting these formulae together gives (ii):
$\displaystyle
\SigmaTheta^*  = \Sigma_{\Delta^*}  
=  
\Expect_\upvarpi[  2 \widehat\MD_0 \MD_0  -  \MD_0^2]   (1+\theta^*)^2 
=  
\frac{1+\beta}{1-\beta}   (1+\theta^*)^2  
$.
\qed

\end{document}